\documentclass[12pt]{article}

\title{From Subfactors to Categories and Topology II. \\
       {\large The Quantum Double of Tensor Categories and Subfactors}}
\author{Michael M\"uger\thanks{Supported by EU through the TMR Networks ``Noncommutative
Geometry'' and  ``Algebraic Lie Representations'', by MSRI through NSF grant
DMS-9701755 and by NWO.} \\ Korteweg-de Vries Institute, Amsterdam, Netherlands \\ email: {\tt mmueger@science.uva.nl}}

\newlength{\dinwidth}
\newlength{\dinmargin}
\usepackage{latexsym, amssymb, amsmath, theorem}

\usepackage{tan-gle}

\input diagrams         

\def\1#1{{\bf #1}}
\def\2#1{{\cal #1}}
\def\3#1{{\sl #1}}
\def\4#1{{\tt #1}}
\def\5#1{{\sf #1}}
\def\6#1{{\mathfrak #1}}
\def\7#1{{\mathbb #1}}

\newcommand{\be}{\begin{equation}}
\newcommand{\ee}{\end{equation}}
\newcommand{\ba}{\begin{array}}
\newcommand{\ea}{\end{array}}
\newcommand{\bea}{\begin{eqnarray}}
\newcommand{\eea}{\end{eqnarray}}
\newcommand{\bean}{\begin{eqnarray*}}
\newcommand{\eean}{\end{eqnarray*}}
\newcommand{\nn}{\nonumber}

\newcommand{\ve}{\varepsilon}

\newcommand{\impl}{\Rightarrow}
\newcommand{\rarr}{\rightarrow}
\newcommand{\restr}{\upharpoonright}
\newcommand{\ol}{\overline}
\newcommand{\ul}{\underline}
\newcommand{\obj}{\mbox{Obj}}
\newcommand{\Hom}{\mbox{Hom}}
\newcommand{\HOM}{{\cal HOM}}
\newcommand{\END}{{\cal END}}
\newcommand{\End}{{\mbox{End}}}
\newcommand{\id}{\mbox{id}}
\renewcommand{\mod}{\mbox{mod}}
\newcommand{\Rep}{\mbox{Rep}}

\newcommand{\qed}{\ \hfill $\blacksquare$\medskip}

\newcommand{\cmp}{Commun. Math. Phys. }

\newcommand{\rmp}{Rev. Math. Phys. }

\def\endexem{\hfill{$\Box$}\medskip}

\theoremheaderfont{\scshape}
\newtheorem{defin}{Definition}[section]
\newtheorem{lemma}[defin]{Lemma}
\newtheorem{prop}[defin]{Proposition}
\newtheorem{theorem}[defin]{Theorem}
\newtheorem{coro}[defin]{Corollary}
\newtheorem{conj}[defin]{Conjecture}
\theorembodyfont{\rmfamily}
\newtheorem{rema}[defin]{Remark}

\newcommand{\bdefin}{\begin{defin}}
\newcommand{\blemma}{\begin{lemma}}
\newcommand{\bprop}{\begin{prop}}
\newcommand{\btheor}{\begin{theorem}}
\newcommand{\bcoro}{\begin{coro}}
\newcommand{\bconj}{\begin{conj}}
\newcommand{\edefin}{\end{defin}}
\newcommand{\elemma}{\end{lemma}}
\newcommand{\eprop}{\end{prop}}
\newcommand{\etheor}{\end{theorem}}
\newcommand{\ecoro}{\end{coro}}
\newcommand{\econj}{\end{conj}}
\newcommand{\brem}{\begin{rema}}
\newcommand{\erem}{\endexem\end{rema}}

\newcommand{\oj}{\overline{J}}
\newcommand{\op}{{\mbox{\scriptsize op}}}
\newcommand{\rev}{{\mbox{\scriptsize rev}}}

\newcommand{\ha}{Hopf algebra}

\newcommand{\qtha}{quasitriangular Hopf algebra}
\newcommand{\qd}{quantum double}
\newcommand{\prf}{{\it Proof. }}
\def\mychi{\raise 2pt\hbox{$\chi$}}

\newcommand{\mcirc}{\,\circ\,}

\def\mobj#1{\raise .4\unitlens\hbox{\put(0,0){$#1$}}}
\def\mmobj#1{\raise .3\unitlens\hbox{\put(0,0){$#1$}}}
\def\mmmobj#1{\raise .5\unitlens\hbox{\put(0,0){$#1$}}}

\begin{document}
\maketitle\noindent

\numberwithin{equation}{section}

\abstract{For every tensor category $\2C$ there is a braided tensor category $\2Z(\2C)$,
the `center' of $\2C$. It is well known to be related to Drinfel'd's notion of the quantum
double of a finite dimensional Hopf algebra $H$ by an equivalence 
$\2Z(H-\mod)\stackrel{\otimes}{\simeq}_{br} D(H)-\mod$ of braided tensor categories. In
the Hopf algebra situation, whenever $D(H)-\mod$ is semisimple (which is the case iff
$D(H)$ is semisimple iff $H$ is semisimple and cosemisimple iff $S^2=\id$ and 
$\mbox{char}\7F\nmid\dim H$) it is modular in the sense of Turaev, i.e.\ its $S$-matrix is
invertible. (This was proven by Etingof and Gelaki in characteristic zero. We give a
fairly general proof in the appendix.) The present paper is concerned with a
generalization of this and other results to the quantum double (center) of more general
tensor categories.

We consider $\7F$-linear tensor categories $\2C$ with simple unit and finitely many
isomorphism classes of simple objects. We assume that $\2C$ is either a $*$-category
(i.e.\ $\7F=\7C$ and there is a positive $*$-operation on the morphisms) or semisimple and
spherical over an algebraically closed field $\7F$. In the latter case we assume
$\dim\2C\equiv\sum_i d(X_i)^2\ne 0$, where the summation runs over the isomorphism
classes of simple objects. We prove that $\2Z(\2C)$ (i) is a semisimple spherical (or
$*$-) category and (ii) is weakly monoidally Morita equivalent (in the sense of
math.CT/0111204) to $\2C\otimes_\7F\2C^\op$. This implies $\dim\2Z(\2C)=(\dim\2C)^2$. (iii) 
We analyze the simple objects of $\2Z(\2C)$ in terms of certain finite dimensional algebras,
of which Ocneanu's tube algebra is the smallest. We prove the conjecture of Gelfand and
Kazhdan according to which the number of simple objects of $\2Z(\2C)$ coincides with the
dimension of the state space $\2H_{S^1\times S^1}$ of the torus in the triangulation TQFT
built from $\2C$. (iv) We prove that $\2Z(\2C)$ is modular and we compute 
$\Delta_\pm(\2Z(\2C))\equiv\sum_i \theta(X_i)^{\pm 1}d(X_i)^2=\dim\2C$.
(v) Finally, if $\2C$ is already modular then 
$\2Z(\2C)\stackrel{\otimes}{\simeq}_{br}\2C\otimes_\7F\tilde{\2C}\stackrel{\otimes}{\simeq}\2C\otimes_\7F\2C^\op$, 
where $\tilde{\2C}$ is the tensor category $\2C$ with the braiding
$\tilde{c}_{X,Y}=c_{Y,X}^{-1}$. 
}

\section{Introduction}
Define the `center' $Z_0(X)$ of a set $X$ to be the monoid of all functions 
$f\colon X\rightarrow X$ (with composition as product and the identity map as unit.)
Then the usual center $Z_1\equiv Z$ of the monoid $Z_0(X)$
is trivial: $Z_1(Z_0(X))=\{\id_X\}$. The cardinality of $Z_0(X)$ is 
given by $\#Z_0(X)=\#X^{\#X}$. The aim of the present work is to prove 1-categorical
analogues of these trivial set theoretic (= 0-categorical) observations.
(I owe the above definition of $Z_0(X)$ to J.\ Baez.)

Given an arbitrary monoidal category (or tensor category) $\2C$ its center $\2Z(\2C)$
is a braided monoidal category which was defined independently by Drinfel'd (unpublished),
Majid \cite{ma1} and Joyal and Street \cite{js2}. (See Section \ref{sec-semisim} for the
definition.) In order to avoid confusion with another notion of center, we will write
$\2Z_1(\2C)$ throughout. In the present work, as in \cite{ma1, js2}, we will assume $\2C$
to be strict, but this is exclusively for notational convenience. The definition of the
center $\2Z_1$ and all results in this paper extend immediately to the non-strict
case. The other assumptions which we must make on $\2C$ are more restrictive, but we are
still left with a class of categories which appears in contexts like low dimensional
topology and subfactor theory. We assume $\2C$ to be linear over a ground field which is
algebraically closed. Furthermore, $\2C$ is semisimple with simple tensor unit and
spherical \cite{bw2}. (A semisimple category is spherical iff it is pivotal 
\cite{bw2} (=sovereign) and every simple object has the same dimension as its dual, cf.\
\cite[Lemma 2.8]{mue09}.) See \cite{bw2} or \cite[Section 2]{mue09} for the precise
definitions.

\bdefin \label{catdim}
Let $\2C$ be a semisimple spherical tensor category with simple unit and let $\Gamma$ be
the set of isomorphism classes of simple objects. If $\Gamma$ is  finite we define 
\[ \dim\2C=\sum_{i\in\Gamma} d(X_i)^2, \]
otherwise we write $\dim\2C=\infty$. 

If $\2C$ is finite dimensional and braided then the Gauss sums of $\2C$ are given by
\[ \Delta_\pm(\2C)=\sum_{i\in\Gamma} \omega(X_i)^{\pm 1}d(X_i)^2, \]
where $\theta(X)=\omega(X)\id_X$ is the twist of the simple object $X$ which is defined by
the spherical structure \cite{y}. 
\edefin

We can now state our Main Theorem:
\btheor \label{main1}
Let $\7F$ be an algebraically closed field and $\2C$ a spherical $\7F$-linear tensor
category with $\End(\11)\cong\7F$. We assume that $\2C$ is semisimple with finitely many
simple objects and $\dim\2C\ne 0$. Then also the center $\2Z_1(\2C)$ has all these
properties and is a modular category \cite{t}. Furthermore, the dimension and the Gauss
sums are given by
\bean \dim\2Z_1(\2C) &=& (\dim\2C)^2, \\
   \Delta_+ (\2Z_1(\2C)) &=& \Delta_-(\2Z_1(\2C)) = \dim\2C. \eean
\etheor
Defining the center $\2Z_2(\2C)$ of a {\it braided} tensor category $\2C$ to be the 
full subcategory whose objects are those $X$ which satisfy
\[ c(X,Y)=c(Y,X)^{-1} \quad\forall Y\in\obj\,\2C, \]
one easily sees that $\2Z_2(\2C)$ is stable w.r.t.\ isomorphisms (thus replete), direct
sums, retractions, tensor products and duals, the inherited braiding obviously being
symmetric. One can show that a braided category satisfying the properties in the theorem 
(i.e.\ a premodular category \cite{brug1}) is modular iff the center $\2Z_2(\2C)$ is
trivial, in the sense that all objects of $\2Z_2(\2C)$ are multiples of the tensor
unit. (This was done in \cite{khr1} for $*$-categories and in \cite{bebl} for spherical 
categories with $\dim\2C\ne0$, see also Corollary \ref{rem-mod} below.) Thus
$\2Z_2(\2Z_1(\2C))$ is trivial for all $\2C$ as in the Main Theorem, which is the promised
analogue of the 0-categorical observation $Z_1(Z_0(X))=\{\11\}$. 

The Main Theorem can be generalized slightly: If $\2C$ is as before except for $\7F$ not
being algebraically closed then there is a finite extension $\7F'\supset\7F$ 
such that $\2Z_1(\2C\otimes_\7F \7F')$ is modular. Concerning the prospects of further  
generalizations the author is not optimistic. There is little hope of proving
semisimplicity of $\2Z_1(\2C)$ without assuming $\dim\2C\ne 0$. (Furthermore, it is known
\cite{t} that the dimension of a modular category must be non-zero.) In the non-semisimple
case one might hope to prove that the center of a spherical noetherian category satisfies
the non-degeneracy condition on the the braiding introduced in \cite{ly2}. But the methods
of this paper will most likely not apply. 

The results of the present work can be considered as generalizations of known results
concerning Hopf algebras and we briefly comment on this in order to put our results into
their context. We recall that the quantum double of a Hopf algebra was introduced, among
many other things, in Drinfel'd's seminal work \cite{drin}. In the following discussion
all Hopf algebras are finite dimensional over some field $\7F$. The quantum double $D(H)$
of a Hopf algebra $H$ is a certain Hopf algebra which contains $H$ and the dual $\hat{H}$
as Hopf subalgebras and it is generated as an algebra by these. We refrain from repeating
the well known definition and refer to \cite{ka} for a nice  treatment. We only remark
that $D(H)\cong H\otimes_\7F\hat{H}$ as a vector space, thus  
\[ \dim_\7F D(H)=(\dim_\7F H)^2. \]
Furthermore, $D(H)$ is quasitriangular, i.e.\ there is an invertible 
$R\in D(H)\otimes D(H)$ such that $\sigma\circ\Delta=R \Delta(\cdot)R^{-1}$ where
$\sigma$ is the flip automorphism of the tensor product. The constructions of the quantum
double of a Hopf algebra and of the center of a monoidal category are linked by the
equivalence 
\[ D(H)\mbox{-mod} \stackrel{\otimes}{\simeq}_{br} \2Z_1(H\mbox{-mod}) \]
of braided monoidal categories, where $H\mbox{-mod}$ and $D(H)\mbox{-mod}$ are the 
categories of finite dimensional left $H$- and $D(H)$-modules, respectively, the 
braiding of $D(H)\mbox{-mod}$ being provided by the R-matrix. Again, see 
\cite[Chapter XIII.4]{ka} for a detailed account. 

Now, the R-matrix of a quantum double $D(H)$ is non-degenerate in a certain sense, 
$D(H)$ being `factorizable' \cite{rst}. If $H$ is semisimple and cosemisimple then $D(H)$ 
is semisimple \cite{rad}. It then turns out to be also modular and the category 
$D(H)\mbox{-mod}$ of finite dimensional left $D(H)$-modules is modular in the sense of 
Turaev \cite{t}. (This was proved in \cite{eg} for algebras over algebraically closed
fields of characteristic zero, but the latter condition can be dropped. In the appendix we
give a general proof.) Furthermore, one clearly has 
\begin{equation}\dim \2Z_1(H\mbox{-mod}) = \dim D(H)\mbox{-mod} = \dim D(H)=(\dim H)^2 =
   (\dim H\mbox{-mod})^2, \label{dims}\end{equation}
where $\dim\2C$ is the dimension of the monoidal category $\2C$ as defined above.

It is now clear that our Main Theorem can be considered as an extension of the above
results to tensor categories which are not necessarily representation categories of Hopf
algebras. Here one remark on the notation is in order. In \cite{kt} Kassel and Turaev
introduced a modified version of the construction of the center $\2Z_1(\2C)$ and called it
the quantum double $\2D(\2C)$, see also \cite{s}. Their category is the categorical
version of a construction of Reshetikhin (which adjoins a certain square root $\theta$ to
a quasitriangular Hopf algebra $H$ in order to turn it into a ribbon algebra $H(\theta)$)
applied to a quantum double, cf.\ \cite[Theorem 5.4.1]{kt}. In the context of \cite{kt}
the starting point was that even if $\2C$ is rigid this need not be true for $\2Z_1(\2C)$,
whereas the category $\2D(\2C)$ is rigid. As we will see, spherical categories (tensor
categories with nice two sided duals) are better behaved in the sense that their centers
$\2Z_1$ are again spherical. In addition, whereas $\2Z_1(\2C)$ is modular for the
categories satisfying the conditions of our Main Theorem, this is never true for
$\2D(\2C)$! This is why we stick to the original definition $\2Z_1(\2C)$. Apart from
writing $\2Z_1(\2C)$ instead of $\2Z(\2C)$, we do not attempt to change the established
symbols, but we use the expression `quantum double' as a synonym for $\2Z_1(\2C)$ rather
than $\2D(\2C)$. 

Unfortunately, the work on Hopf algebras mentioned above provides no clues on how to prove 
Theorem \ref{main1}. This is where subfactor theory enters the present story. Starting from
an inclusion $N\subset M$ of hyperfinite type $II_1$ factors of finite index and depth,
Ocneanu \cite{ocn2} defined an `asymptotic subfactor' $B\subset A$: 
\[ B=M \vee (M_\infty\cap M')\subset M_\infty=A. \]
(Here $N\subset M\subset M_1\subset M_2\subset\ldots$ is the Jones tower associated with
$N\subset M$ and $M_\infty=\vee_i M_i$.) In \cite{ocn3} he argued that a certain monoidal
category associated with $B\subset A$ is braided, concluding that the asymptotic subfactor
is an `analogue' of Drinfel'd's quantum double of a Hopf algebra. In fact, Ocneanu does
not use category language and does not refer to the quantum double (center) of monoidal
categories. In \cite{ek1} Evans and Kawahigashi published proofs for most of Ocneanu's
announcements. In the paper \cite{lre}, which otherwise has little to do with the
asymptotic subfactor, Longo and Rehren then constructed a subfactor $B\subset A$ 
from an infinite factor $M$ and a -- in our language -- finite dimensional full monoidal 
subcategory $\2C$ of $\End(M)$ and conjectured that it is related to Ocneanu's
construction. This conjecture was made precise and proven in \cite{mas}.
The author's involvement in the present story began when in 1998 he received
a copy of a short preprint \cite{iz1} by M. Izumi. In the meantime a full account 
of Izumi's results has appeared in \cite{iz2}. In \cite{iz1,iz2} Izumi gives an in-depth
analysis of the LR-subfactor, in particular its $B-B$ sectors. Seeing \cite{iz1} the
present author was struck by the fact that its main theorem implicitly contained the
definition of the center of a monoidal category. In fact properly formulated, Izumi's
results provide a precise and completely general form of Ocneanu's `analogy' between
the asymptotic subfactor (or the LR-subfactor) and the quantum double, albeit the 
categorical one instead of the one for Hopf algebras. In Section \ref{ss-subf} we will
rephrase Izumi's results in categorical language to make this evident. Yet, this is not
the main purpose of the present work.

In \cite{ocn3} it has been argued that the braided monoidal category associated with 
$B\subset A$ is modular and a complete proof has been provided in \cite{iz2}, where 
it was also shown that the dimension of the category in question is given by $(\dim\2C)^2$.
As in our discussion of the Hopf algebra quantum double, it is again natural to ask
whether a purely categorical version of these results can be proven. Here we have to face
the problem that finite-index subfactors have a lot of `in-built' categorical structure
which is not a priori available in a purely categorical setting. (In particular, most of
\cite{iz2} strongly relies on this structure.) Yet this problem can be overcome once one
realizes that the more algebraic part of subfactor theory can be cast into the language of
2-categories. This is the content of \cite{mue09}, which in a sense can be considered a
continuation of \cite{lro}, though in a somewhat more general setting.

The paper is organized as follows. In Section 2 we first recall some of the less standard 
definitions from \cite{mue09}. We then summarize the main results of \cite{mue09} on
Frobenius algebras in tensor categories, related 2-categories and the notion of weak
monoidal Morita equivalence of tensor categories. This section can by no means replace
\cite{mue09}. Our study of the quantum double $\2Z_1(\2C)$ begins in Section
\ref{sec-semisim}, where we show that it preserves the closedness w.r.t.\ direct sums and
subobjects and sphericity. Most importantly and least trivially, we prove the
semisimplicity of $\2Z_1(\2C)$. These results do not yet rely on the machinery of
\cite{mue09}. In Section \ref{sec-morita} we prove the weak monoidal Morita equivalence  
$\2Z_1(\2C)\approx\2C\otimes_\7F\2C^\op$, which in particular implies that the double
construction squares the dimension of the category. Section \ref{sec-modular} is devoted
to the proof of modularity of $\2Z_1(\2C)$, equivalent to triviality of the category
$\2Z_2(\2Z_1(\2C))$. As an important first step we analyze the structure of the simple
objects of the double, providing an explanation for Ocneanu's `tube algebra'. 
The next two sections consider the case of categories with a positive $*$-structure
($C^*$-categories or unitary categories) and the special case where $\2C$ is already
braided. In Section \ref{sec-appl} we consider applications to the invariants of
3-manifolds, proving a conjecture of Gelfand and Kazhdan and speculating about a far
stronger result. Finally, we establish the link with subfactor theory, relying heavily on
Izumi's work, improving on it only slightly.


\section{Preliminaries}
\subsection{Some definitions and notations}
We refer to \cite[Section 2]{mue09} for our general conventions and recall only a few less
standard notations. A retract $Y\prec X$ is also called a subobject. A category has
subobjects if all idempotents split, and every category $\2A$ has a canonical completion 
$\ol{\2A}^p$ for which this is the case. A $\7F$-linear category is semisimple if it has
direct sums and subobjects and every object is a finite direct sum of simple objects, $X$
being simple iff $\End(X)\cong\7F$. For monoidal categories we require in addition that
$\11$ is simple. A subcategory of a semisimple category is called semisimple if it is
closed w.r.t.\ direct sums and retractions, thus in particular replete (stable under
isomorphisms).

Since all categories in question are $\7F$-linear we understand the product 
$\2K\otimes_\7F\2L$ of (tensor) categories in the sense of enriched category theory. Thus 
\[ \obj\,\2K\otimes_\7F\2L = \{ K\boxtimes L,\ K\in \obj\,\2K, L\in \obj\,\2L\}, \]
where $X\boxtimes Y$ stands for the pair $(X,Y)$, and
\[ \Hom_{\2K\otimes_\7F\2L}(K_1\boxtimes L_1, K_2\boxtimes L_2))= 
   \Hom_\2K(K_1,K_2)\otimes_\7F\Hom_\2L(L_1,L_2) \]
with the obvious composition laws. 
We denote by $\2K\boxtimes\2L=\ol{\2K\otimes_\7F\2L}^\oplus$ the completion w.r.t.\
finite direct sums. If $\2X, \2Y$ are monoidal categories the same holds for
$\2X\boxtimes\2Y$. In order to save brackets we declare $\boxtimes$ to bind stronger than 
$\otimes$ but weaker than juxtaposition $XY$ of objects (which abbreviates $X\otimes
Y$). Note that $\otimes$ and $\boxtimes$ commute: 
\[ X_1\boxtimes Y_1 \,\otimes\, X_2\boxtimes Y_2=
  (X_1\otimes X_2) \,\boxtimes\, (Y_1\otimes Y_2)=
   X_1X_2 \,\boxtimes\, Y_1 Y_2. \]


\subsection{Frobenius algebras and 2-categories}
\bdefin \label{d-Frob}
Let $\2A$ be a (strict) monoidal category. A Frobenius algebra in $\2A$ is a
quintuple $\5Q=(Q,v,v',w,w')$, where $Q$ is an object in $\2A$ and 
$v:\11\rightarrow Q, v': Q\rightarrow\11, w: Q\rightarrow Q^2, w':Q^2\rightarrow Q$ 
are morphisms satisfying the following conditions:
\begin{equation} \label{W1}
   w\otimes\id_Q\mcirc w=\id_Q\otimes w\mcirc w, 
\end{equation}
\begin{equation} \label{W1'}
   w' \mcirc w'\otimes\id_Q=w' \mcirc \id_Q\otimes w', 
\end{equation}
\begin{equation} \label{W2}
   v'\otimes\id_Q\mcirc w = \id_Q = \id_Q\otimes v'\mcirc w,
\end{equation}
\begin{equation} \label{W2'}
   w' \mcirc v\otimes\id_Q = \id_Q = w' \mcirc \id_Q\otimes v,
\end{equation}
\begin{equation} \label{W3}
   w'\otimes\id_Q \mcirc \id_Q\otimes w = w\mcirc w' = \id_Q\otimes w'\mcirc w\otimes \id_Q. 
\end{equation} 
A Frobenius algebra $\5Q$ in a $\7F$-linear category is canonical if
\bea w'\circ w &=& \lambda_1 \id_Q, \label{W4} \\
   v'\mcirc v &=& \lambda_2 \id_\11\label{W5}\eea
with $\lambda_1,\lambda_2\in\7F^*$. $\5Q$ is normalized if $\lambda_1=\lambda_2$.
\edefin

Let $X$ be an object in a spherical category $\2A$. Then the quintuple
\[ (X\ol{X},\ \ve(X),\ \ol{\ve}(X),\ \id_X\otimes\ve(\ol{X})\otimes\id_{\ol{X}},\ 
   \id_X\otimes\ol{\ve}(\ol{X})\otimes\id_{\ol{X}}) \]
is easily seen to be a normalized canonical Frobenius algebra in $\2A$. The following
theorem, which combines the Theorems 3.12 and 5.13 from \cite{mue09}, shows that in fact
every canonical Frobenius algebra in a tensor category arises in this way, provided one is
ready to embed the category as a corner into a bicategory.

\btheor \label{Main}
Let $\2A$ be a strict $\7F$-linear tensor category and $\5Q=(Q,v,v',w,w')$ a canonical
Frobenius algebra in $\2A$. Then:
\begin{itemize}
\item[(i)] There is a bicategory $\2E$ such that
\begin{enumerate}
\item The sets of 2-morphisms in $\2E$ are finite dimensional $\7F$-vector spaces and the
horizontal and vertical compositions are bilinear.
\item Idempotent 2-morphisms in $\2E$ split. 
\item $\displaystyle \obj\,\2E=\{\6A,\6B\}$.
\item There is an equivalence
$\displaystyle\HOM_\2E(\6A,\6A)\stackrel{\otimes}{\simeq}\ol{\2A}^p$ of tensor categories,
thus an equivalence
$\displaystyle\HOM_\2E(\6A,\6A)\stackrel{\otimes}{\simeq}\2A$ if $\2A$ has subobjects.
\item There are 1-morphisms $J: \6B\rarr \6A$ and $\oj : \6A\rarr \6B$ such that
$Q=J\oj $.  
\item $J$ and $\oj $ are mutual two-sided duals, i.e.\ there are 2-morphisms
\[ e_J: \11_\6A\rarr J\oj,\ \ \ve_J: \11_\6B\rarr\oj J,\ \ d_J: \oj J\rarr\11_\6B,\ \ 
   \eta_J: J\oj\rarr\11_\6A \]
satisfying the usual equations.
\item We have
\bean v &=& e_J:\ \11_\6A\rarr Q=J\oj , \\
   v' &=& \eta_J:\ Q=J\oj \rarr\11_\6A, \\
   w &=& \id_J\otimes\ve_J\otimes\id_{\oj }:\ Q=J\oj \rarr J\oj J\oj =Q^2, \\
   w' &=& \id_J\otimes d_J\otimes\id_{\oj }:\ Q^2=J\oj J\oj \rarr J\oj =Q
\eean
and therefore $d_J\circ\ve_J=\lambda_1\id_{\11_\6B}$, $\eta_J\circ e_J=\lambda_2\id_{\11_\6A}$.
\item $\2E$ is uniquely determined up to equivalence by the above properties. Isomorphic
Frobenius algebras $\5Q, \tilde{\5Q}$ give rise to isomorphic bicategories 
$\2E, \tilde{\2E}$. 
\end{enumerate}
\item[(ii)] If $\2A$ has direct sums then $\2E$ has direct sums of 1-morphisms.
\item[(iii)] If the multiplicity of $\11$ in $Q$ is exactly one (it is at least one due to
the existence of $v,v'$) then $J, \oj, \11_\6B$ are simple. (There is a weaker condition
implying only simplicity of $\11_\6B$.)
\item[(iv)] If $\7F=\7C$ and $\2A$ has a positive $*$-operation then $\2E$ has a positive
$*$-operation and is semisimple.
\item[(v)] If $\2A$ is strict spherical and $\5Q$ satisfies (iii) and is normalized 
then $\2E$ is spherical. If, furthermore, $\2A$ is semisimple and $\7F$ is algebraically
closed then $\2E$ is semisimple. 
\item[(vi)] If (iv) or (v) apply then the tensor category $\2B=\HOM_\2E(\6B,\6B)$
satisfies $\dim\2B=\dim\2A$. 
\end{itemize}
\etheor

\brem 1. If two tensor categories $\2A, \2B$ are `corners' of a 2-category as above we
call them weakly monoidally Morita equivalent. This is an equivalence relation which is
considerably weaker than the usual equivalence, yet it implies that $\2A$ and $\2B$ have
the same dimension and define the same triangulation invariant \cite{bw1, gk} for
3-manifolds. See \cite{mue09} for the details. 

2. Unfortunately, the above statement of the theorem will not be sufficient for our
purposes since beginning in Subsection \ref{fff} we will make use of the concrete structure
of the bicategory $\2E$, which is explicitly constructed in the proof of Theorem
\ref{Main}. This is not the place to explain the latter which occupies the larger part of 
\cite{mue09}. We can only hope that the above statement of the theorem and its r\^{o}le in
this paper are sufficient to motivate the reader to acquire some familiarity with
\cite{mue09}. 
\erem


\section{The Quantum Double of a Tensor Category}\label{sec-semisim}
\subsection{On half braidings}
We begin with the definition of the quantum double $\2Z_1(\2C)$ of a (strict) monoidal
category $\2C$. 
\bdefin Let $\2C$ be a strict monoidal category and let $X\in\2C$. A half braiding
$e_X$ for $X$ is a family $\{e_X(Y)\in\Hom_\2C(XY,YX),\ Y\in\2C\}$ of morphisms satisfying 
\begin{itemize}
\item[(i)] Naturality w.r.t.\ the argument in brackets, i.e.
\begin{equation}t\otimes\id_X\mcirc e_X (Y)= e_X(Z)\mcirc\id_X
  \otimes t \quad\forall t: Y\rarr Z.
\label{hb-i}\end{equation}
\item[(ii)] The braid relation
\begin{equation} e_X(Y\otimes Z)=\id_{Y}\otimes e_X(Z) \mcirc
    e_X(Y)\otimes\id_{Z} \quad \forall Y,Z\in\2C. \label{hb-ii}\end{equation}
\item[(iii)] All $e_X(Z)$ are isomorphisms.
\item[(iv)] Unit property:
\begin{equation} e_X(\11)=\id_X. \label{hb-iv}\end{equation}
\end{itemize}
\label{hb}\edefin

\blemma \label{iii-iv}
Let $\{ e_X(Y), Y\in\2C \}$ satisfy (i) and (ii). Then (iii)$\impl$(iv). If (iv) holds 
and $Y$ has a right dual $Y^*$ then $e_X(Y)$ is invertible.
\elemma
\prf Considering (\ref{hb-ii}) with $Y=Z=\11$ gives $e_X(\11)=e_X(\11)^2$. Thus (iii)
implies $e_X(\11)=\id_X$. Let $Y^*$ a right dual of $Y$ with
$\ve_Y:\11\rarr Y^*\otimes Y, \eta_Y: Y\otimes Y^*\rarr\11$. Then using (i) and (iv) we
find
\[
\begin{tangle}
\coev\mobj{\eta_Y}\Step\xx\mobj{e_X(Y)}\\
\id\Step\xx\step[-.3]\mobj{e_X(Y^*)}\step[2.3]\id\\
\id\Step\id\Step\ev\step[-.5]\mobj{\ve_Y}\\
\object{Y}\Step\object{X}
\end{tangle}
\quad=\quad
\begin{tangle}
\coev\Step\id\step[1.5]\id\\
\hh\id\Step\id\mobj{Y^*}\Step\id\step\frabox{e_X(\11)}\\
\id\Step\ev\step[1.5]\id\\
\object{Y}\step[5.5]\object{X}
\end{tangle}
\quad=\quad
\begin{tangle}
\id\Step\id\\
\id\Step\id\\
\object{Y}\Step\object{X}
\end{tangle}
\]
Thus $e_X(Y)$ has a right inverse, which by a similar computation is seen to be also a
left inverse.
\qed

For later use we record the following alternative characterization of half braidings.
\blemma \label{halfbr}
Let $\2C$ be semisimple and $\{X_i, i\in \Gamma\}$ a basis of simple objects. Let $Z\in\2C$.
Then there is a one-to-one correspondence between (i) families of morphisms 
$\{ e_Z(X_i)\in\Hom_\2C(Z X_i,X_iZ),\ i\in \Gamma\}$ such that
\begin{multline} t \otimes\id_Z \mcirc  e_Z(X_k)=\id_{X_i}\otimes e_Z(X_j)
   \mcirc e_Z(X_i)\otimes\id_{X_j}\mcirc\id_Z\otimes t \\
   \forall i,j,k\in \Gamma,\ t\in\Hom_\2C(X_k,X_iX_j), 
\label{hb-v}\end{multline}  
and (ii) families of morphisms $\{ e_Z(X)\in\Hom_\2C(Z X,XZ),\ X\in\2C\}$ satisfying 1.\
and 2.\ from Definition \ref{hb}. All $e_Z(X), X\in\2C$ are isomorphisms iff all $e_Z(X_i),
i\in \Gamma$ are isomorphisms. 
\label{l17}\elemma
\prf (ii)$\impl$(i). Obvious: just restrict $e_Z(\cdot)$ to 
$X\in\{X_i, i\in \Gamma\}$. Then (\ref{hb-i}, \ref{hb-ii}) imply (\ref{hb-v}).\\
(i)$\impl$(ii). Let $X\cong\bigoplus_i n_i X_i$ and let
$\{ x_i^\alpha, \alpha=1,\ldots, n_i\}$, $\{ {x'_i}^\alpha, \alpha=1,\ldots, n_i\}$ be 
dual bases in $\Hom_\2C(X_i,X)$ and $\Hom_\2C(X, X_i)$, respectively. Then define
\[  e_Z(X)=\sum_{i\in \Gamma}\sum_{\alpha=1}^{n_i} x_i^\alpha\otimes\id_Z
   \mcirc e_Z(X_i)\mcirc\id_Z\otimes {x'}_i^\alpha. \]
Independence of  $e_Z(X)\in\Hom_\2C(Z X,XZ)$ of the choice of the $x_i^\alpha$ follows
from duality of the bases $\{x_i^\alpha\}, \{{x'_i}^\alpha\}$. In order to prove
naturality (\ref{hb-ii}) consider $Y\cong\bigoplus_i m_i X_i$ and corresponding
intertwiners $y_i^\alpha, {y'_i}^\alpha$ and let $t\in \Hom_\2C(X,Y)$. Then 
${y'_j}^\beta tx_i^\alpha\in\Hom(X_i,X_j)$, which vanishes if $i\ne j$. Thus 
\[ t=\sum_{i\in \Gamma}\sum_{\alpha=1}^{n_i}\sum_{\beta=1}^{m_i}
   c(i,\alpha,\beta) y_i^{\beta} {x'}_i^\alpha, \]
where $c(i,\alpha,\beta)\in\7F$. Therefore,
\[ t\otimes\id_Z\mcirc e_Z(X)=  \sum_{i\in \Gamma}\sum_{\alpha=1}^{n_i}\sum_{\beta=1}^{m_i}
  c(i,\alpha,\beta) \ y_i^{\beta}\otimes\id_Z\mcirc e_Z(X_i)\mcirc\id_Z \otimes
  {x'}_i^\alpha, \] 
which coincides with $e_Z(Y)\mcirc\id_Z\otimes t$.
If now $t\in\Hom_\2C(X_k, X_iX_j)$ then naturality implies
$t\otimes\id_Z\circ e_Z(X_k)= e_Z(X_iX_j)\mcirc\id_Z\otimes t$.
Together with (\ref{hb-v}) this implies 
\[  e_Z(X_iX_j)\mcirc\id_Z\otimes t=\id_{X_i}\otimes e_Z
  (X_j)\mcirc e_Z(X_i)\otimes\id_{X_j}\mcirc\id_Z\otimes t, \]
and since this holds for all $t\in\Hom_\2C(X_k,X_iX_j)$ (\ref{hb-ii}) follows.
(This is a consequence of
\[ \sum_{k\in \Gamma}\sum_{\alpha=1}^{N_{ij}^k} t_k^\alpha\circ {t'}_k^\alpha=
  \id_{X_iX_j}, \]
where the $\{t_k^\alpha, \alpha=1,\ldots,N_{ij}^k\}$ are  bases in 
$\Hom_\2C(X_k,X_iX_j)$.) \qed


\subsection{Elementary properties of the quantum double}
\bdefin The center $\2Z_1(\2C)$ of a strict monoidal category $\2C$ has as objects pairs
$(X,  e_X)$, where $X\in\2C$ and $ e_X$ is a half braiding. The morphisms are given by 
\bea \lefteqn{ \Hom_{\2Z_1(\2C)}((X, e_X),(Y, e_Y)= \{ t\in\Hom_\2C(X,Y) \ | } \nn \\
 && \quad\quad \id_X\otimes t \mcirc e_X(Z)= e_Y(Z)\mcirc t\otimes\id_X \quad 
\forall Z\in\2C \}. \label{dmor}\eea
The tensor product of objects is given by 
$(X, e_X)\otimes(Y, e_Y)=(XY, e_{XY})$, where
\begin{equation} e_{XY}(Z)= e_X(Z)\otimes\id_Y\mcirc\id_X\otimes e_Y(Z).
\label{dprod}\end{equation}
The tensor unit is $(\11,e_\11)$ where $e_\11(X)=\id_X$. The composition and tensor
product of morphisms are inherited from $\2C$. The braiding is given by 
\[  c((X, e_X),(Y, e_Y))=  e_X(Y). \]
\label{defdouble}
\edefin

For the proof that $\2Z_1(\2C)$ is a strict braided tensor category we refer to \cite{ka}. 
The following is immediate from the definition of the center $\2Z_1(\2C)$:

\blemma If $\2C$ is $\7F$-linear then so is $\2Z_1(\2C)$. 
If the unit $\11$ of $\2C$ is simple, then $\11_{\2Z_1(\2C)}$ is simple.
\elemma

In \cite[Proposition 1]{s} it is proven that the center of an abelian monoidal category is
abelian. In this paper we do not use the language of abelian categories since the notions
of (co)kernels are not really needed. (Yet semisimple categories are abelian if we assume
existence of a zero object.) Therefore, we prove two lemmas which show that the center
construction behaves nicely w.r.t.\ direct sums and subobjects. The first result is
contained in \cite{s}, but we repeat it for the sake of completeness.

\blemma If $\2C$ has direct sums then also $\2Z_1(\2C)$ has direct sums. \elemma 
\prf Let $(Y, e_Y), (U, e_U)$ be objects in $\2Z_1(\2C)$. Let
$\2C\ni Z\cong Y\oplus U$ with morphisms $v\in\Hom_\2C(Y,Z), w\in\Hom_\2C(U,Z)$,
$v'\in\Hom_\2C(Z,Y), w'\in\Hom_\2C(Z,U)$ satisfying 
$v'\circ v=\id_Y,w'\circ w=\id_U, v\circ v'+w\circ w'=\id_Z$.
Defining $ e_Z(X)\in\Hom_\2C(Z X,XZ)$ for all $X\in\2C$ by
\[  e_Z(X)= \id_X\otimes v \mcirc e_Y(X)\mcirc v'\otimes\id_X +
   \id_X\otimes w \mcirc e_U(X)\mcirc w'\otimes\id_X, \]
we claim that $(Z,  e_Z)$ is an object of $\2Z_1(\2C)$ and 
\begin{equation}(Z, e_Z)\cong(Y, e_Y)\oplus (U, e_U). 
\label{sum}\end{equation}
Naturality of $ e_Z(X)$ w.r.t.\ $X$ is obvious, and (\ref{hb-ii}) is very
easily verified using $v'\circ w=0$. Finally, we have
\[  e_Z(X)\mcirc v\otimes\id_X= \id_X\otimes v \mcirc e_Y(X), \]
which is just the statement that 
$v\in\Hom_{\2Z_1(\2C)}((Y, e_Y),(Z, e_Z)$. The analogous statement holding
for $v', w, w'$, (\ref{sum}) follows. \qed

\blemma If $\2C$ has subobjects then also $\2Z_1(\2C)$ has subobjects.  
\label{subobj}\elemma
\prf Let $(Y, e_Y)\in \2Z_1(\2C)$ and let $e$ be an idempotent in 
$\End_{\2Z_1(\2C)}((Y,e_Y))$. By definition of $\2Z_1(\2C)$ this means that $e$ is an
idempotent in $\End_\2C(Y)$ such that 
\begin{equation}\id_X\otimes e\mcirc e_Y(X)= e_Y(X)\mcirc e\otimes\id_X \quad 
  \forall X\in\2C. 
\label{ez}\end{equation}
Since $\2C$ has subobjects there are $U\in\2C$ and $v\in\Hom_\2C(U,Y), v'\in\Hom_\2C(Y,U)$ 
such that $v\circ v'=e$ and $v'\circ v=\id_U$. Defining
\[  e_U(X)=\id_X\otimes v'\mcirc  e_Y(X)\mcirc v\otimes\id_X
  \in\Hom_\2C(U X,XU),   \quad X\in\2C, \]
naturality w.r.t.\ $X$ is again obvious. Now,
\bean \lefteqn{ e_U(X_1X_2) =\id_{X_1X_2}\otimes v'\mcirc 
   e_Y(X_1X_2)\mcirc v\otimes\id_{X_1X_2}} \\
  &&= \id_{X_1X_2}\otimes v'  \mcirc\id_{X_1}\otimes e_Y(X_2) \mcirc
    e_Y(X_1)\otimes\id_{X_2}\mcirc v\otimes\id_{X_1X_2} \\
  &&= \id_{X_1X_2}\otimes v'  \mcirc\id_{X_1}\otimes e_Y(X_2) \mcirc
 \id_{X_1}\otimes v\otimes\id_{X_2} \\
  && \quad\mcirc\id_{X_1}\otimes v'\otimes\id_{X_2} \mcirc
    e_Y(X_1)\otimes\id_{X_2}\mcirc v\otimes\id_{X_1X_2} \\
  &&= \id_{X_1}\otimes e_U(X_2) \mcirc   e_U(X_1)\otimes\id_{X_2},
\eean
whereby $ e_U$ is a half braiding and $(U, e_U)$ an object in $\2Z_1(\2C)$. 
We used $v\mcirc v'=e$, (\ref{ez}) and $e\mcirc v=v\mcirc v'\mcirc v=v$.
Using the same facts we finally compute
\[ \id_X\otimes v\mcirc e_U(X)=\id_X\otimes v\mcirc\id_X\otimes v'
  \mcirc  e_Y(X)\mcirc v\otimes\id_X=
   e_Y(X)\mcirc v\otimes\id_X. \]
Thus $v\in\Hom_{\2Z_1(\2C)}((U, e_U),(Y, e_Y))$ and we have 
$(U, e_U)\prec(Y, e_Y)$. \qed

\blemma \label{l3}
Let $\2C$ be pivotal and $e_Y$ a half braiding satisfying (i-iv). Then
\begin{equation} e_Y(\ol{X})=\id_{\ol{X}Y}\otimes \ol{\ve}(X) \mcirc
   \id_{\ol{X}}\otimes e_Y(X)^{-1}\otimes\id_{\ol{X}}\mcirc \ve(\ol{X})
   \otimes \id_{Y\ol{X}}. \label{vebar}\end{equation}
\elemma
\prf By naturality and the braid relation we have
\[ \ve(\ol{X})\otimes\id_Y = e_Y(\ol{X}X) \mcirc\id_Y\otimes \ve(\ol{X})=
  \id_{\ol{X}}\otimes e_Y(X)\mcirc e_Y(\ol{X})\otimes\id_X\mcirc
  \id_Y\otimes \ve(\ol{X}) \]
and using the invertibility of $e_Y(X)$ we get
\[ \id_{\ol{X}}\otimes e_Y(X)^{-1}\mcirc \ve(\ol{X})\otimes\id_Y  = 
   e_Y(\ol{X})\otimes\id_X \mcirc \id_Y\otimes \ve(\ol{X}). \]
Now (\ref{vebar}) follows by a use of the duality equations, see, e.g., 
\cite[Subsection 2.3]{mue09}.
\qed 

\bprop \label{hb-conj}
Let $\2C$ be (strict) pivotal. Then also $\2Z_1(\2C)$ is (strict) pivotal, the
dual $\ol{(Y, e_Y)}$ being given by $(\ol{Y}, e_{\ol{Y}})$, where $e_{\ol{Y}}(X)$ is
defined by 
\[
\begin{diagram}
\ol{Y}\otimes X & \rTo^{\id_{\ol{Y}X}\otimes \ve(Y)} & \ol{Y}\otimes X\otimes Y\otimes\ol{Y} &
\rTo^{\id_{\ol{Y}}\otimes e_Y(X)^{-1}\otimes\id_{\ol{Y}}} & \ol{Y}\otimes Y\otimes X\otimes\ol{Y} & \rTo 
\end{diagram}
\]
\be
\begin{diagram}
   & \rTo^{\ol{\ve}(\ol{Y})\otimes\id_{X\ol{Y}}} & X\otimes\ol{Y}
\end{diagram}
\label{eybar}\end{equation}
The evaluation and coevaluation maps are inherited from $\2C$:
\[ \ve((Y,  e_Y))=\ve(Y), \quad \ol{\ve}((Y,  e_Y))=\ol{\ve}(Y). 
\]
If $\2C$ is spherical then also $\2Z_1(\2C)$ is spherical.
\eprop
\prf We begin by showing that $e_{\ol{Y}}(\cdot)$ is a half braiding for $\ol{Y}$. By
construction we have $ e_{\ol{Y}}(X)\in\Hom_\2C(\ol{Y}X,X\ol{Y})$, and naturality
w.r.t.\ $X$ follows easily from the corresponding property for $ e_Y$. Now
\bean \lefteqn{ e_{\ol{Y}}(X_1X_2) = \ol{\ve}(\ol{Y})\otimes\id_{X_1X_2\ol{Y}}
  \mcirc \id_{\ol{Y}}\otimes e_Y(X_1X_2)^{-1}\otimes\id_{\ol{Y}} 
  \mcirc \id_{\ol{Y}X_1X_2}\otimes \ve(Y) } \\
 &&= \ol{\ve}(\ol{Y})\otimes\id_{X_1X_2\ol{Y}}\mcirc \id_{\ol{Y}}\otimes
   e_Y(X_1)^{-1}\otimes\id_{X_2\ol{Y}} \\
  && \quad\quad\mcirc\id_{\ol{Y}X_1}
  \otimes e_Y(X_2)^{-1}\otimes\id_{\ol{Y}}\mcirc
  \id_{\ol{Y}X_1X_2}\otimes \ve(Y) \\
 &&= \id_{X_1}\otimes \ol{\ve}(\ol{Y})\otimes\id_{X_2\ol{Y}} \mcirc
  \id_{X_1\ol{Y}}\otimes e_Y(X_2)^{-1}\otimes\id_{\ol{Y}}\mcirc
  \id_{X_1\ol{Y}X_2}\otimes \ve(Y) \\
  && \quad\quad\mcirc  \ol{\ve}(\ol{Y})\otimes\id_{X_1\ol{Y}X_2}\mcirc
  \id_{\ol{Y}}\otimes e_Y(X_1)^{-1}\otimes\id_{\ol{Y}X_2}\mcirc
  \id_{\ol{Y}X_1}\otimes\ve(Y)\otimes\id_{X_2}\\
 &&=\id_{X_1}\otimes e_{\ol{Y}}(X_2)\mcirc e_{\ol{Y}}(X_1)\otimes
   \id_{X_2}.
\eean
In the third equality we have used the duality equation 
$\id_Y\otimes \ol{\ve}(\ol{Y})\mcirc \ve(Y)\otimes\id_Y=\id_Y$ and the interchange law. 

In view of the definition (\ref{eybar}) of $e_{\ol{Y}}(X)$ together with
$e_Y(\11)=\id_Y$ and the duality equation we have $e_{\ol{Y}}(\11)=\id_{\ol{Y}}$. Now
Lemma \ref{iii-iv} implies invertibility of $e_{\ol{Y}}(X)$ for all $X$.

It remains to show that $\ve(Y):\11_\2C\rarr Y\otimes\ol{Y}$ is actually in
\[ \Hom_{\2Z_1(\2C)}(\11_{\2Z_1(\2C)},(Y,e_Y)\otimes(\ol{Y},e_{\ol{Y}}))=
   \Hom_{\2Z_1(\2C)}((\11,\id),(Y\ol{Y},e_{Y\ol{Y}})), \]
which in view of (\ref{dmor}) amounts to 
\begin{equation} \label{conjd}
\begin{tangle}
\object{Z}\Step\object{Y}\Step\object{\ol{Y}}\\
\id\Step\ev\\
\object{Z}
\end{tangle}
\quad=\quad
\begin{tangle}
\object{Z}\Step\object{Y}\Step\object{\ol{Y}}\\
\xx\mobj{e_Y(Z)}\Step\id\\
\id\Step\xx\mobj{e_{\ol{Y}}(Z)}\\
\ev\Step\id\\
\step[4]\object{Z}
\end{tangle}
\end{equation}
With the definition (\ref{eybar}) of $e_{\ol{Y}}$ the right hand side equals
\[
\begin{tangle}
\object{Z}\Step\object{Y}\Step\object{\ol{Y}}\\
\xx\mobj{e_Y(Z)}\Step\id\\
\id\Step\d\step\id\\
\hh\id\step\hcoev\step\id\step\id\\
\id\step\id\step\hx\step\id\mobj{e_Y(Z)^{-1}}\\
\hh\hev\step\id\step\hev\\
\Step\object{Z}
\end{tangle}
\quad\quad=\quad\quad
\begin{tangle}
\object{Z}\Step\object{Y}\Step\object{\ol{Y}}\\
\xx\mobj{e_Y(Z)}\Step\id\\
\x\step[-.2]\mmobj{e_Y(Z)^{-1}}\step[2.2]\id\\
\id\Step\ev\\
\object{Z}
\end{tangle}
\]
which coincides with the left hand side of (\ref{conjd}) as desired. That $\ol{\ve}(Y)$ is
a morphism in $\2Z_1(\2C)$ is shown analogously. The composition of morphisms being the
same in $\2Z_1(\2C)$ as in $\2C$, $\ve(X), \ol{\ve}(X)$ inherit from $\2C$ all equations
needed to make $\2Z_1(\2C)$ pivotal (spherical). If the pivotal structure of $\2C$ is strict
then the same clearly holds for $\2Z_1(\2C)$.
\qed\\


\subsection{Semisimplicity of $\2Z_1(\2C)$}
\blemma Let $\2C$ be semisimple spherical with simple unit. We assume that there are
only finitely many simple objects and that $\dim\2C\ne 0$. Let
$(X,e_X),(Y,e_Y)\in\2Z_1(\2C)$. Then the map
$E_{X,Y}: \Hom_\2C(X,Y)\rightarrow\Hom_\2C(X,Y)$ defined by 
\[ E_{X,Y}(t)=(\dim\2C)^{-1}\sum_{i\in \Gamma} d_i \ \ 
\begin{tangle}
\step\object{\ol{\ve}(X_i)}\Step\object{Y} \\
\step\hh\hcoev\step\id\\
\step\id\step\hxx\mobj{e_Y(\ol{X_i})}\\
\obj{X_i}\step\id\step\O t\step\id\hstep\obj{\ol{X_i}} \\
\step[-1.2]\mobj{e_X(X_i)}\step[2.2]\hxx\step\id \\
\step\hh\id\step\hev\\
\step\object{X}\step{\ve(X_i)}
\end{tangle}
\]
is a projection onto $\Hom_{\2Z_1(\2C)}((X,e_X),(Y,e_Y))\subset\Hom_\2C(X,Y)$. Here
$\{ X_i, i\in \Gamma\}$ is a basis of simple objects and we abbreviate $d_i=d(X_i)$.
The family of maps $E_{X,Y}$ is a conditional expectation in the sense that
\begin{equation}E_{X,T}(c\circ b\circ a)=c\circ E_{Y,Z}(b)\circ a
\label{cexp}\end{equation} 
if $a\in\Hom_{\2Z_1(\2C)}((X,e_X),(Y,e_Y)), b\in\Hom_\2C(Y,Z),
c\in\Hom_{\2Z_1(\2C)}((Z,e_Z),(T,e_T))$.
\label{condexp1}\elemma
\prf We compute
\[ \dim\2C \cdot \id_Z\otimes E_{X,Y}(t) \mcirc e_X(Z) = \sum_{i\in \Gamma} d_i \ \
\begin{tangle}
\object{Z}\step[3]\object{Y}\\
\hh\id\hstep\obj{\ol{\ve}(X_i)}\step[2.5]\id\\
\hh\id\step\hcoev\step\id\\
\id\step\id\step\hxx\\
\id\obj{X_i}\step\id\step\O t\step\id\obj{\ol{X_i}} \\
\id\step\hxx\step\id \\
\hh\id\step\id\step\hev\obj{\ve(X_i)}\\
\hxx\\
\object{X}\step\object{Z}
\end{tangle}
\quad=\ \ \sum_i\sum_{j,\alpha} d_i \quad \quad
\begin{tangle}
\object{Z}\step[3]\object{Y}\\
\hh\id\step\hcoev\step\id\\
\step[-0.5]\obj{X_j}\hstep\hcu\step\hxx\\
\step[-1]\mobj{{p'}_i^{j,\alpha}}\step\hcd\step\O t\step\id\obj{\ol{X_i}} \\
\id\step\hxx\step\id \\
\hh\id\step\id\step\hev\\
\hxx\\
\object{X}\step\object{Z}
\end{tangle}
\]

\[ = 
\sum_i\sum_{j,\alpha} d_i \ \ 
\begin{tangle}
\object{Z}\step[3]\object{Y}\\
\hh\id\step\hcoev\step\id\\
\step[-0.5]\obj{X_j}\hstep\hcu\step\hxx\\
\dh\step\O t\step\hd \\
\step\hxx\step[1.5]\id\obj{\ol{X_i}} \\
\hdd\hstep\hcd\step\id\\
\hh\hstep\id\step\id\step\hev\\
\hstep\object{X}\step\object{Z}
\end{tangle}
\quad = \ \ 
\sum_{i,j,\alpha} \frac{d_id_j}{d(Z)}\ \ 
\begin{tangle}
\step\hstep\object{Z}\step[1.5]\object{Y}\\
\step[-.2]\mobj{{q'}_{i,j}^\alpha}\step[1.2]\hcd\step\id\\
\obj{X_j}\step\id\step\hxx\\
\step\id\step\O t\step\id\obj{\ol{X_i}} \\
\step\hxx\step\id \\
\hh\step\id\step\hcu\mobj{q_{i,j}^\alpha}\\
\hh\step\id\step[1.5]\id\\
\step\object{X}\step[1.5]\object{Z}
\end{tangle}
\quad=\quad
\sum_j\sum_{i,\alpha} d_j \ \ 
\begin{tangle}
\step[2.5]\object{Z}\step\object{Y}\\
\hh\hstep\hcoev\step\id\step\id\\
\hstep\id\step\hcu\hstep\ddh\\
\hstep\id\step[1.5]\hxx\\
\obj{X_j}\dh\step\O t\step\hd\obj{\ol{X_i}} \\
\step\hxx\step\hcd\mobj{{r'}_j^{i\alpha}}\\
\hh\step\id\step\hev\step\id\\
\step\object{X}\step[3]\object{Z}
\end{tangle}
\]

\[ =
\sum_j\sum_{i,\alpha} d_j \ 
\begin{tangle}
\step[3]\object{Z}\step\object{Y}\\
\step[3]\hxx\\
\hh\step\hcoev\step\id\step\id\\
\step\id\step\hxx\step\id\\
\obj{X_j}\step\id\step\O t\step\hcu\obj{\ol{X_i}}\\
\step\hxx\step\hcd\\
\hh\step\id\step\hev\step\id\\
\step\object{X}\step[3]\object{Z}
\end{tangle}
\quad = \quad
\sum_i d_j\ 
\begin{tangle}
\step[3]\object{Z}\step\object{Y}\\
\step[3]\hxx\\
\hh\step\hcoev\step\id\step\id\\
\step\id\step\hxx\step\id\\
\obj{X_j}\step\id\step\O t\step\id\obj{\ol{X_j}}\step\id \\
\step\hxx\step\id\step\id \\
\hh\step\id\step\hev\step\id\\
\step\object{X}\step[3]\object{Z}
\end{tangle}
\quad = \dim\2C \cdot e_Y(Z)\mcirc E_{X,Y}(t)\otimes\id_{Z} 
\]
Here $\{ p_i^{j,\alpha}, \alpha=1,\ldots,N_{Z,X_i}^{X_j}\}$ is, for every $j\in \Gamma$, a
basis in $\Hom_\2C(X_j,ZX_i)$ with dual basis $\{{p'}_i^{j,\alpha}\}$ such that 
${p'}_i^{j,\alpha}\circ p_i^{k,\beta}=\delta_{j,k}\delta_{\alpha\beta}\,\id_{X_j}$ and
$\id_{ZX_i}=\sum_{j,\alpha} p^{j,\alpha}_i\circ {p'}^{j,\alpha}_i$.
We used the fact that $e_X(\cdot),e_Y(\cdot)$ are half-braidings, i.e.\ natural w.r.t.\
the second argument. Furthermore, the basis $\{ q_{i,j}^\alpha \}$ in
$\Hom_\2C(Z,X_jX_{\ol{\imath}})$ and its dual basis are normalized such that
$tr_Z({q'}_{i,j}^\beta\circ q_{i,j}^\alpha)=d(Z)\delta_{\alpha\beta}$. We used that a
basis together with its dual can be replaced  by another one provided the normalizations
are the same.

Since the above computation holds for all $Z\in\2C$ we conclude that $E_{X,Y}(t)$ is
in $\Hom_{\2Z_1(\2C)}((X,e_X),(Y,e_Y))$. The property (\ref{cexp}) for morphisms $a, c$ in
$\2Z_1(\2C)$ is obvious since by (\ref{dmor}) $a,c$ can be pulled through the half
braidings, changing the subscript of the conditional expectation $E$ appropriately. In
order to show that $E_{X,Y}$ is idempotent it thus suffices to show
$E_{X,X}(\id_X)=\id_X$, which follows from the definition of $\dim\2C$. 
\qed

\brem 1. Since the conditional expectations depend also on the half braidings we should
in principle denote them $E_{(X,e_X),(Y,e_Y)}$. We stick to $E_{X,Y}$ in order to keep the 
formulae simple.

2. The r\^{o}le of the assumption on the dimension is obvious: If $\dim\2C=0$ then
the map $E_{X,X}$ with the factor $(\dim\2C)^{-1}$ removed is identically zero on
$\End_{\2Z_1(\2C)}(X)$, thus we cannot use it to obtain a conditional expectation.

3. The proof uses a special instance of the `handle sliding' which is formalized in
\cite{brug3}. 
\erem

\blemma \label{condexp2}
For every $X\in\2C$ we have $tr_X\circ E_{X,X}=tr_X$, where $tr_X$ is the trace on 
$\End_\2C(X)$ provided by the spherical structure. 
\elemma
\prf Let $t\in\Hom_\2C(X,X)$. Using the fact that the spherical structure of $\2Z_1(\2C)$
is induced from $\2C$ we compute
\[ \dim\2C\,tr\circ E_{X,X} (t)=\sum_i d_i \quad
\begin{tangle}
\step[1.5]\object{\ol{\ve}(\ol{X})}\\
\Coev\Step \\
\hh\id\step\hcoev\step\id\\
\id\step\id\step\hxx\\
\id\obj{X_i}\step\id\step\O t\step\id\obj{\ol{X_i}} \\
\id\step\hxx\step\id \\
\hh\mobj{\ol{X}}\id\step\mobj{X}\id\step\hev\mobj{\ve(X_i)}\\
\hh\hev\obj{\ve(\ol{X})}
\end{tangle}
\quad = \sum_i d_i \quad
\begin{tangle}
\Coev \\
\hh\id\step\hcoev\step\id\\
\id\step\id\step\hxx\\
\id\obj{X_i}\step\id\step\O t\step\id\obj{\ol{X_i}} \\
\id\step\id\step\hx\\
\hh\id\step\hev\step\id\obj{X}\\
\Ev\Step\mobj{\ve(\ol{X})}
\end{tangle}
\ \ = \sum_i d_i \sum_{j,\alpha}\quad
\begin{tangle}
\coev\\
\id\step[1.5]\hcd \\
\id\step[1.5]\hxx\\
\id\step[1.5]\O t\step\id\hstep\obj{\ol{X_i}} \\
\id\step[1.5]\hx\obj{X}\\
\hh\id\step[1.5]\hcu\mobj{t^\alpha} \\
\hh\id\obj{\ol{X_j}}\Step\id\obj{X_j}\\
\ev\step[-1]\obj{\ve(\ol{X_j})}
\end{tangle}
\]
\[ \ \ = \ \sum_i d_i \sum_{j,\alpha} \quad
\begin{tangle}
\coev\\
\hh\id\step[1.5]\hcd \\
\obj{X_i}\id\step[1.5]\O t\step\id\step[.3]\obj{\ol{X_i}} \\
\hh\id\step[1.5]\hcu \\
\ev
\end{tangle}
\quad\quad =\ \sum_i d_i \quad\quad
\begin{tangle}
\Coev \\
\hh\id\step\hcoev\step\id\\
\step[-1]\obj{X_i}\step\id\step[.3]\obj{\ol{X}}\step[.7]\id\step\O t\step\id\step[.3]\obj{\ol{X_i}} \\
\hh\id\step\hev\step\id\\
\Ev\Step\mobj{\ve(X_i)}
\end{tangle}
\ \quad = \ \dim\2C \,tr(t).
\]
In the first step we have used Proposition \ref{hb-conj}, the second is based on standard
properties of categories with duals. In the next step we use that, given a basis 
$\{ t^\alpha \}$ in $\Hom_\2C(X_j,\ol{X_i}X)$ with dual basis $\{ \hat{t}^\alpha \}$,  
$\{ e_X(\ol{X_i}) \circ t^\alpha \}$ is a basis in $\Hom_\2C(X_j,X\ol{X_i})$ 
with dual basis $\{ \hat{t}^\alpha\circ e_X(\ol{X_i})^{-1} \}$. Replacing one basis by 
the other leaves the expression invariant. \qed\\

A trace on a finite dimensional $\7F$-algebra $A$ is a $\7F$-linear map $A\rarr\7F$ such
that $tr(ab)=tr(ba)$. It is non-degenerate if for every $a\ne 0$ there is $b$ such that
$tr(ab)\ne 0$. 

\blemma \label{semi}
Let $A$ be a finite dimensional $\7F$-algebra and $tr: A\rarr \7F$ a non-degenerate
trace. If $tr$ is vanishes on nilpotent elements then $A$ is semisimple. Conversely, every
trace (not necessarily non-degenerate) on a semisimple algebra vanishes on nilpotent
elements.
\elemma 
\prf Well known, but see, e.g., \cite{mue09}. \qed

\blemma \label{semisim1}
Let $A$ be a finite dimensional semisimple algebra over $\7F$ with a non-de\-ge\-ne\-rate
trace $tr: A\rightarrow \7F$. Let $B$ be a subalgebra containing the unit of $A$ and
assume there is a conditional expectation $E: A\rightarrow B$ (i.e.\ a linear map such
that $E(bab')=bE(a)b'$ for $a\in A, b,b'\in B$) such that $tr\circ E=tr$. Then $B$ is
semisimple. 
\elemma
\prf Let $0\ne x\in B$. By non-degeneracy of $tr$ there is $y\in A$ such that $tr(xy)\ne 0$. 
Now, using the properties of $E$ we compute $0\ne tr(xy)=tr\circ E(xy)=tr(xE(y))$. Since
$E(y)\in B$ we conclude that the restriction $tr_B\equiv tr\restr B$ is non-degenerate,
too. By Lemma \ref{semi} $tr$ vanishes on nilpotent elements, thus the same trivially
holds for $tr_B$. Now the other half of Lemma \ref{semi} applies and $B$ is semisimple.
\qed 

\brem Algebra extensions $A\supset B$ admitting a conditional expectation $E: A\rarr B$
(satisfying certain conditions) are well known as Frobenius extensions, cf., e.g.,
\cite{kad} and are called Markov extensions if there is an $E$-invariant trace on $A$.
\erem

Now we can put everything together:
\btheor \label{tsemisim}
Let $\7F$ be algebraically closed and $\2C$ a $\7F$-linear, spherical and
semisimple tensor category. We assume that there are only finitely many simple objects and
that $\dim\2C\ne 0$. Then the quantum double $\2Z_1(\2C)$ is spherical and semisimple. 
\etheor 
\prf Recall that by our definition of semisimplicity, $\2C$ has direct sums, subobjects
and a simple unit. By our earlier results also $\2Z_1(\2C)$ has these properties and is
spherical. It therefore only remains to show that the endomorphism algebra of every object
of $\2Z_1(\2C)$ is a multi matrix algebra.

Let $(X,e_X)\in\2Z_1(\2C)$. Then $\End_\2C(X)$ is a finite dimensional multi matrix
algebra by semisimplicity of $\2C$. The trace on $\End_\2C(X)$ provided by the duality
structure is non-degenerate, cf.\ e.g.\ \cite[Lemma 3.1]{gk}, and Lemmas \ref{condexp1},
\ref{condexp2} provide us with a trace preserving conditional expectation 
$E_X: \End_\2C(X)\rightarrow \End_{\2Z_1(\2C)}((X,e_X))$. Thus $\End_{\2Z_1(\2C)}((X,e_X))$ is
semisimple by Lemma \ref{semisim1} and therefore a multi matrix algebra since $\7F$ is
assumed algebraically closed. 
\qed

\brem 1. Even if $(X,e_X)$ is simple as an object of $\2Z_1(\2C)$ there is no reason why $X$ 
should be simple in $\2C$. Usually it is not. Since we do not know a priori which
non-simple objects of $\2C$ appear in the simple objects of $\2Z_1(\2C)$ we cannot dispense
with the assumption that $\2C$ has all finite direct sums as done, e.g., in \cite{t}.

2. Note that we do net yet know that $\2Z_1(\2C)$ has finitely many isomorphism classes of
simple objects. To show this will be our next aim.
\erem


\section{Weak Morita equivalence of $\2Z_1(\2C)$ and $\2C\boxtimes\2C^\op$}\label{sec-morita}
\subsection{A Frobenius Algebra in $\2C\boxtimes\2C^\op$}\label{ss-mcata}
Throughout this section $\2C$ will be a strict spherical tensor category with simple unit
over an algebraically closed field $\7F$. We require that $\2C$ is semisimple with
finite set $\Gamma$ of isomorphism classes of simple objects and $\dim\2C\ne 0$. The set
$\Gamma$ has a distinguished element $0$ representing the tensor unit and an involution
$i\mapsto\ol{\imath}$ which associates with every class the class of dual objects. We
choose representers $\{X_i,\ i\in \Gamma\}$ of these classes, which are arbitrary except
that we require $X_0=\11$. (We emphasize that we do not require $\ol{X_i}=X_{\ol{\imath}}$.
This can be achieved by a suitable strictification of the category if and only if all
self-dual objects are orthogonal \cite{bw1}. (The terms real vs.\ pseudo-real do not seem
appropriate if $\7F\ne\7C$.)) We choose once and for all square roots of the $d_i=d(X_i)$,
as well as $\lambda=\sqrt{\dim\2C}$  and $(\dim\2C)^{1/4}=\sqrt{\lambda}$. 
Let $N_{ij}^k$ be the dimension of the space $\Hom(X_k,X_iX_j)$, let 
$\{t_{ij}^{k\alpha},\ \alpha=1,\ldots,N_{ij}^k\}$ be a basis in $\Hom(X_k,X_iX_j)$ and
let $\{ {t'}_{ij}^{k\alpha} \}$ be the basis in $\Hom(X_iX_j,X_k)$ which is dual in the
sense of ${t'}_{ij}^{k\alpha}\circ t_{ij}^{k,\beta}=\delta_{\alpha\beta}$. Note that
this normalization of the dual basis differs from the one provided by the trace by a
factor of $d_k$. The present choice is more convenient since otherwise the dimensions
would appear in the equation 
\[ \sum_{k,\alpha} t_{ij}^{k\alpha}\circ {t'}_{ij}^{k\alpha}=\id_{X_iX_j}. \]
The choice of the of the square roots, the $X_i$ and of the bases $\{t_{ij}^{k\alpha}\}$
is immaterial but will be kept fixed throughout the rest of the paper, and the symbols  
$\Gamma, X_i, N_{ij}^k, t_{ij}^{k\alpha}$ will keep the above meanings. 

With these preparations we can embark on the 2-categorical approach to the quantum 
double. We define
\bea \2A & = & \2C\boxtimes\2C^\op, \\
 \hat{X}_i &=& X_i\boxtimes X_i^\op\in\obj\,\2A. \eea
By \cite[Lemma 2.9]{mue09}, $\2C^\op, \2C\otimes_\7F\2C^\op$ and $\2A$ are strict spherical
in a canonical way. Every $\hat{X}_i,\ i\in \Gamma$ is simple and if it is self-dual (i.e.\
if $i=\ol{\imath}$) then it is orthogonal irrespective of whether $X_i$ is orthogonal or
symplectic. 

The following is a very slight generalization of \cite[Proposition 4.10]{lre}.
\bprop Let $\7F$ be quadratically closed and let $\2C$ be $\7F$-linear semisimple
spherical with $\dim\2C\ne 0$. There is a normalized canonical Frobenius algebra
$\5Q=(Q,v,v',w,w')$ in $\2A=\2C\boxtimes\2C^\op$ (with $\lambda_1=\lambda_2=\lambda$) 
such that  
\begin{equation}Q\cong\bigoplus_{i\in \Gamma} \, \hat{X}_i. \label{gamma}\end{equation}
\eprop
\prf Clearly, $d(Q)=\dim\2C$. By definition of $Q$ there are morphisms 
\[ v_i\in\Hom(\hat{X}_i,Q),\quad v'_i\in\Hom(Q,\hat{X}_i),\ \ i\in \Gamma, \]
such that
\begin{equation}v'_i \circ v_j=\delta_{ij}\id_{\hat{X}_i},\quad \sum_i v_i\circ v'_i=\id_Q. 
\label{cuntz}\end{equation}
Defining $v=\lambda^{1/2} v_0, v'=\lambda^{1/2} v'_0$, (\ref{W5}) is trivial. With 
${t'}_{ij}^{k\alpha} \in\Hom_\2C(X_iX_j,X_k)\equiv\Hom_{\2C^\op}(X^\op_k,X^\op_iX^\op_j)$ 
the morphisms
\[ \hat{t}_{ij}^k= \sum_{\alpha=1}^{N_{ij}^k} t_{ij}^{k\alpha} \boxtimes 
   {t'}_{ij}^{k\alpha} \ \in \ \Hom_\2A(\hat{X}_k,\hat{X}_i\hat{X}_j) \]
are independent of the choices of the bases $\{t_{ij}^{k\alpha}\}$. Then
\begin{equation}w= \lambda^{-1/2} \ \sum_{i,j,k\in \Gamma} \sqrt{\frac{d_id_j}{d_k}} \,
   v_i \otimes v_j \mcirc \hat{t}_{ij}^k \mcirc v'_k \label{defW}\end{equation}
is in $\Hom_\2A(Q,Q^2)$, and $w'\in\Hom_\2A(Q^2,Q)$ is defined dually.
Eqs. (\ref{W2}, \ref{W2'}) of Definition \ref{d-Frob} are almost obvious. (Use
$N_{0i}^j=\delta_{ij}$). The proof that $w, w'$ satisfy (\ref{W1}, \ref{W1'}) and
(\ref{W3}) is omitted since it is entirely analogous to the one in \cite[p. 591]{lre}.
Finally, $w'\circ w=\lambda\id_Q$ is proven by a simple computation observing 
${\hat{t'}}^k_{ij} \circ \hat{t}_{ij}^k=N_{ij}^k\,\id_{\hat{X}_k}$ and using
\begin{equation}\label{ddN}
 \sum_{i,j} d_id_j N_{ij}^k=\sum_{i,j} d_id_j N_{j\ol{k}}^{\ol{\imath}}=\sum_j d_k d_j^2=
  d_k\dim\2C=d_k\lambda^2. \end{equation}
Thus $(Q,v,v',w,w')$ is a canonical Frobenius algebra in $\2A$. 
\qed

Thus Theorem \ref{Main} applies and yields a spherical bicategory $\2E$. ($\2E$ is strict as a
bicategory except for the existence of non-trivial unit constraints for $\11_\6B$ and
strict pivotal \cite{bw2} except for isomorphisms 
$\gamma_{X,Y}: \ol{Y}\circ\ol{X}\rarr\ol{X\circ Y}$ which are  non-trivial whenever
$\mbox{Ran}(Y)=\mbox{Src}(X)=\6B$.) In
particular, we have a spherical tensor category $\2B=\END_\2E(\6B)$.
In the rest of the paper $\2A, \5Q, \2E$ and $\2B$ will have the above meanings. By 
construction $Q$ contains the identity object of $\2A$ with multiplicity 1, thus $J, \oj $
and $\11_B$ are simple by \cite[Proposition 5.3]{mue09} and $d(J)=d(\oj)=\lambda$.
(The condition (iii) of that proposition can also easily be verified directly.) 

\blemma $\dim\2B=(\dim\2C)^2$. \elemma
\prf Follows from $\dim\2B=\dim\2A$ and $\dim\2A=(\dim\2C)^2$. The former is
\cite[Proposition 5.16]{mue09} and the latter is obvious since the simple objects of
$\2C\boxtimes\2C^\op$ are those of the form $X\boxtimes Y^\op$ with $X, Y$ simple. \qed 

In the sequel we will write $\11$ instead of $\11^\op$ in order to alleviate the
notation.


\subsection{A fully faithful tensor functor $F: \2Z_1(\2C)\rightarrow\2B$} \label{fff}
In this subsection we will construct a functor $F: \2Z_1(\2C)\rightarrow\2B$ and prove
that it is fully faithful and monoidal. This already implies that $\2Z_1(\2C)$ has finitely 
many isomorphism classes of simple objects, which is not at all obvious from
Definition \ref{defdouble}.

\blemma \label{11mor}
Let $X,Y\in\2C$. There is a one-to-one correspondence between morphisms 
$u\in \Hom_\2A((X\boxtimes\11)Q,Q(Y\boxtimes\11))\equiv\Hom_\2E(\oj (X\boxtimes\11)J,\oj (Y\boxtimes\11)J)$
and families $\{ u[i]\in\Hom_\2C(XX_i,X_iY),\ i\in \Gamma\}$. With $Z\in\2C$ and 
$v\in \Hom_\2A((Y\boxtimes\11)Q,Q(Z\boxtimes\11))\equiv\Hom_\2E(\oj (Y\boxtimes\11)J,\oj (Z\boxtimes\11)J)$ we have
\begin{equation}(v\bullet u)[k]= (d_k\,\lambda)^{-1} \sum_{i,j\in \Gamma}
   \sum_{\alpha=1}^{N_{ij}^k} d_i d_j \  
\begin{tangle}
\step[1.5]\hstep\object{X_k}\step[1.5]\object{Z}\\
\hh\mobj{{t'}_{ij}^{k\alpha}}\step[1.5]\hcd\step\id\\
\hh\step[1.5]\id\mobj{X_i}\step\id\mobj{X_j}\step\id\\
\hh\step[1.5]\id\step\frabox{v[j]} \\
\hh\step[1.5]\id\step\id\obj{Y}\step\id\\
\hh\step[1.5]\frabox{u[i]}\step\id \\
\hh\step[1.5]\id\step\id\mobj{X_i}\step\id\mobj{X_j}\\
\hh\step[1.5]\id\step\hcu\obj{t_{ij}^{k\alpha}}\\
\step[1.5]\object{X}\step[1.5]\object{X_k}
\end{tangle}
\label{mul}\end{equation}
\elemma
\prf Let $u\in \Hom_\2A((X\boxtimes\11)Q,Q(Y\boxtimes\11))$. Then
\[ v'_j\otimes\id_{Y\boxtimes\11}\mcirc u\mcirc \id_{X\boxtimes\11} \otimes v_i \]
is in
\[\Hom_\2A((X\boxtimes\11)\hat{X}_i,\hat{X}_j(Y\boxtimes\11)) 
  =\Hom_\2C(X X_i, X_jY)\otimes_\7F \Hom_{\2C^\op}(X^\op_i, X^\op_j), \]
which vanishes if $i\ne j$. Thus 
\[ v'_i\otimes\id_{Y\boxtimes\11} \mcirc u\mcirc \id_{X\boxtimes\11}
  \otimes v_i= u[i]\boxtimes\id_{X^\op_i} \]
defines $u[i]\in\Hom_\2C(X X_i,X_iY)$. Conversely, given
$\{ u[i]\in\Hom_\2C(X X_i,X_iY), i\in \Gamma\}$, 
\begin{equation}u=\sum_i v_i\otimes\id_{Y\boxtimes\11} \mcirc
   u[i]\boxtimes\id_{X^\op_i} \mcirc\id_{X\boxtimes\11} \otimes v'_i 
\label{U-ve}\end{equation}
defines a morphism $u\in\Hom_\2A((X\boxtimes\11)Q,Q(Y\boxtimes\11))$. 

Eq. (\ref{mul}) follows easily from (\ref{U-ve}), the definition 
\cite[Proposition 3.8]{mue09} of the $\bullet$-multiplication in $\2E$ and the formula
(\ref{defW}) for $w, w'$.
\qed

\blemma \label{uulemma}
Let $u\in\End_\2B(J(X\boxtimes\11)\oj)$. Then the associated family $\{ u[i]\}$ satisfies 
the braiding fusion equation 
\be
\begin{tangle}
\object{X_i}\step\object{X_j}\step\object{X}\\
\hh\id\step\id\step\id\\
\hh\id\step\frabox{u[j]}\\
\hh\id\step\id\obj{X}\step\id\\
\hh\frabox{u[i]}\step\id\\
\hh\id\step\id\mobj{X_i}\step\id\mobj{X_j}\\
\hh\id\step\frabox{t}\\
\hh\id\step[1.5]\id\\
\object{X}\step[1.5]\object{X_k}
\end{tangle} \ \ = \ \ \ \
\begin{tangle}
\object{X_i}\step\object{X_j}\step\object{X}\\
\hh\id\step\id\step\id\\
\hh\frabox{t}\step\id\\
\hh\step\id\obj{X_k}\step\id\\
\hh\step\frabox{u[k]}\\
\hh\step\id\step\id\\
\hh\step\object{X}\step\object{X_k}
\end{tangle}
\label{bfe}\end{equation}
for all $i,j,k\in \Gamma$ and all $t\in\Hom_\2C(X_k,X_iX_j)$ iff $u$ satisfies
\be
\begin{tangle}
\object{Q}\step\object{Q}\Step\object{X\boxtimes\11}\\
\hh\id\step\id\step\id\\
\hh\id\step\frabox{u}\\
\hh\id\step\id\step\id\\
\hh\frabox{u}\step\id\\
\hh\id\step\hcu\obj{w}\\
\hh\id\step[1.5]\id\\
\step[-.5]\object{X\boxtimes\11}\Step\object{Q}
\end{tangle} \ \ = \ \ \
\begin{tangle}
\object{Q}\step\object{Q}\step[1.7]\object{X\boxtimes\11}\\
\hh\id\step\id\hstep\id\\
\hh\step[-.3]\obj{w}\step[.3]\hcu\hstep\id\\
\hh\hstep\id\step\id\\
\hh\hstep\frabox{u}\\
\hh\hstep\id\step\id\\
\hh\hstep\id\step\mobj{Q}\\
\hstep\object{X\boxtimes\11}
\end{tangle}
\label{uu2}\end{equation}
\elemma
\prf In view of the definition (\ref{defW}) of $w\in\Hom_\2A(Q,Q^2)$ and of 
(\ref{U-ve}), the left hand side of (\ref{uu2}) is seen to equal
\begin{multline*}
 \lambda^{-1/2} \sum_{i,j,k}\sqrt{\frac{d_id_j}{d_k}} v_i\otimes
  v_j\otimes\id_{X\boxtimes\11} 
   \mcirc \id_{\hat{X}_i}\otimes u[j]\boxtimes\id_{X^\op_j}\mcirc
    u[i]\boxtimes\id_{X^\op_i}\otimes\id_{\hat{X}_j} \\
 \mcirc\id_{X\boxtimes\11}\otimes\hat{t}^k_{ij}\mcirc
  \id_{X\boxtimes\11}\otimes v'_k, 
\end{multline*}
whereas the right hand side equals
\[ \lambda^{-1/2}\sum_{i,j,k}\sqrt{\frac{d_id_j}{d_k}} v_i\otimes
   v_j\otimes\id_{X\boxtimes \11}\mcirc \hat{t}^k_{ij}\otimes\id_{X\boxtimes\11}\mcirc 
   u[k]\boxtimes \id_{X^\op_k}\mcirc\id_{X\boxtimes\11}\otimes
   v'_k. \]
In view of the orthogonality relation satisfied by the $v's$, these two expressions
are equal iff
\begin{multline*} \lefteqn{\id_{\hat{X}_i}\otimes u[j]\boxtimes\id_{X^\op_j}\mcirc
    u[i]\boxtimes\id_{X^\op_i}\otimes\id_{\hat{X}_i}\mcirc
   \id_{X\boxtimes\11}\otimes\hat{t}^k_{ij}} \\
   =\hat{t}^k_{ij}\otimes\id_{X\boxtimes\11}\mcirc u[k]\boxtimes
   \id_{X^\op_k} \quad\quad\forall i,j,k\in \Gamma. 
\end{multline*}
Inserting 
$\hat{t}^k_{ij}=\sum_\alpha t^{k\alpha}_{ij}\boxtimes {t'}^{k\alpha}_{ij}$,
this becomes 
\bea \sum_{\alpha=1}^{N_{ij}^k} \left( \id_{X_i}\otimes u[j]\mcirc
   u[i]\otimes\id_{X_j}\mcirc \id_X\otimes t_{ij}^{k\alpha}\right)
  &\boxtimes& {t'}^{k\alpha}_{ij} \nn\\ = \ \ 
  \sum_{\alpha=1}^{N_{ij}^k} \left( t_{ij}^{k\alpha}\otimes\id_X\mcirc u[k] \right) 
  &\boxtimes& {t'}^{k\alpha}_{ij}. \label{cond1}\eea
Multiplying from the right with 
$\id_{XX_k}\boxtimes {t}^{k\alpha}_{ij}$, we arrive at the condition
(\ref{bfe}). Conversely, $\boxtimes$-tensoring (\ref{bfe}) with 
${t'}^{k\alpha}_{ij}$ and summing over $\alpha$ we obtain (\ref{cond1}). 
\qed

\bprop \label{funcf}
There is a faithful functor $F: \2Z_1(\2C)\rightarrow\2B$.
\eprop
\prf Let $(X,e_X)\in\2Z_1(\2C)$. By Lemma \ref{11mor} the half braiding 
$\{ e_X(Z), Z\in\2C\}$ provides us with an element $p^0_X$ in 
$\End_\2B(\oj (X\boxtimes\11)J)\equiv\Hom_\2A((X\boxtimes\11)Q,Q(X\boxtimes\11))$.
Since $e_X(\cdot)$ satisfies the braiding fusion relation (\ref{bfe}), $p^0_X$ 
satisfies (\ref{uu2}). Now, multiplying (\ref{uu2}) from the left with 
$w'\otimes\id_{X\boxtimes\11}$ and using (\ref{W4}) we obtain 
\begin{equation}w'\otimes\id_{X\boxtimes\11}\mcirc\id_Q\otimes p^0_X\mcirc p^0_X\otimes 
  \id_Q\mcirc\id_{X\boxtimes\11}\otimes w = \lambda p^0_X, \label{idemp}\end{equation}
which is just $p^0_X\bullet p^0_X=\lambda p^0_X$ in $\End_\2E(\oj (X\boxtimes\11)J)$. 
Thus with $p_X=\lambda^{-1}p^0_X$, 
\begin{equation}F((X,e_X)):=(\oj (X\boxtimes\11)J, p_X) \label{deff}\end{equation}
is an object in $\2B$, which defines the functor $F$ on the objects. We will mostly write
$F(X,e_X)$ instead of $F((X,e_X))$. Let $(X,e_X)\in\2Z_1(\2C)$ with the above idempotent 
$p_X\in\Hom_\2A((X\otimes\11)Q,Q(X\otimes\11))\equiv\End_\2E(\oj (X\boxtimes\11)J)$
and similarly $(Y,e_Y),\ p_Y$. Consider now
$s\in\Hom_{\2Z_1(\2C)}((X,e_X),(Y,e_Y))\subset\Hom_\2C(X,Y)$. Then condition (\ref{dmor}) 
implies 
\begin{equation}\id_Q\otimes (s\boxtimes\id_{\11}) \circ p_X =
    p_Y \circ (s\boxtimes\id_{\11})\otimes\id_Q. \label{emor}\end{equation}
The element of $u\in\Hom_\2A((X\boxtimes\11)Q,Q(Y\boxtimes\11))\equiv\Hom_\2E(\oj
(X\boxtimes\11)J,\oj (Y\boxtimes\11)J)$ defined by (\ref{emor}) clearly satisfies
$p_Y\bullet u\bullet p_X=u$ and is therefore a morphism in 
$\Hom_\2E((\oj (X\boxtimes\11)J,p_X),(\oj (Y\boxtimes\11)J,p_Y)$. That the map 
$s\mapsto u$ is faithful follows from the first term in (\ref{emor}) and the fact that the
$e_X(X_i),\ i\in \Gamma$, and thus $p_X$ (as a morphism in $\2A$) are invertible. This
defines $F$ on the morphisms, and $F$ is faithful. The simple argument proving that $F$
respects the composition of morphisms is left to the reader. 
\qed

\bprop The functor $F$ is full. \label{fullf}\eprop
\prf We must show that every morphism in $\Hom_\2B(F(X,e_X),F(Y,e_Y))$, where
$(X,e_X),$ $(Y,e_Y)\in\2Z_1(\2C)$, is of the form $F(s)$ with 
$s\in\Hom_{\2Z_1(\2C)}((X,e_X),(Y,e_Y))$.
Now, the morphisms in $\Hom_\2B((\oj (X\boxtimes\11)J,p_X),(\oj (Y\boxtimes\11)J,p_Y))$
are those elements $s$ in $\Hom_\2B(\oj (X\boxtimes\11)J,\oj (Y\boxtimes\11)J)$ which
satisfy $s=p_Y\bullet s\bullet p_X$. $p_X, p_Y$ being idempotents, every such $s$
obviously is of the form $s=p_Y\bullet t\bullet p_X$ for some
$t\in\Hom_\2B(\oj (X\boxtimes\11)J,\oj (Y\boxtimes\11)J)$. By definition of $\2B$ and
by Lemma \ref{11mor}, $s$ and $t$ are represented by elements $\{s[i]\},\ \{t[i]\}$ of
$\bigoplus_{i\in \Gamma} \Hom_\2C(XX_i,X_iY)$.  Given arbitrary $t$ and setting 
$s=p_Y\bullet t\bullet p_X$ we will show that 
$s[0]\in\Hom_\2C(X,Y)$ is in fact in $\Hom_{\2Z_1(\2C)}((X,e_X),(Y,e_Y))$ and that
\[ s[m]= \id_{X_m}\otimes s[0] \mcirc e_X(X_m) \quad\forall m\in \Gamma, \]
which is equivalent to 
\[ s= \id_Q\otimes (s[0]\boxtimes\id_{\11}) \circ p_X = F(s[0]).
\]
Starting from the explicit statement of $s=p_Y\bullet t\bullet p_X$ we compute:
\[ d_m\lambda^4\cdot s[m]=\sum_{i,j,k,l\in \Gamma}\sum_{\alpha,\beta} d_id_jd_k \ \ 
\begin{tangle}
\hstep\object{X_m}\step[2.5]\object{Y}\\
\hh\step[-1]\mobj{{t'}_{kl}^{m\alpha}}\step\hcd\Step\id\\
\id\step\hdcd\step\id\\
\id\step\id\step\hxx\\
\hh\id\step\frabox{t[i]}\step\id \\
\hh\step[-1]\mobj{X_k}\step\id\step\id\step\id\mobj{X_i}\step\id\mobj{X_j}\\
\hxx\step\hddcu\step[-0.5]\mobj{X_l}\\
\hh\id\step\hcu\obj{t_{kl}^{m\alpha}}\\
\object{X}\step[1.5]\object{X_m}
\end{tangle}
\]
\[ =  \sum_{{i,j,k,l\in \Gamma}\atop{\alpha,\beta}} d_id_jd_k\frac{d_m}{d_k}\frac{d_l}{d_j}
\quad\quad
\begin{tangle}
\object{X_m}\step[2.5]\object{\ol{\ve}(\ol{X_l})}\step[2.5]\object{Y}\\
\id\step\Coev\Step\mobj{X_l}\step\id\\
\hh\id\step\id\step\hcoev\step\id\step\id\\
\id\step\id\step\id\step\hcu\step[-.4]\mobj{s_{\ol{\imath}l}^{j\beta}}\step[.9]\ddh\\
\id\step\id\step\id\mobj{X_i}\step[1.5]\hxx\\
\step[-.7]\mobj{s_{m\ol{l}}^{k\alpha}}\step[.7]\hcu\step\id\hstep\hdd\step\id\\
\hh\hstep\id\step[1.5]\frabox{t[i]}\step[1.5]\id\mobj{X_j} \\
\step[-.3]\mobj{X_k}\step[.3]\dh\step\id\step\id\mobj{X_i}\step\hcd\mobj{{s'}_{\ol{\imath}l}^{j\beta}}\\
\step\hxx\step\hev\step[-1.2]\mobj{\ve(X_i)}\step[2.2]\id\mobj{X_l}\\
\step\id\step\hdcd\mobj{{s'}_{m\ol{l}}^{k\alpha}}\Step\id\\
\step\id\step\id\step\ev\mobj{\ve(\ol{X_l})}\\
\step\object{X}\step\object{X_m}
\end{tangle} 
= \ \sum_{{i,j,k,l\in \Gamma}\atop{\alpha,\beta}} d_id_md_l \quad\quad
\begin{tangle}
\object{X_m}\step[2.5]\object{\ol{\ve}(\ol{X_l})}\step[2.5]\object{Y}\\
\id\step\Coev\step[3]\id\\
\id\step\id\step\hcoev\step\hxx\\
\id\step\id\step\id\step\hxx\step\id\mobj{X_l}\\
\hh\step[-1.2]\mmobj{s_{m\ol{l}}^{k\alpha}}\step[1.2]\hcu\step\id\step\id\step\hcu\mmobj{s_{\ol{\imath}l}^{j\beta}}\\
\hh\step[-0.5]\mobj{X_k}\step\id\step[1.5]\frabox{t[i]}\step[1.5]\id\mobj{X_j} \\
\hh\step[-1.2]\mmobj{{s'}_{m\ol{l}}^{k\alpha}}\step[1.2]\hcd\step\id\step\id\mobj{X_i}\step\hcd\mmobj{{s'}_{\ol{\imath}l}^{j\beta}}\\
\id\step\hxx\step\hev\step\id\mobj{X_l}\\
\hxx\step\Ev\step[1.5]\mobj{\ve(\ol{X_l})}\\
\object{X}\step\object{X_m}
\end{tangle}
\]
\[ = d_m \sum_{i,l\in \Gamma} d_id_l \ 
\begin{tangle}
\object{X_m}\step[2.5]\object{\ol{\ve}(\ol{X_l})}\step[2.5]\object{Y}\\
\id\step\Coev\step[3]\id\\
\id\step\id\step\hcoev\step\hxx\\
\id\step\id\step\id\step\hxx\step\id\\
\hh\id\step\id\step\frabox{t[i]}\step\id\step\id\mobj{X_l} \\
\hh\id\step\id\step\id\step\id\mobj{X_i}\step\id\step\id\\
\id\step\hxx\step\hev\step[-1]\mobj{\ve(X_i)}\Step\id\\
\hxx\step\Ev\step[1.5]\mobj{\ve(\ol{X_l})}\\
\object{X}\step\object{X_m}
\end{tangle}
\quad = \ \ d_m  \ \ 
\begin{tangle}
\object{X_m}\Step\object{Y}\\
\id\Step\id\\
\hh\id\step[1.5]\frabox{E_{X,Y}(u)}\\
\xx\\
\hh\id\Step\id\\
\object{X}\Step\object{X_m}
\end{tangle}
\]
where $u\in\Hom_\2C(X,Y)$ does not depend on $m$. (On the very left, $p_X, p_Y$ and each
of the $\bullet$-operations contribute one factor $\lambda$. Furthermore,
$s_{m\ol{l}}^{k\alpha}, \alpha=1,\ldots,N_{m\ol{l}}^k$ is a basis in 
$\Hom(X_k,X_m\ol{X_l})$ with dual basis $s'$. We do not use $t_{m\ol{l}}^{k\alpha}$ since
we cannot assume $\ol{X_l}=X_{\ol{l}}$ without losing generality.)
For $m=0$ we have $X_m=\11$ and thus
$s[0]=\lambda^{-4}E_{X,Y}(u)$, thus $s[0]\in\Hom_{\2Z_1(\2C)}((X,e_X), (Y,e_Y))$.
Plugging this into the above equation for $m\ne 0$ we obtain
$s[m]=\id_{X_m}\otimes s[0]\mcirc e_X(X_m)$ and therefore $s=F(\lambda s[0])$. 
We conclude that the functor $F$ is full. \qed

\bprop The functor $F$ is strong monoidal. \eprop
\prf First we observe that $F(\11_{\2Z_1(\2C)})=F(\11_\2C,\id_X)=(\oj J,\lambda^{-1}\id_Q)$,
which follows from (\ref{deff}) by putting $X=\11$ and $e_X(X_i)=\id_{X_i}\ \forall i$.
Comparing with \cite[Theorem 3.12]{mue09} we see $F(\11_{\2Z_1(\2C)})=\11_\6B$.

Now we have to show that $F(X,e_X)F(Y,e_Y)$ and $F((X,e_X)(Y,e_Y))$ are naturally
isomorphic. We compute 
\bean F(X,e_X)F(Y,e_Y) &=& (\oj (X\boxtimes\11)J,p_X)(\oj (Y\boxtimes\11)J,p_Y) \\
   &=&    (\oj (X\boxtimes\11)Q(Y\boxtimes\11)J,u_1(X,Y)), \\
   F((X,e_X)(Y,e_Y)) &=& F(XY,e_{XY}) = (\oj (XY\boxtimes\11)J,u_2(X,Y)), \eean
where
\[ u_1(X,Y) = \ \ 
\begin{tangle}
\step[1.2]\object{X\boxtimes\11}\step[3]\object{Y\boxtimes\11}\\
\hh\obj{Q}\step\id\step[1.5]\obj{Q}\step[1.5]\id\\
\hh\id\step\id\step\hcd\step\id\\
\hh\frabox{p_X}\step\id\step\frabox{p_Y}\\
\hh\id\step\hcu\step\id\step\id\\
\hh\id\step[1.5]\obj{Q}\step[1.5]\id\step\obj{Q}\\
\step[.2]\object{X\boxtimes\11}\step[3]\object{Y\boxtimes\11}
\end{tangle} \quad\quad\quad\quad\quad
u_2(X,Y) = \lambda \quad
\begin{tangle} 
\step[1.2]\object{X\boxtimes\11}\\
\hh\obj{Q}\step\id\Step\object{Y\boxtimes\11}\\
\hh\id\step\id\step\id\\
\hh\frabox{p_X}\step\id\\
\hh\id\step\id\mobj{Q}\step\id\\
\hh\id\step\frabox{p_Y}\\
\hh\id\step\id\step\id\\
\hh\step[-0.3]\object{X\boxtimes\11}\step[1.3]\id\step\obj{Q}\\
\step[1.2]\object{Y\boxtimes\11}
\end{tangle}
\]
Obviously, $F$ is not strict. But with 
$s(X,Y)\in\Hom_\2E((X\boxtimes\11)(Y\boxtimes\11), (X\boxtimes\11)Q(Y\boxtimes\11))$ and
$t(X,Y)\in\Hom_\2E((X\boxtimes\11)Q(Y\boxtimes\11),(X\boxtimes\11)(Y\boxtimes\11))$
defined by 
\[ s(X,Y)= \quad
\begin{tangle}
\step[1.2]\object{X\boxtimes\11}\step[3]\object{Y\boxtimes\11}\\
\hh\obj{Q}\step\id\step[1.5]\obj{Q}\step[1.5]\id\\
\hh\id\step\id\step\hcd\step\id\\
\hh\frabox{p_X}\step\id\step\frabox{p_Y}\\
\hh\id\step\hcu\step\id\step\id\\
\id\step[1.5]\counit\step[1.5]\id\step\id\\
\hh\id\step[3]\id\step\mobj{Q}\\
\step[.2]\object{X\boxtimes\11}\step[3]\object{Y\boxtimes\11}
\end{tangle} \quad\quad\quad\quad\quad
 t(X,Y)= \lambda\quad
\begin{tangle}
\step[1.2]\object{X\boxtimes\11}\step[3]\object{Y\boxtimes\11}\\
\obj{Q}\step\id\step[1.5]\unit\step[1.5]\id\\
\hh\id\step\id\step\hcd\step\id\\
\hh\frabox{p_X}\step\id\step\frabox{p_Y}\\
\hh\id\step\hcu\step\id\step\id\\
\hh\id\step[1.5]\obj{Q}\step[1.5]\id\step\obj{Q}\\
\step[.2]\object{X\boxtimes\11}\step[3]\object{Y\boxtimes\11}
\end{tangle}
\]
we compute
\[ t(X,Y)\bullet s(X,Y)= \lambda\quad
\begin{tangle}
\hstep\object{Q}\step[1.7]\object{X\boxtimes\11}\step[3]\object{Y\boxtimes\11}\\
\hcd\step\id\step[1.5]\unit\step[1.5]\id\\
\hh\id\step\id\step\id\step\hcd\step\id\\
\hh\id\step\frabox{p_X}\step\id\step\frabox{p_Y}\\
\hh\id\step\id\step\hcu\step\id\step\id\\
\hh\id\step\id\step\hcd\step\id\step\id\\
\hh\frabox{p_X}\step\id\step\frabox{p_Y}\step\id\\
\hh\id\step\hcu\step\id\step\id\step\id\\
\id\step[1.5]\counit\step[1.5]\id\step\hcu\\
\step[.2]\object{X\boxtimes\11}\step[2.7]\object{Y\boxtimes\11}\step[1.5]\object{Q}
\end{tangle} \ \ = \lambda\quad
\begin{tangle}
\hstep\object{Q}\step[1.7]\object{X\boxtimes\11}\step[5]\object{Y\boxtimes\11}\\
\hstep\id\step[1.5]\id\step[1.5]\coev\step[1.5]\id\\
\hh\hcd\step\id\step[1.5]\id\step[1.5]\hcd\step\id\\
\hh\id\step\frabox{p_X}\step[1.5]\id\step[1.5]\id\step\frabox{p_Y}\\
\hh\id\step\id\step\id\step[1.5]\id\step[1.5]\id\step\id\step\id\\
\hh\frabox{p_X}\step\id\step[1.5]\id\step[1.5]\frabox{p_Y}\step\id\\
\hh\id\step\hcu\step[1.5]\id\step[1.5]\id\step\hcu\\
\id\step[1.5]\ev\step[1.5]\id\step[1.5]\id\\
\step[.2]\object{X\boxtimes\11}\step[5]\object{Y\boxtimes\11}\step[1.5]\object{Q}
\end{tangle}
\]
Here the second equality follows (verify!) by repeated use of the equations
(\ref{W1}-\ref{W3}). Using $p_X\bullet p_X=p_X, p_Y\bullet p_Y=p_Y$ and the duality
equation for $Q$ we obtain $t(X,Y)\bullet s(X,Y)=u_2(X,Y)$. Now,
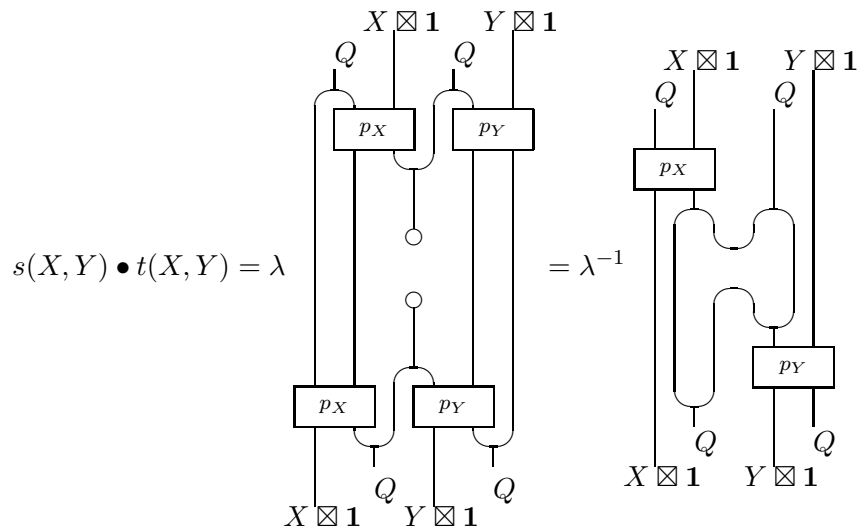
\[ s(X,Y)\bullet t(X,Y)= \lambda\quad
\begin{tangle}
\step\step[1.2]\object{X\boxtimes\11}\step[3]\object{Y\boxtimes\11}\\
\hh\hstep\obj{Q}\step[1.5]\id\step[1.5]\obj{Q}\step[1.5]\id\\
\hh\hcd\step\id\step\hcd\step\id\\
\hh\id\step\frabox{p_X}\step\id\step\frabox{p_Y}\\
\hh\id\step\id\step\hcu\step\id\step\id\\
\id\step\id\step[1.5]\counit\step[1.5]\id\step\id\\
\id\step\id\step[1.5]\unit\step[1.5]\id\step\id\\
\hh\id\step\id\step\hcd\step\id\step\id\\
\hh\frabox{p_X}\step\id\step\frabox{p_Y}\step\id\\
\hh\id\step\hcu\step\id\step\hcu\\
\hh\id\step[1.5]\obj{Q}\step[1.5]\id\step[1.5]\obj{Q}\\
\step[.2]\object{X\boxtimes\11}\step[3]\object{Y\boxtimes\11}
\end{tangle} \ \ = \lambda^{-1}\quad
\begin{tangle}
\step[1.2]\object{X\boxtimes\11}\step[3]\object{Y\boxtimes\11}\\
\hh\obj{Q}\step\id\Step\obj{Q}\step\id\\
\hh\id\step\id\Step\id\step\id \\
\hh\frabox{p_X}\Step\id\step\id\\
\hh\id\hstep\hcd\step\hcd\hstep\id\\
\hh\id\hstep\id\step\hev\step\id\hstep\id\\
\hh\id\hstep\id\step\hcoev\step\id\hstep\id\\
\hh\id\hstep\id\step\id\step\hcu\hstep\id\\
\hh\id\hstep\id\step\id\step[1.5]\frabox{p_Y}\\
\hh\id\hstep\hcu\step[1.5]\id\step\id\\
\hh\id\step\mobj{Q}\Step\id\step\mobj{Q}\\
\step[.2]\object{X\boxtimes\11}\step[3]\object{Y\boxtimes\11}
\end{tangle}
\]
Here we used that $p_X, p_Y$ come from half-braidings, implying that we have (\ref{uu2})
and its dual version by Lemma \ref{uulemma}. (Take into account two factors of $\lambda$
which come from the normalization of $p_{X/Y}$.) It is easy to see that the last
expression equals $u_1(X,Y)$. 

It remains to verify that the functor $F$ is coherent in the sense of
\cite[XI.2]{cwm}. The computations present no difficulties and are simplified by the fact
that the categories $\2Z_1(\2C)$ and $\2B$ are strict except for the unit in
$\2B=\End_\2E(\6B)$. We refrain from spelling out the details. \qed


\subsection{$F$ is essentially surjective}\label{esssurj}
In order to conclude that $F$ establishes an equivalence 
$\2B\stackrel{\otimes}{\cong}\2Z_1(\2C)$ of tensor categories it remains to prove that $F$
is essentially surjective, viz. that for every object $Y$ of $\2B$ there is
$(X,e_X)\in\2Z_1(\2C)$ such that $F(X,e_X)\cong Y$. We begin with a result due to Izumi
\cite{iz2}. 

\blemma Let $Y\in\2C$ be simple. Then the 1-morphisms 
$(Y\boxtimes\11)J,\, (\11\boxtimes Y^\op)J: \ \6B\rarr\6A$ and 
$\oj(Y\boxtimes\11),\, \oj(\11\boxtimes Y^\op): \6A\rarr\6B$ are
simple. Furthermore, 
\bean 
  (Y\boxtimes\11)J &\cong & (\11\boxtimes\ol{Y}^\op)J, \\
  \oj(Y\boxtimes\11) &\cong & \oj(\11\boxtimes\ol{Y}^\op). 
\eean
\label{iz1}\elemma
\prf Let $Y,Z\in\2C$. By duality be have the isomorphism
\[ \Hom_\2E((Y\boxtimes\11)J, (Z\boxtimes\11)J)\cong
  \Hom_\2E((Y\boxtimes\11)J\oj, Z\boxtimes\11)=
  \Hom_\2A((Y\boxtimes\11)Q, Z\boxtimes\11) \]
of vector spaces. In view of $Q\cong\bigoplus_i X_i\boxtimes X^\op_i$ this
implies 
\[ \Hom_\2E((Y\boxtimes\11)J, (Z\boxtimes\11)J)\cong \Hom_\2C(Y,Z). \]
In particular, if $Y\in\2C$ is simple then 
$(Y\boxtimes\11)J\in\Hom_\2E(\6B,\6A)$ is simple, and so is
$(\11\boxtimes Y^\op)J$ by a similar argument. Furthermore,
\[ \Hom_\2E((Y\boxtimes\11)J, (\11\boxtimes\ol{Y}^\op)J)
  \cong \Hom_\2E(Y\boxtimes Y^\op, J\oj) = \Hom_\2A(Y\boxtimes Y^\op, Q). \]
Now, $Y\boxtimes Y^\op$ is simple and contained in $Q$ with
multiplicity one, thus these spaces are one dimensional and
\[ (Y\boxtimes\11)J \cong (\11\boxtimes\ol{Y}^\op)J. \]
Similar arguments apply to the $\6A-\6B$-morphisms. \qed

\bcoro \label{cor}
Let $X, Y\in\2C$. Then there is $Z\in\2C$ such that
\[ \oj(X\boxtimes Y^\op)J\cong \oj(Z\boxtimes\11)J \]
and such that the isomorphisms 
\bean e &\in& \Hom_\2E(\oj(X\boxtimes Y^\op)J, \oj(Z\boxtimes\11)J)
  \equiv\Hom_\2A((X\boxtimes Y^\op)Q, Q(Z\boxtimes\11)), \\
  f &\in& \Hom_\2E(\oj(Z\boxtimes\11)J, \oj(X\boxtimes Y^\op)J)\equiv
  \Hom_\2A((Z\boxtimes\11)Q, Q(X\boxtimes Y^\op)) \eean
can be chosen such that
\[ e=v\otimes \tilde{e},\quad f=v\otimes \tilde{f} \]
with $\tilde{e}\in\Hom_\2A((X\boxtimes Y^\op)Q, Z\boxtimes\11)$ and 
$\tilde{f}\in\Hom_\2A((Z\boxtimes\11)Q, X\boxtimes Y^\op)$. (Alternatively, one can 
find morphisms of the form $e=\tilde{e}\otimes v', f=\tilde{f}\otimes v'$.)
\ecoro
\prf Using the lemma we compute 
\[ \oj(X\boxtimes Y^\op)J=\oj(X\boxtimes\11)(\11\boxtimes Y^\op)J\cong
  \oj(X\boxtimes\11)(\ol{Y}\boxtimes\11)J=\oj(X\ol{Y}\boxtimes\11)J. \]
We put $Z=X\ol{Y}$ and denote by $\hat{e}$ the isomorphism 
$(\11\boxtimes Y^\op)J\rarr(\ol{Y}\boxtimes\11)J$ provided by the preceding lemma.
Now the claim follows with $\tilde{e}=\id_{X\boxtimes\11}\times\hat{e}$ if we keep in mind
that tensoring $\tilde{e}$ with $\id_{\oj}$ (in $\2E$) amounts to tensoring with $v$ in
$\2A$, as follows from the definition of $\2E$. $\tilde{f}$ is defined
similarly. Alternatively, using the isomorphism 
$\oj(X\boxtimes\11)\cong\oj(\11\boxtimes\ol{X}^\op)$ one obtains a solution with
$e=\tilde{e}\otimes v'$ etc.
\qed

The lemma implies that every object of $\2B$ is isomorphic to one of the form
$(\oj(X\otimes\11)J, p_X)$. This looks quite promising since also $F(X,e_X)$ has this
form. In fact, by Lemma \ref{11mor} and Lemma \ref{halfbr} we obtain a family of morphisms 
$\{ e_X(Y): XY\rarr YX,\ Y\in\2C\}$ natural w.r.t.\ $Y$. Yet, in order to conclude that
this is a half braiding (and therefore $(\oj(X\otimes\11)J, p)=F(X, e_X)$) we need that
$p$ satisfies (\ref{uu2}) and $p[0]=\id_X$. Not every object of $\2B$ satisfies these
conditions as is exemplified, e.g., the object
$\oj(X\boxtimes\11)J=(\oj(X\boxtimes\11)J,p)$ where 
\[ p=\id_{\oj(X\boxtimes\11)J}=v\otimes\id_{X\boxtimes\11}\otimes v'
  \in\Hom_\2A((X\boxtimes\11)Q,Q(X\boxtimes\11)). \]
One easily verifies that $p$ does not satisfy (\ref{uu2}). In view of
$p[i]=\delta_{i0}\id_X$ it is also clear that the corresponding $e_X(Y)$ fails to be
invertible for all $Y$.

The following result on the 2-category $\2E$ is quite general in that it does not rely on 
$\2A=\2C\boxtimes\2C^\op$.
\blemma \label{qq-lemma}
Let $X\in\2A$ and $p=p\bullet p\in\End_\2E(\oj XJ)$. Then there is $Y\in\2A$, 
$q=q\bullet q\in\End_\2E(\oj YJ)$ such that $(\oj YJ, q)\cong(\oj XJ, p)$ and in addition
\begin{equation} \label{qq}
\begin{tangle}
\object{Q}\step\object{Q}\step[.7]\object{Y}\\
\hh\id\step\id\hstep\id\\
\hcu\hstep\id\\
\hh\hstep\frabox{q}\\
\hstep\id\step\id\\
\hstep\object{Y}\step\object{Q}
\end{tangle}
\quad=\lambda\quad
\begin{tangle}
\object{Q}\step\object{Q}\step\object{Y}\\
\hh\id\step\id\step\id\\
\hh\id\step\frabox{q}\\
\hh\id\step\id\step\id\\
\hh\frabox{q}\step\id\\
\id\step\hcu\\
\object{Y}\step[1.5]\object{Q}
\end{tangle} 
\end{equation}
\elemma

\brem The condition (\ref{qq}) implies $q\bullet q=q$ as is seen by multiplication with
$w'\otimes\id_Y$ from the left.  
\erem
\prf Using the (non-strict) unit $\6B-\6B$ morphism $\11_\6B=(\oj J, \lambda^{-1}\id_Q)$
we put $(Y,q)=\11_\6B(X,p)=(QX, \lambda^{-1}\id_Q\times p)$. The isomorphism 
$(\oj YJ, q)\cong(\oj XJ, p)$ was proven in \cite[Theorem 3.12]{mue09}. We claim that $q$
satisfies (\ref{qq}). In terms of $(Y,q)$ (and keeping in mind that $Y=QX$ !) the left
hand side of (\ref{qq}) is given by 
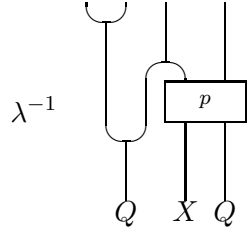
\[ \lambda^{-1}\quad
\begin{tangle}
\hh\hcu\step\id\step[1.5]\id\\
\hh\hstep\id\step\hcd\step\id\\
\hh\hstep\id\step\id\step\frabox{p}\\
\hstep\hcu\step\id\step\id\\
\hstep\hstep\object{Q}\step[1.5]\object{X}\step\object{Q}
\end{tangle}
\]
For the right hand side we compute
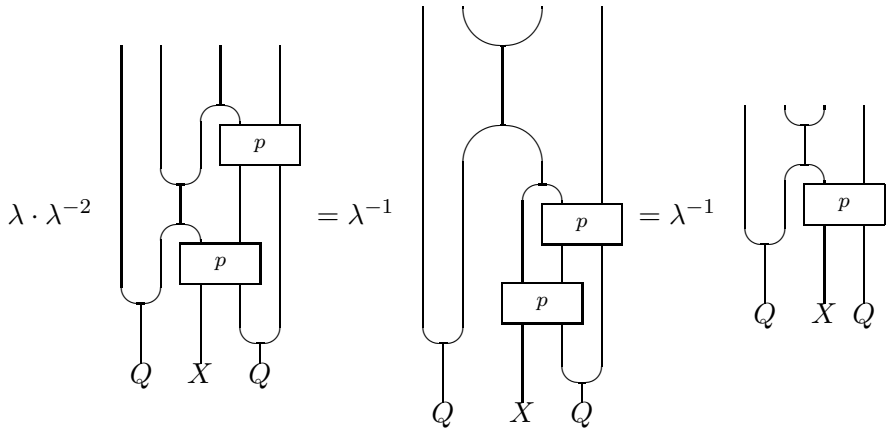
\[ \lambda\cdot\lambda^{-2}\quad
\begin{tangle}
\id\hstep\hstep\id\step\hcd\step\id\\
\hh\id\hstep\hstep\id\step\id\step\frabox{p}\\
\hh\id\hstep\hstep\hcu\step\id\step\id\\
\hh\id\step\hcd\step\id\step\id\\
\hh\id\step\id\step\frabox{p}\step\id\\
\hh\hcu\step\id\step\id\step\id\\
\hh\hstep\id\step[1.5]\id\step\hcu\\
\hstep\object{Q}\step[1.5]\object{X}\step[1.5]\object{Q}
\end{tangle}
\quad=\lambda^{-1}\quad
\begin{tangle}
\id\step\cu\step[1.5]\id\\
\id\step\cd\step[1.5]\id\\
\hh\id\step\id\step[1.5]\hcd\step\id\\
\hh\id\step\id\step[1.5]\id\step\frabox{p}\\
\hh\id\step\id\step[1.5]\id\step\id\step\id\\
\hh\id\step\id\step[1.5]\frabox{p}\step\id\\
\hh\hcu\step[1.5]\id\step\id\step\id\\
\hh\hstep\id\Step\id\step\hcu\\
\hstep\object{Q}\Step\object{X}\step[1.5]\object{Q}
\end{tangle}
\quad=\lambda^{-1}\quad
\begin{tangle}
\hh\id\step\cu\step\id\\
\hh\id\step\cd\step\id\\
\hh\id\step\id\step\frabox{p}\\
\hcu\step\id\step\id\\
\hstep\object{Q}\step[1.5]\object{X}\step\object{Q}
\end{tangle}
\]
In the last step we have used $p\bullet p=p$. That the result coincides with the left hand
side follows now from a standard computation using the properties of a Frobenius algebra.
\qed

\bprop Every object of $\2B$ is isomorphic to one of the form
$(\oj(Z\boxtimes\11)J, q)$ where $q\in\End_\2B(\oj(Z\boxtimes\11)J)$ satisfies (\ref{qq})
(with $Y=Z\boxtimes\11$).
\eprop
\prf By the preceding lemma every object $(\oj(X\boxtimes Y^\op)J,p)$ of $\2B$ is
isomorphic to one which satisfies (\ref{qq}), which allows us to assume this property in
the rest of the proof. By Corollary \ref{cor} there is $Z\in\2C$ such that 
$\oj(X\boxtimes Y^\op)J\cong\oj(Z\boxtimes\11)J$. Let
$e: \oj(X\boxtimes Y^\op)J\rarr\oj(Z\boxtimes\11)J$, 
$f: \oj(Z\boxtimes\11)J\rarr\oj(X\boxtimes Y^\op)J$ be a pair of mutually inverse
isomorphisms. Then with $q=e\bullet p\bullet f$ we have
$(\oj(Z\boxtimes\11)J, q)\cong(\oj(X\boxtimes Y^\op)J, p)$. If we can show that 
also $q$ satisfies (\ref{qq}) Lemma \ref{uulemma} applies and the claim follows.
Now by Corollary  \ref{cor} $Z, e, f$ can be chosen such that $e=v\otimes\tilde{e},
f=v\otimes\tilde{f}$, where $\tilde{e}, \tilde{f}$ are mutually inverse 2-morphisms
between $(X\otimes Y)J$ and $(Z\otimes\11)J$. Therefore, 
\[ q=e\bullet p\bullet f =\quad
\begin{tangle}
\step\object{Q}\Step\object{Z}\\
\hh\step\id\Step\id\\
\sw1\id\se1\step\id\\
\hh\id\step\id\step\frabox{e}\\
\hh\id\step\id\step\id\mmobj{X}\step\id\\
\hh\id\step\frabox{p}\step\id\\
\hh\id\step\id\mmobj{X}\step\id\step\id\\
\hh\frabox{f}\step\id\step\id\\
\id\step\nw1\id\ne1\\
\hh\id\Step\id\\
\object{Z}\Step\object{Q}
\end{tangle}
\quad=\quad
\begin{tangle}
\step\object{Q}\Step\object{Z}\\
\hh\step\id\Step\id\\
\hh\step\id\step\frabox{\tilde{e}}\\
\hh\step\id\step\id\mmobj{X}\step\id\\
\hh\step\frabox{p}\step\id\\
\hh\step\id\mmobj{X}\step\id\step\id\\
\hh\frabox{\tilde{f}}\step\id\step\id\\
\id\step\nw1\id\ne1\\
\hh\id\Step\id\\
\object{Z}\Step\object{Q}
\end{tangle}
\]
where the four-fold vertices denote triple (co)products.
That $q$ satisfies (\ref{qq}) is now obvious from the respective property of $p$ and 
$\tilde{f}\bullet\tilde{e}=\id$.
\qed

\bprop \label{prop-ess}
The preceding proposition remains true if one adds the requirement that
$q[0]=\lambda^{-1}\id_Z$ (notation of Lemma \ref{11mor}).
\eprop
\prf Let $Z, q$ be as in the preceding proposition. Multiplying (\ref{qq}) with 
$v'\otimes v'\otimes\id_{Z\boxtimes\11}$ and using 
$v'\otimes\id_{Z\boxtimes\11}\circ q=q[0]\otimes v$ we obtain 
$\lambda q[0]^2=q[0]\in\End_\2C(Z)$. 
Let $f: \tilde{Z}\rarr Z,\, g: Z\rarr\tilde{Z}$ be a splitting of the idempotent
$\lambda q[0]$. Then it is easy to verify that with 
\[ \tilde{q}=\id_Q\otimes (g\boxtimes\id_\11) \mcirc q\mcirc (f\boxtimes\id_\11)\otimes\id_Q \]
we have $(\oj(\tilde{Z}\boxtimes\11)J, \tilde{q})\cong(\oj(Z\boxtimes\11)J, q)$. This
$\tilde{q}$ still verifies (\ref{qq}) and, in addition,
$\tilde{q}[0]=\lambda^{-1}\id_{\tilde{Z}}$.
\qed

Now we are ready to state our first main result.

\btheor The tensor categories $\2B$ and $\2Z_1(\2C)$ are equivalent as spherical categories,
thus we have the weak monoidal Morita equivalence (in the sense of \cite{mue09})
$\2Z_1(\2C)\approx\2C\boxtimes\2C^\op$. In particular,
\[ \dim\2Z_1(\2C)=(\dim\2C)^2. \]
\etheor
\prf We have shown that every simple object in $\2B$ isomorphic to the image under $F$ of
a simple object in $\2Z_1(\2C)$. Since $\2Z_1(\2C)$ and $\2B$ are both semisimple (in
particular closed under direct sums and subobjects) we conclude that $F$ is essentially
surjective. Since $F$ is also fully faithful we have an equivalence of categories by
\cite[Theorem IV.4.1]{cwm}. $F$ being monoidal we have an equivalence of monoidal categories
by \cite[I.4.4]{sr}. (This already implies that $\2B$ and $\2Z_1(\2C)$ have the same
dimension, since by \cite[Proposition 2.4]{mue09} the latter are well-defined independently of
the chosen spherical or $*$-structure and, of course, invariant under monoidal equivalence.)

It remains to show that the spherical structures are compatible. As to the conjugation
maps, we have
\bean \ol{F(X,e_X)} &=& \ol{(\oj(X\boxtimes\11)J,p_{e_X})}=(\oj(\ol{X}\boxtimes\11)J,\ol{p_{e_X}}), \\
   F(\ol{(X,e_X)}) &=& F((\ol{X},e_{\ol{X}}))=(\oj(\ol{X}\boxtimes\11)J, p_{e_{\ol{X}}}). \eean
Putting together Proposition \ref{hb-conj} and Lemma \ref{iii-iv}, $e_{\ol{X}}(X_i)$ is as in
Fig.\ \ref{fig1}, where the pair of unlabeled morphisms is any solution of the duality
equations. 
\begin{figure}
\[ e_{\ol{X}}(X_i)=\quad
\begin{tangle}
\step[1.5]\object{\ol{\ve}(\ol{X})}\step[2.5]\object{X_i}\step\object{\ol{X}}\\
\mcoev\step\id\step\id\\
\hh\id\step\hcoev\step\id\step\id\step\id\\
\hh\id\step\id\step\frabox{e_X(X_{\ol{\imath}})}\step\id\step\id\\
\hh\id\step\id\step\id\step\hev\step\id\\
\id\step\id\step\mev\\
\object{\ol{X}}\step\object{X_i}\step[2.5]\object{\ve(X)}
\end{tangle}
\quad\quad\quad\quad
\ol{p_{e_X}}=\quad
\begin{tangle}
\step[1.5]\object{\ol{\ve}(\ol{X}\boxtimes\11)}\step[2.5]\object{Q}\step[1.5]\object{\ol{X}\boxtimes\11}\\
\mcoev\step\id\step\id\\
\hh\id\step[.8]\mobj{r'}\step[.2]\hcoev\step\id\step\id\step\id\\
\hh\id\step\id\step\frabox{p_{e_X}}\step\id\step\id\\
\hh\id\step\id\step\id\step\hev\mobj{r}\step\id\\
\id\step\id\step\mev\\
\step[-.5]\object{\ol{X}\boxtimes\11}\step[1.5]\object{Q}\step[2.5]\object{\ve(X\boxtimes\11)}
\end{tangle}
\]
\caption{$e_{\ol{X}}(X_i)$ and $\ol{p_{e_X}}$}
\label{fig1}
\end{figure}
In view of the definition of $p_X$ in Proposition \ref{funcf} and of $\ol{p_{e_X}}$ in
\cite[Theorem 5.14]{mue09}, cf.\ Fig.\ \ref{fig1}, it is clear that
$\ol{p_{e_X}}=p_{e_{\ol{X}}}$ and therefore 
\[ \ol{F(X,e_X)}=F(\ol{(X,e_X)}). \]

Now by Proposition \ref{hb-conj}, the spherical structure of $\2Z_1(\2C)$ is inherited from 
$\2C$, concretely $\ve_{\2Z_1(\2C)}(X,e_X)=\ve_\2C(X)$. Considering how the spherical
structures of $\2E_0$ and $\2E$ arise from that of $\2A$ in \cite[Theorem 5.13]{mue09} it
is essentially obvious that $F: \2Z_1(\2C)\rarr\2B$ is an equivalence of spherical
categories irrespective of the fact that the latter is neither strict monoidal nor strict
spherical. We omit the easy details.
\qed

\brem 1. In the case where $\2C$ is the representation category of a finite dimensional
involutive semisimple and cosemisimple Hopf algebra $H$, $\2Z_1(\2C)$ is equivalent
\cite[Theorem XIII.5.1]{ka} to the representation category of the quantum double $D(H)$
and our result is just the fact $\dim D(H)=(\dim H)^2$.

2. It seems likely that a simpler proof of the theorem can be given using the
interpretation of the tensor category $\2B=\END_\2E(\6B)$ as bimodule category
$Q-\mbox{Mod}-Q$, together with the recent work \cite{ostr}.
\erem


\section{Modularity of the Quantum Double}\label{sec-modular}
\subsection{The `Tube algebra'} \label{ss-tube}
By definition \cite{mue09} of $\2E$, every simple $\6B-\6B$ morphism is contained in
$\oj(X\boxtimes Y^\op)J$ for some simple $X, Y$. In view of Lemma \ref{iz1} every simple
$\6B-\6B$-morphism is  in fact contained in  $\oj (X\boxtimes\11)J$ for some simple
$X\in\2C$ (as well as in $\oj (\11\boxtimes \ol{X}^\op) J$). Defining
\[ \hat{Y}_L = \bigoplus_{i\in \Gamma} X_i\boxtimes \11, \quad\quad
   \hat{Y}_R = \bigoplus_{i\in \Gamma} \11\boxtimes X_i^\op, \]
we conclude that either of $\oj \hat{Y}_LJ,\ \oj \hat{Y}_RJ$ 
contains all simple $\6B-\6B$-morphisms. With
\[ \Xi_L = \End_\2E(\oj \hat{Y}_LJ), \quad\quad  \Xi_R = \End_\2E(\oj \hat{Y}_R J) \]
we thus have a one-to-one correspondence between isomorphism classes of simple 
$\6B-\6B$-morphisms and minimal central idempotents in $\Xi_L$ or, equivalently, in
$\Xi_R$. From now on we will stick to $\Xi_L$. By construction of $\2E$ we have 
$\Xi_L=\Hom_\2A(\hat{Y}_LQ,Q\hat{Y}_L)$ as a vector space. Thus
\bean \Xi_L &\cong& \bigoplus_{i,j,k,l} \, \Hom_\2A(X_iX_j\boxtimes X_j^\op,
   X_kX_l\boxtimes X_k^\op) \\
  &\cong& \bigoplus_{i,j,k,l} \, \Hom_\2C(X_iX_j,X_kX_l) \,\otimes_\7F \,
  \Hom_{\2C^\op}(X_j^\op, X_k^\op)  \\
  &\cong& \bigoplus_{i,j,l} \, \Hom_\2C(X_iX_j,X_jX_l). \eean
We therefore have 
\begin{equation}\End_\2E(\oj \hat{Y}_L J)\equiv\Hom_\2A(\hat{Y}_LQ,Q\hat{Y}_L)\cong
   \bigoplus_{i,j,k} \, \Hom_\2C(X_iX_j,X_jX_k), \label{tube1}\end{equation}
and in complete analogy to the proof of (\ref{mul}) one shows that the multiplication in
$\Xi_L$ is given by
\begin{equation}(v\bullet u)[i,j,k]= (d_j\lambda)^{-1} \sum_{l,m,n\in \Gamma} 
  \sum_{\alpha=1}^{N_{mn}^j} d_md_n \ 
\begin{tangle}
\step[1.5]\hstep\object{X_j}\step[1.5]\object{X_k}\\
\hh\mobj{{t'}_{mn}^{j\alpha}}\step[1.5]\hcd\step\id\\
\hh\step[1.5]\id\mobj{X_m}\step\id\mobj{X_n}\step\id\\
\hh\step[1.5]\id\step\frabox{v[l,n,k]} \\
\hh\step[1.5]\id\step\id\obj{X_l}\step\id\\
\hh\step[1.5]\frabox{u[i,m,l]}\step\id \\
\hh\step[1.5]\id\step\id\obj{X_m}\step\id\obj{X_n}\\
\hh\step[1.5]\id\step\hcu\obj{t_{mn}^{j\alpha}}\\
\step[1.5]\object{X_i}\step[1.5]\object{X_j}
\end{tangle}
\label{tube2}\end{equation}

\brem 1. We observe that up to a different normalization (\ref{tube1}) and (\ref{tube2})
coincide with Ocneanu's definition of the `tube algebra', cf. \cite{ocn3,ek1,iz2}. 
(The $(ij|X|jk)$ of Izumi corresponds to our $u[i,j,k]\cdot\frac{d_j}{\lambda}$.)
Note, however, that we {\it derive} (\ref{tube1}, \ref{tube2}) from an intrinsic
definition of the algebra $\Xi_L=\End_\2E(\oj \hat{Y}_LJ)$, which makes the
correspondence between the minimal central idempotents of $\Xi_L$ and the isomorphism
classes of simple objects in $\2B$ completely obvious. (Compare this to the laborious
proof in \cite{iz2}.) The above considerations therefore completely clarify the r\^{o}le
of the tube algebra. We suspect that Ocneanu arrived at his definition of the tube algebra
by similar considerations.

2. Note that in the definition of $\Xi_L$ we could replace $\hat{Y}_L$ by
\[ \hat{Y}^N_L = \bigoplus_{i\in \Gamma} N_i\,(X_i\boxtimes \11) \] 
with arbitrary $\{N_i\}\in\7N^\Gamma$. The algebras $\End_\2E(\oj \hat{Y}^N_L J)$, of
which the tube algebra happens to be the smallest, all have the same center, thus are
Morita equivalent. We emphasize that only this common center has an invariant meaning, and
in fact it has a well-known interpretation in terms of TQFTs, see Subsection \ref{TQFT}.
\erem

\blemma \label{adj}
Let $(X,e_X)\in \2Z_1(\2C)$ and $Y\in\2C$. Then there is a isomorphism between the 
vector spaces $\Hom_\2C(X,Y)$ and $\Hom_\2B(F(X,e_X), \oj (Y\boxtimes\11)J)$. 
\elemma
\prf The proof is similar to the one of Proposition \ref{fullf}, but simpler. Let
$s\in\Hom_\2C(X,Y)$ and let 
$\ul{t}\in\Hom_\2B(\oj (X\boxtimes\11)J,\oj (Y\boxtimes\11)J)$ be defined (using Lemma
\ref{11mor}) by $\ul{t}[i]=\delta_{i0}s$. Then the map 
$\pi: s\mapsto\ul{s}=\ul{t}\bullet p_{(X,e_X)}\in\Hom_\2B(F(X,e_X),\oj (Y\boxtimes\11)J)$ 
is injective since 
$(\ul{t}\bullet p_{(X,e_X)})[i]=\lambda^{-1}\id_{X_i}\otimes s\circ e_X(X_i)$, in
particular $(\ul{t}\bullet p_{(X,e_X)})[0]=\lambda^{-1}s$. Let, conversely, 
$\ul{s}\in\Hom_\2B(F(X,e_X),\oj (Y\boxtimes\11)J)$, i.e.\
$\ul{s}=\ul{s}\bullet p_{(X,e_X)}\in\Hom_\2B(\oj (X\boxtimes\11)J,\oj (Y\boxtimes\11)J)$. 
Then 
\[ d_k\lambda^2\cdot \ul{s}[k]= 
\sum_{m,l\in \Gamma}\sum_{\alpha=1}^{N_{ml}^k} \ d_md_l \ \ 
\begin{tangle}
\hstep\object{X_k}\step[1.5]\object{Y}\\
\hh\hstep\id\step[1.5]\id\\
\hh\step[-1.5]\mmobj{t_{m,l}^{'k\alpha}}\step[1.5]\hcd\step\id\\
\hh\step[-1.2]\mobj{X_m}\step[1.2]\id\step\id\step\id\\
\hh\id\step\frabox{\ul{s}[l]} \\
\hxx\step\id\mobj{X_l}\\
\hh\id\step\hcu\mobj{t_{m,l}^{k\alpha}}\\
\hh\id\step[1.5]\id\\
\object{X}\step[1.5]\object{X_k}
\end{tangle}
\ \ = \ \sum_{m,l\in \Gamma}\sum_\alpha \ d_md_l\frac{d_k}{d_m}
\begin{tangle}
\object{X_k}\step[1.5]\object{\ol{\ve}(\ol{X_l})}\step[1.5]\object{Y}\\
\hh\id\step\hcoev\step\id\\
\hh\step[-1.2]\mmobj{s_{k\ol{l}}^{m\alpha}}\step[1.2]\hcu\step\id\step\id\\
\hh\hstep\id\step[1.5]\frabox{\ul{s}[l]} \\
\hstep\id\step\ddh\hstep\id\\
\hstep\hxx\mobj{X_m}\step[1.5]\id\mobj{X_l}\\
\step[-0.5]\hdd\step[-0.5]\mmobj{{s'}_{k\ol{l}}^{m\alpha}}\step\hcd\step\id\\
\hh\id\step\id\step\hev\mmobj{\ve(\ol{X_l})}\\
\object{X}\step\object{X_k}
\end{tangle}
\]
\[ = d_k \ \sum_{m,l\in \Gamma}\sum_\alpha  d_l \ \ 
\begin{tangle}
\object{X_k}\step[3]\object{Y}\\
\hh\id\step\hcoev\step\id\\
\hh\hcu\step\id\step\id\\
\hh\step[-0.7]\mobj{X_m}\step[1.2]\id\step[1.5]\frabox{\ul{s}[l]} \\
\hh\hcd\step\id\step\id\\
\id\step\hxx\step\id\mobj{X_l}\\
\hxx\step\hev\\
\object{X}\step\object{X_k}
\end{tangle}
\ \ = d_k \ \sum_{l\in \Gamma} \ d_l \ \
\begin{tangle}
\object{X_k}\step[1.5]\object{\ol{\ve}(\ol{X_l})}\step[1.5]\object{Y}\\
\hh\id\step\hcoev\step\id\\
\hh\id\step\id\step\id\step\id\\
\hh\id\step\id\step\frabox{\ul{s}[l]} \\
\hh\id\step\id\step\id\step\id\\
\id\step\hxx\step\id\mobj{X_l}\\
\hxx\step\hev\step[-.5]\mmobj{\ve(\ol{X_l})}\\
\object{X}\step\object{X_k}
\end{tangle}
\]
Putting $k=0$ we obtain
\[ \lambda^2\cdot \ul{s}[0]= \sum_l d_l \ \  
\begin{tangle}
\Step\object{Y}\\
\hh\hcoev\step\id\\
\hh\id\step\id\step\id\\
\hh\id\step\frabox{\ul{s}[l]} \\
\hxx\step\id\mobj{X_l}\\
\hh\id\step\hev\\
\object{X}
\end{tangle}
\]
and plugging this back into the preceding equation we obtain
\[ \ul{s}[k]=\id_{X_i}\otimes \ul{s}[0]\circ e_X(X_i). \]
Thus $\ul{s}$ is in the image of $\pi$, which proves that $\pi$ is an isomorphism.
\qed

\bprop Let $(X,e_X)\in \2Z_1(\2C)$ be simple and let $N^X_i=\dim\Hom_\2C(X_i,X)$.
Then the (simple) object $F(X,e_X)\in\2B$ is contained in $\oj (X_i\boxtimes\11)J$
with multiplicity $N^X_i$. Let $\{ p_i^\alpha \}, \{ {p'}_i^\alpha \}$ be bases in 
$\Hom_\2C(X_i,X), \Hom_\2C(X, X_i)$, respectively, normalized by
${p'}_i^{\alpha}\circ p_i^\beta=\delta_{\alpha\beta}\id_{X_i}$. Then
$q_i^\alpha\in\Hom_\2E(F(X,e_X), \oj (X_i\boxtimes\11)J)$,
${q'}_i^{\alpha}\in\Hom_\2E(\oj (X_i\boxtimes\11)J,F(X,e_X))$ defined by 
\bean q_i^\alpha[k] &=& \left(\frac{d(X)}{\lambda^2 d_i}\right)^{1/2} \id_{X_k}\otimes
  {p'}_i^{\alpha} \mcirc e_X(X_k),\\ 
   {q'}_i^{\alpha}[k] &=& \left(\frac{d(X)}{\lambda^2 d_i}\right)^{1/2} e_X(X_k) \mcirc
   p_i^{\alpha}\otimes \id_{X_k}  
\eean
satisfy ${q'}_i^{\alpha}\bullet q_i^\beta= \delta_{\alpha\beta} \, \id_{F(X,e_X)}$. 
The idempotent $z^i_{(X,e_X)}=\sum_{\alpha=1}^{N^X_i}q_i^\alpha\bullet {q'}_i^\alpha$ in 
$\End_\2E(\oj (X_i\boxtimes\11)J)$ corresponding to the isotypic component of $(X,e_X)$
is given by 
\begin{equation}z^i_{(X,e_X)}[k]= \frac{d(X)}{\lambda d_i} \sum_{\alpha=1}^{N^X_i} \id_{X_k}\otimes
  {p'_i}^\alpha \circ e_X(X_k) \circ p_i^\alpha \otimes\id_{X_k}. \end{equation}
\label{dualb}\eprop

\brem The choice of square root of $d(X)$ is immaterial, but it must be the same in the 
equations defining $q_i^\alpha$ and ${q'}_i^\alpha$.
\erem
\prf In view of the preceding lemma all that remains to be verified is the normalization. 
Since $(X,e_X)$ is simple we have 
${q'}_i^{\alpha}\bullet q_i^\beta=c_{\alpha\beta}\,\id_{F(X,e_X)}$. Plugging
$q_i^\alpha[k], {q'}_i^\alpha[k]$ into (\ref{mul}) and comparing with the middle term of
the computation in Lemma \ref{condexp1} (with $Z=X_k$) we see that 
\bean ({q'}_i^{\alpha}\bullet q_i^\beta)[k] &=& \frac{d(X)}{\lambda d_i} \
   \id_{X_k}\otimes E_{X,X}(p_i^{\alpha}\circ{p'}_i^\beta)\mcirc e_X(X_k) \\
 &=& \frac{tr_X\circ E_{X,X}(p_i^{\alpha}\circ{p'}_i^\beta)}{\lambda d_i} \  e_X(X_k), 
\eean
since $E_{X,X}(p_i^{\alpha}\circ{p'}_i^\beta)$ is a scalar multiple of $\id_X$ due to the
simplicity of $(X,e_X)$. Now, by definition of the functor $F$ we have
$\id_{F(X,e_X)}[k]=\lambda^{-1}e_X(X_k)$, thus by comparison we find
$c_{\alpha\beta}=d_i^{-1} tr_X\circ E_{X,X}(p_i^{\alpha}\circ{p'}_i^\beta)$. Computing
\[ tr_X\circ E_{X,X}(p_i^{\alpha}\circ{p'}_i^\beta)=tr_X(p_i^{\alpha}\circ{p'}_i^\beta)=
  tr_{X_i}({p'}_i^{\alpha}\circ p_i^\beta)=d_i   \delta_{\alpha\beta}, \] 
where we used the invariance of the trace under the conditional expectation and cyclic
permutations, we obtain $c_{\alpha\beta}=\delta_{\alpha\beta}$ as claimed.

Now we can compute $z^i_{(X,e_X)}=\sum_\alpha q_i^\alpha\bullet {q'}_i^\alpha$ as follows:
\[ z^i_{(X,e_X)}[k]= \frac{d(X)}{\lambda^2 d_i} \sum_{\alpha=1}^{N^X_i} 
    \sum_{l,m\in \Gamma} \sum_{\beta=1}^{N_{lm}^k} \frac{d_ld_m}{\lambda d_k} \ 
\begin{tangle}
\step[1.5]\hstep\object{X_k}\step[1.5]\object{X_i}\\
\hh\mobj{{t'}_{l,m}^{k,\beta}}\step[1.5]\hcd\step\id\\
\hh\step[1.5]\id\mobj{X_l}\step\id\mobj{X_m}\step\id\\
\step[1.5]\id\step\id\step\O {{p'}_i^\alpha} \\
\step[1.5]\id\step\hxx \\
\hh\step[1.5]\id\step\id\obj{X}\step\id\\
\step[1.5]\hxx\step\id \\
\step[1.5]\O {p_i^\alpha}\step\id\step\id \\
\hh\step[1.5]\id\step\id\obj{X_l}\step\id\obj{X_m}\\
\hh\step[1.5]\id\step\hcu\obj{t_{l,m}^{k,\beta}}\\
\step[1.5]\object{X_i}\step[1.5]\object{X_k}
\end{tangle}
\ \ = \ \ 
\frac{d(X)}{\lambda d_i} \sum_{\alpha=1}^{N^X_i} \ \ 
\begin{tangle}
\hstep\object{X_k}\step\object{X_i}\\
\hh\hstep\id\step\id\\
\hstep\id\step\O {{p'}_i^\alpha} \\
\mobj{X}\hstep\hxx \\
\hstep\O {p_i^\alpha}\step\id \\
\hh\hstep\id\step\id\\
\hstep\object{X_i}\step\object{X_k}\\
\end{tangle}
\]
We have pulled ${t'}_{lm}^{k\beta}$ through the braiding and used (\ref{ddN}).
\qed

\bprop Let $(X,e_X)\in \2Z_1(\2C)$ be simple. The minimal central idempotent $z_{(X,e_X)}$
in $\Xi_L$ corresponding to $F(X,e_X)$ is given by
\begin{equation}z_{(X,e_X)}[i,j,k]= \delta_{ik} \frac{d(X)}{\lambda d_i}
  \sum_{\alpha=1}^{N^X_i} \id_{X_j}\otimes {p'_i}^\alpha \circ 
   e_X(X_j) \circ p_i^\alpha \otimes\id_{X_j}, \label{isotyp}\end{equation}
where the $\{ p_i^\alpha \}, \{ {p'_i}^\alpha \},\ i\in \Gamma, \alpha=1,\ldots,N^X_i$
are bases as in Proposition \ref{dualb}. \eprop
\prf Since $\oj \hat{Y}_LJ$ is a direct sum $\bigoplus_i \oj (X_i\boxtimes\11)J$ we
only need to add up the idempotents in $\End_\2E(\oj (X_i\boxtimes\11)J)$, which we
identified in Proposition \ref{dualb}, inside $\Xi_L=\End_\2E(\oj \hat{Y}_LJ)$. With the
isomorphism (\ref{tube1}) the claimed identity follows. \qed

As a first application of the tube algebra we can give an easy bound on the `size' of the
quantum double:
\bcoro The number $\#\2Z_1(\2C)$ of isomorphism classes of simple objects of $\2Z_1(\2C)$
satisfies 
\[ \#\2Z_1(\2C) \le \sum_{i,j\in \Gamma} \dim\Hom_\2C(X_iX_j,X_jX_i). \]
\ecoro
\prf By the equivalence $\2Z_1(\2C)\stackrel{\otimes}{\simeq}\2B$ and the above
considerations we have $\#\2Z_1(\2C)=\dim Z(\Xi_L)$. Since the center of $\Xi_L$ is spanned
by the $z_{(X,e_X)}$ constructed above and since $z_{(X,e_X)}[i,j,k]=0$ if $i\ne k$ we
have 
\[ Z(\Xi_L)\subset\bigoplus_{i,j}\Hom_\2C(X_iX_j,X_jX_i), \]
which implies the bound.
\qed

\brem If $G$ is a finite abelian group whose order is non-zero in $\7F$ then
$G-\mbox{mod}$ is semisimple, symmetric and all simple objects have dimension one. Thus
the right hand side of the above inequality equals $|G|^2$. In view of
$\2Z_1(G-\mbox{mod})\simeq D(G)-\mbox{mod}$ we have $\#\2Z_1(G-\mbox{mod})=|G|^2$, which
proves that the bound is optimal.
\erem

The next two subsections, which do not pretend much originality, will follow \cite{iz2}
quite closely except for shortcuts in the proofs.


\subsection{Invertibility of the $S$-matrix}\label{modular}
In this subsection we will prove that the $S$-matrix 
\[ S((X,e_X),(Y,e_Y))=\quad
\begin{tangle}
\coev\step\coev\\
\id\Step\hxx\Step\id\\
\id\step[.2]{\scriptsize \mobj{(X,e_X)}}\step[1.8]\hxx\step[-.2]{\scriptsize \mobj{(Y,e_Y)}}\step[2.2]\id\\
\step[-3]\mobj{\ve(\ol{(X,e_X)})}\step[3]\ev\step\ev\mobj{\ve((Y,e_Y))}
\end{tangle}
\]
of $\2Z_1(\2C)$ is invertible, thus $\2Z_1(\2C)$ is modular in the sense of Turaev
\cite{t}. The strategy will be to define a vector space isomorphism $\6S$ of the
subspace 
\[ \Xi_0=\bigoplus_{i,j\in \Gamma}\Hom_\2C(X_iX_j,X_jX_i) \]
of $\Xi_L$ which we have seen to contain the center of $\Xi_L$. We will prove that $\6S$
leaves $Z(\Xi_L)$ stable and that the $S$-matrix of $\2Z_1(\2C)$ is the matrix
representation of $\6S\restr Z(\Xi_L)$ w.r.t.\ the basis 
$\{ d(X)^{-1}z_{(X,e_X)},\ (X,e_X)\ \mbox{simple}\}$.

\blemma 
The application $\6S: \Xi_0\rarr\Xi_0$ defined by
\begin{equation}
\begin{tangle}
\object{X_j}\step\object{X_i}\\
\hh\id\step\id\\
\hh\frabox{s}\\
\hh\id\step\id\\
\object{X_i}\step\object{X_j}
\end{tangle} \quad\ \  \mapsto \quad
\begin{tangle}
\Step\object{X_i}\step\object{X_{\ol{\jmath}}}\\
\hh\hcoev\step\id\step\id\\
\hh\id\step\id\mobj{X_j}\step\id\step\id\\
\hh\id\step\frabox{s}\step\id\\
\hh\id\step\id\step\mobj{X_j}\id\step\id\\
\hh\id\step\id\step\hev\\
\object{X_{\ol{\jmath}}}\step\object{X_i}
\end{tangle}
\label{slocal}\end{equation}
on the direct summands, where $\begin{tangle}\hh\hcoev\end{tangle}$,
$\begin{tangle}\hh\hev\end{tangle}$ are any solution of the duality equation for
$X_j, X_{\ol{\jmath}}$, is a vector space isomorphism of order four. 
\elemma 
\prf The above map
$\Hom_\2C(X_iX_j,X_jX_i)\rarr\Hom_\2C(X_{\ol{\jmath}}X_i,X_iX_{\ol{\jmath}})$ is an
isomorphism by duality. The same holds for $\6S$ which is just a direct sum of such
isomorphisms, since the map $(i,j)\mapsto(\ol{\jmath},i)$ is a permutation of 
$\Gamma\times \Gamma$. That $\6S$ has order four is an obvious consequence of sphericity 
of $\2C$. 
\qed 

\blemma \label{slemma}Let $(X,e_X), (Y,e_Y)$ be simple objects in $\2Z_1(\2C)$. Then
\begin{equation}\label{smat}
  z_{(Y,e_Y)}\6S(z_{(X,e_X)})= \frac{d(X)}{d(Y)\lambda^2}S(\ol{(X,e_X)},(Y,e_Y))\cdot
  z_{(Y,e_Y)}.  
\end{equation}
\elemma
\prf With (\ref{isotyp}), (\ref{slocal}), and (\ref{tube2}) we compute
\bean \lefteqn{ \left( z_{(Y,e_Y)}\6S(z_{(X,e_X)}) \right)[i,j,i]}  \\
 && = (d_j\lambda)^{-1}
\sum_{k,l}\sum_{\alpha=1}^{N_{kl}^j}\sum_{\beta=1}^{N_k^X}\sum_{\gamma=1}^{N_i^Y} 
d_kd_l \frac{d(X)}{\lambda d_k}\frac{d(Y)}{\lambda d_i} \quad
\begin{tangle}
\step[4]\object{X_j}\Step\object{X_i}\\
\step[2]\mobj{{t'}_{kl}^{j\alpha}}\step\cd\step\id\\
\step[3]\id\Step\id\step\O {{q'}_i^\gamma}\\
\step[2]\mobj{X_k}\step\id\step[1.5]\mobj{Y}\hstep\hxx\\
\step[3]\id\Step\O {q_i^\gamma}\step\id\\
\coev\step\id\Step\id\step\id\\
\id\Step\id\step\O {{p'}_k^\beta}\Step\id\step\id\\
\id\step[1.5]\mobj{X}\hstep\hxx\step\mobj{X_i}\step\id\step\id\\
\id\Step\O {p_k^\beta}\step\id\Step\id\step\id\\
\id\Step\id\mobj{X_k}\step\ev\step\id\mobj{X_l}\\
\id\Step\Cu\step[-.2]\mobj{t_{kl}^{j\alpha}}\\
\object{X_i}\step[4]\object{X_j}
\end{tangle}
\eean
\[ =(d_j\lambda)^{-1}\sum_{k,l}\sum_{\alpha=1}^{N_{kl}^j}\sum_{\beta=1}^{N_k^X}
\sum_{\gamma=1}^{N_i^Y} d_kd_l \frac{d(X)}{\lambda d_k}\frac{d(Y)}{\lambda d_i} 
\quad\quad\quad
\begin{tangle}
\step\object{X_j}\step[3]\object{X_i}\\
\step[-1]\mobj{{t'}_{kl}^{j\alpha}}\step\cd\Step\O {{q'}_i^\gamma}\\
\id\mobj{X_k}\Step\mobj{Y}\xx\mobj{e_Y(X_l)}\\
\O {{p'}_k^\beta}\Step\O {q_i^\gamma}\Step\id\\
\step[-2.5]\mmobj{e_X(X_i)^{-1}}\step[2.5]\x\mobj{X}\Step\id\\
\id\Step\O {p_k^\beta}\Step\id\mobj{X_l}\\
\id\Step\cu\mobj{t_{kl}^{j\alpha}}\\
\object{X_i}\step[3]\object{X_j}
\end{tangle}
\]
where we have used Lemma \ref{l3}.
Replacing 
\[
\begin{tangle}
\object{X_k}\Step\object{X_l}\\
\cu\mobj{t_{kl}^{j\alpha}}\\
\step\object{X_j}
\end{tangle}
\quad\quad\quad \mbox{by}\quad\quad \left(\frac{d_j}{d_l}\right)^{1/2}
\begin{tangle}
\object{X_k}\step[3]\object{X_l}\\
\id\Step\cd\mmmobj{{s'}_{\ol{k}j}^{l\alpha}}\\
\ev\step[-.5]\mobj{\ve(X_k)}\step[2.5]\id\\
\step[4]\object{X_j}
\end{tangle} 
\]
where $\{{s'}_{\ol{k}j}^{l\alpha}\}$ is a basis in $\Hom(\ol{X_k}X_j,X_l)$, 
and correspondingly for the dual basis, pulling ${s'}_{\ol{k}j}^{l\alpha}$
through the half braiding $e_Y(\cdot)$ and summing over $l, \alpha$ we obtain 
\[ = \frac{d(X)d(Y)}{d_i\lambda^3} \sum_k\sum_{\beta=1}^{N_k^X}\sum_{\gamma=1}^{N_i^Y}
\quad\quad 
\begin{tangle}
\step[4]\object{X_j}\Step\object{X_i}\\
\mobj{\ol{\ve}(X_k)}\step[4]\id\Step\O {{q'}_i^\gamma}\\
\coev\Step\xx\\
\step[-1]\mobj{X_k}\step\id\Step\mobj{Y}\xx\Step\id\\
\O {{p'}_k^\beta}\Step\O {q_i^\gamma}\Step\id\Step\id\\
\x\mobj{X}\Step\id\Step\id\\
\id\Step\O {p_k^\beta}\Step\id\mobj{\ol{X_k}}\Step\id\\
\id\Step\ev\step[-2]\mobj{\ve(X_k)}\step[4]\id\\
\object{X_i}\step[6]\object{X_j}
\end{tangle}
\]
By sphericity of $\2C$, naturality of the $e_Y(\cdot)$ and 
$\sum_{k,\beta} p_k^\beta\circ {p'}_k^\beta=\id_X$ this equals
\[ = \frac{d(X)d(Y)}{d_i\lambda^3} \sum_{\gamma=1}^{N_i^Y} \quad\quad
\begin{tangle}
\step[4]\object{X_j}\Step\object{X_i}\\
\mobj{\ol{\ve}(X)}\step[4]\id\Step\O {{q'}_i^\gamma}\\
\coev\Step\xx\\
\step[-.6]\mobj{X}\step[.6]\id\Step\mobj{Y}\xx\Step\id\\
\id\Step\O {q_i^\gamma}\Step\id\Step\id\\
\x\mobj{X}\Step\id\Step\id\mobj{\ol{X}}\\
\id\Step\ev\step[-2]\mobj{\ve(X)}\step[4]\id\\
\object{X_i}\step[6]\object{X_j}
\end{tangle}
\]
Using naturality of $e_X(\cdot)$ we can pull $q_i^\gamma$ through
$e_X(X_i)^{-1}$. Furthermore, since $(Y,e_Y)$ is simple we have 
\[
\begin{tangle}
\object{\ol{\ve}(X)}\Step\object{Y}\\
\hh\hcoev\step\id\\
\step[-.6]\mobj{X}\step[.6]\id\step\hxx\mobj{e_Y(\ol{X})}\\
\step[-2.5]\mobj{e_X(Y)^{-1}}\step[2.5]\hx\step\id\\
\hh\id\step\hev\\
\object{Y}\step\mmobj{\ve(X)}
\end{tangle} 
\quad\quad\quad = \quad\quad
\begin{tangle}
\object{\ol{\ve}(X)}\Step\object{Y}\\
\hh\hcoev\step\id\\
\id\step\hxx\mobj{e_Y(\ol{X})}\\
\step[-.6]\mobj{X}\step[.6]\id\step\hxx\mobj{e_{\ol{X}}(Y)}\\
\hh\hev\step\id\\
\object{\ve(X)}\Step\object{Y}
\end{tangle}
\quad\quad\quad =\frac{S(\ol{(X,e_X)},(Y,e_Y))}{d(Y)}\,\id_Y, \]
and (\ref{smat}) follows by comparison with (\ref{isotyp}).
\qed

\bprop $\6S$ maps the center of $\Xi_L$ into itself. The modular matrix $S$ is invertible.
\eprop
\prf Summing (\ref{smat}) over all classes of simple $(Y,e_Y)$ and using 
$\sum_{(X,e_X)} z_{(X,e_X)}=1_{\Xi_L}$ we obtain
\begin{equation}\6S \left(\frac{z_{(X,e_X)}}{d(X)}\right)= 
\sum_{(Y,e_Y)} \lambda^{-2}S(\ol{(X,e_X)},(Y,e_Y)) \, \frac{z_{(Y,e_Y)}}{d(Y)}, 
\label{sss}\end{equation}
whence the first claim. Therefore the isomorphism $\6S: \Xi_0\rarr\Xi_0$ restricts to
$Z(\Xi_L)$ and the matrix $\lambda^{-1}S(\ol{\cdot},\cdot)$ expresses the action of
$\6S\restr Z(\Xi_L)$ in terms of the basis $\{ d(X)^{-1}z_{(X,e_X)} \}$. Thus $S$ is
invertible. 
\qed 

\brem 1. Note that $\lambda^2=\dim\2C=\sqrt{\dim\2Z_1(\2C)}$. This is the correct
normalization since $(\dim\2M)^{-1/2}S$ is known to be of order four in every modular
category $\2M$ \cite{t,khr1}.

2. An alternative proof of the modularity of $\2Z_1(\2C)$ and of
$\dim\2Z_1(\2C)=(\dim\2C)^2$ could be given as follows. If $\2C$ satisfies the assumptions
of our Theorem \ref{main1}, there exists a finite dimensional quantum groupoid $H$ such
that $\2C$ is monoidally equivalent to the category $H-\mbox{Mod}$ of left modules over
$H$. In \cite{ntv}, the quantum double of finite dimensional quantum groupoids was
defined, and the category $D(H)-\mbox{Mod}$ was shown to be modular. Modularity of
$\2Z_1(\2C)$ follows, provided one proves the equivalence $\2Z_1(\2C)\simeq
D(H)-\mbox{Mod}$ of braided tensor categories. Proceeding in analogy to the Hopf algebra
case \cite{ka}, this should not present any serious difficulty. Yet, we think that a
direct categorical proof which avoids weak Hopf algebras is more satisfactory.

3. The tensor category $\2B=\END_\2E(\6B)$ defined in \cite{mue09} is known to be
equivalent to the category of $Q-Q$-bimodules, cf.\ \cite[Remark 3.18]{mue09}. Combining
this with the ideas of \cite{ostr}, it should be possible to give a considerably simpler
proof of the braided equivalence $\2Z_1(\2C)\simeq\2B$.
\erem

In order to give the promised analogue of the (rather trivial) observation
$Z_1(Z_0(S))=\{\id_S\}$ from the Introduction we need the following
\bdefin The center $\2Z_2(\2C)$ of a braided monoidal category $\2C$ is the full
subcategory defined by
\[ \obj\,\2Z_2(\2C)=\{ X\in\obj\,\2C\ | \ c(X,Y)=c(Y,X)^{-1} \ \forall Y\in\obj\,\2C \}. \]
\label{Z2}\edefin
Obviously the subcategory $\2Z_2(\2C)$ is symmetric, contains the monoidal unit and is
stable w.r.t. \ direct sums, retractions (in particular isomorphisms, thus replete) and
duals. 

\bcoro The category $\2Z_2(\2Z_1(\2C))$ is trivial, i.e. all objects are direct multiples
of the monoidal unit. 
\ecoro
\prf It is well known that a semisimple braided category $\2A$ containing a simple object
$X\not\cong\11$ in $\2Z_2(\2A)$ is not modular. ($X\in\2Z_2(\2A)$ implies
$S(X,Y)=d(X)d(Y)$ for all $Y$. This is colinear to $S(\11,Y)=d(Y)$.) 
\qed

\brem One can in fact prove \cite{bebl} that $\2C$ is modular iff $\dim\2C\ne 0$ and
the center $\2Z_2(\2C)$ consists only of the direct multiples of the unit or,
equivalently, iff all simple objects of $\2Z_2(\2C)$ are isomorphic to the unit object. We
will show this in Subsection \ref{moddoub} as a byproduct of a more general computation. 
\erem

\brem There is little doubt that a more conceptual understanding of the above proof (and
of the subsequent subsection) can be gained by looking at them in the light of
Lyubashenko's works \cite{ly1,ly2}. The latter also raise the question whether there is a
generalization to non-semisimple Noetherian categories. We hope to pursue this elsewhere. 
\erem

\brem It is natural to ask whether there are higher dimensional analogues to the above
result in $d=1$ and the trivial case $d=0$ mentioned in the Introduction. (See
\cite{baez2} for a review of the theory of n-categories.) Thus, considering the center
constructions in $d=2$ \cite{baez,crans}, can one show that
$\2Z^{(2)}_3(\2Z^{(2)}_2(\2C))$ is trivial? Here
$\2C, \2Z^{(2)}_2(\2C), \2Z^{(2)}_3(\2Z^{(2)}_2(\2C))$ are (semisimple spherical) braided,
sylleptic and symmetric 2-categories, respectively.
\erem


\subsection{Computation of the Gauss sums}\label{gss}
If $\2C$ is a braided spherical tensor category a theorem of Deligne, cf.\ \cite[Proposition
2.11]{y}, implies that $\2Z_1(\2C)$ is a ribbon category (or balanced). Namely, 
\[
\theta_X=\quad\quad
\begin{tangle}
\object{\ol{\ve}(\ol{X})}\Step\object{X}\\
\hh\hcoev\step\id\\
\step[-.8]\mobj{\ol{X}}\step[.8]\id\step\hxx\\
\hh\hev\step\id\\
\object{\ve(\ol{X})}\Step\object{X}
\end{tangle}
\quad\quad \mbox{for } X\in\2C
\]
defines a natural automorphism $\{ \theta_X, X\in\obj\,\2C \}$ of the identity functor
which satisfies 
\[ \theta_{XY}=\theta_X\otimes\theta_Y\mcirc c(Y,X)\mcirc c(X,Y) \quad
   \forall X,Y, \]
\[ \theta_{\ol{X}}=\ol{\theta_X}, \quad\forall X. \]
(A similar results hold for $*$-categories, cf.\ \cite{lro}.) For the simple objects we
have $\theta_X=\omega_X\id_X$ with $\omega_X\in\7F^*$. 

The quantum double $\2Z_1(\2C)$ is braided and by the arguments in Section \ref{sec-semisim}
we know that it has a spherical structure which is induced by the one on $\2C$. We will
show that the numbers $\omega_{(X,e_X)}$ can be computed in terms of the tube algebra and
will compute the Gauss sum 
\[ \Delta_\pm(\2Z_1(\2C))= \sum_{(X,e_X)} \omega^{\pm 1}_{(X,e_X)} d(X,e_X)^2, \]
which plays an important role in the construction of topological invariants.

Following \cite{iz2} we consider the element $t\in\Xi_L$ defined by
\begin{equation}t[i,j,k]= \frac{\lambda}{d_i} \delta_{ik}\delta_{ij} \id_{X_i^2}. 
\label{deft}\end{equation}
It will turn out that $t$ is in the center of $\Xi_L$.

\blemma For simple $(X,e_x)\in\2Z_1(\2C)$ we have
\[ tz_{(X,e_X)} = \omega^{-1}_{(X,e_X)}z_{(X,e_X)}. \]
\elemma
\prf From (\ref{tube2}), (\ref{isotyp}) and (\ref{deft}) we obtain
\[ (tz_{(X,e_X)})[i,j,k] = 
\delta_{ik} \frac{\lambda}{d_i}\frac{d(X)}{\lambda d_i}\sum_{m\in \Gamma}
   \sum_{\alpha=1}^{N_{mi}^j} \sum_{\beta=1}^{N^X_i} \frac{d_md_i}{d_j\lambda}
\begin{tangle}
\step[1.5]\hstep\object{X_j}\step[1.5]\object{X_i}\\
\hh\mobj{{t'}_{mi}^{j\alpha}}\step[1.5]\hcd\step\id\\
\hh\hstep\mobj{X_m}\step\id\step\id\obj{X_i}\step\id\\
\step[1.5]\id\step\O {{p'}^\beta_i}\step\id \\
\step\obj{X}\hstep\hxx\step\id \\
\step[1.5]\O {p^\beta_i}\step\id\step\id \\
\hh\step[1.5]\id\step\id\obj{X_m}\step\id\\
\hh\step[1.5]\id\step\hcu\obj{t_{mi}^{j\alpha}}\\
\step[1.5]\object{X_i}\step[1.5]\object{X_j}
\end{tangle}
\]
\[ = \delta_{ik} \frac{d(X)}{\lambda d_i}\sum_{m\in \Gamma} \sum_{\alpha=1}^{N_{mi}^j}
  \sum_{\beta=1}^{N^X_i} 
\begin{tangle}
\hstep\object{X_j}\step[1.5]\object{\ol{\ve}(\ol{X_i})}\step[2.5]\object{X_i}\\
\hh\hstep\id\step\hcoev\Step\id\\
\hh\step[-1]\mobj{s_{j\ol{\imath}}^{m\alpha}}\step[1.5]\hcu\step\id\Step\id\\
\step[0.5]\dh\step\O {{p'}^\beta_i}\Step\id \\
\step\obj{X}\hstep\hxx\Step\id \\
\step[1.5]\O {p^\beta_i}\step\hd\step[-0.5]\obj{m}\step[1.5]\id \\
\hh\step[1.5]\id\step\hcd\step\id\\
\hh\step[1.5]\id\step\id\step\hev\mobj{\ve(\ol{X_i})}\\
\step[1.5]\object{X_i}\step\object{X_j}
\end{tangle}
\  =\delta_{ik} \frac{d(X)}{\lambda d_i}\sum_{\beta=1}^{N^X_i} 
\begin{tangle}
\hstep\object{X_j}\step[1.5]\object{\ol{\ve}(\ol{X_i})}\step[1.5]\object{X_i}\\
\hh\hstep\id\step\hcoev\step\id\\
\step[0.5]\id\step\id\step\O {{p'}^\beta_i}\step\id \\
\step[0.5]\id\step\hxx\step\id\\
\obj{X}\hstep\hxx\step\hev\step[-.5]\mobj{\ve(\ol{X_i})}\\
\hstep\O {p^\beta_i}\step\id\\
\hh\hstep\id\step\id\\
\hstep\object{X_i}\step\object{X_j}
\end{tangle}
\]
Now the claim is a consequence of the following computation
\[ 
\begin{tangle}
\step[2]\object{X_i}\\
\hh\hcoev\step\id\\
\id\step\O {{p'}^\beta_i}\step\id \\
\hxx\step\id\\
\hh\id\step\hev\\
\object{X}
\end{tangle}
\ \ = \ \
\begin{tangle}
\step[2]\object{X_i}\\
\hh\hcoev\step\id\\
\O {p^\beta_i}\step\id\step\id \\
\hxx\step\id\\
\hh\id\step\hev\\
\object{X}
\end{tangle}
\ \ = \ \
\begin{tangle}
\step[2]\object{X_i}\\
\hh\hcoev\step\id\\
\hxx\step\id\\
\id\step\O {p^\beta_i}\step\id \\
\hh\id\step\hev\\
\object{X}
\end{tangle}
\ \ = \ \
\begin{tangle}
\step[2]\object{X_i}\\
\Step\O {{p'}^\beta_i}\\
\hh\hcoev\step\id\\
\hxx\step\id\\
\hh\id\step\hev\\
\object{X}
\end{tangle}
\ \ = \ \
\begin{tangle}
\Step\object{X_i}\\
\Step\O {{p'}^\beta_i}\\
\hh\hcoev\step\id\\
\id\step\hx\\
\hh\hev\step\id\\
\Step\object{X}
\end{tangle}
\ = \omega^{-1}_{(X,e_X)}\,{p'}^\beta_i,
\]
which is justified by the same arguments as in the proof of Lemma \ref{slemma}. Here we
used standard properties of the spherical structure in the first and third equalities and 
naturality of the half braiding $e_X(\cdot)$ w.r.t.\ the second argument in the second
equality. The rest follows since $(X,e_X)$ is simple. \qed

\bprop We have
\[ \Delta_\pm(\2Z_1(\2C)) = \dim\2C. \]
\eprop
\prf In view of $\sum z_{(X,e_X)}=\11_{\Xi_L}$ the lemma implies
\[ t=\sum_{(X,e_X)} \omega_{(X,e_X)}^{-1} z_{(X,e_X)}, \]
which proves that $t$ is central in $\Xi_L$.
To this equation we apply the linear form $\phi\in\Xi_{L}^*$
\[  \phi(x)= \lambda\sum_{i\in \Gamma} d_i \, tr_{X_i}(x[i,0,i]). \]
One one hand by (\ref{deft}) we clearly have $\phi(t)=\lambda^2$. On the other hand with
(\ref{isotyp}) we compute
\[ \phi(z_{(X,e_X)})= d(X) \sum_{i\in \Gamma} tr_{X_i}\left( \sum_{\alpha=1}^{N^X_i} 
  {p'}^\alpha_i \circ p^\alpha_i \right) =d(X)\sum_{i\in \Gamma} d_i N^X_i= d(X)^2. \] 
Putting everything together we obtain $\Delta_-(\2Z_1(\2C))=\lambda^2=\dim\2C$. The equality
for $\Delta_+$ follows from $\dim\2Z_1(\2C)=(\dim\2C)^2$ and the fact 
$\Delta_+(\2M)\Delta_-(\2M)=\dim\2M$, which holds for every modular category $\2M$
\cite{t,khr1}. 
\qed\\

This completes the proof of Theorem \ref{main1}.

\brem 1. A modular category satisfying $\Delta_+(\2C)=\Delta_-(\2C)$ gives rise to an
anomaly-free surgery TQFT, cf.\ \cite{t}. Thus for quantum doubles the construction of the 
associated TQFTs simplifies considerably. 

2. The representation category of a rational conformal quantum field theory is a 
braided $*$-category and the central charge $c\in\7R$ of the CQFT is related, cf.\ e.g.\ 
\cite{khr1}, to the Gauss sums $\Delta_-(\2C)$ by 
\[ \frac{\Delta_-(\2C)}{|\Delta_-(\2C)|}= \exp \left(\frac{2\pi ic_\2C}{8} \right). \]
Since the Gauss sum of a quantum double is given by $\Delta_-(\2Z_1(\2C))=\dim\2C$, thus
positive, we conclude that the `central charge' of a double satisfies
\[ c_{\2Z_1(\2C)}\equiv 0 \ (\mbox{mod}\ 8). \]
\erem


\section{The quantum double of a $*$-category}
Consider the quantum double $\2Z_1(\2C)$ of a $*$-category $\2C$. If 
$s\in\Hom_{\2Z_1(\2C)}((X,e_X),(Y,e_Y))\subset\Hom_\2C(X,Y)$ then clearly
$s^*\in\Hom_\2C(Y,X)$. It does, not, however, follow that
$s^*\in\Hom_{\2Z_1(\2C)}((Y,e_Y),(X,e_X))$. But there is a suitable full subcategory of
$\2Z_1(\2C)$ which is a $*$-category. 

\bdefin Let $\2C$ be a tensor $*$-category. Then the unitary quantum double 
$\2Z_1^*(\2C)$ is defined as $\2Z_1(\2C)$ except that the half braidings $e_X(Y)$ are 
required to be unitary, not just invertible. \edefin
\blemma Let $\2C$ be a tensor $*$-category. Then $\2Z_1^*(\2C)$ is a $*$-category. \elemma
\prf For $s\in\Hom_{\2Z_1(\2C)}((X,e_X),(Y,e_Y))\subset\Hom_\2C(X,Y)$ we have
\[ \id_Z\otimes s\mcirc e_X(Z)= e_Y(Z)\mcirc s\otimes\id_Z\quad\forall Z. \]
Starring this equation and using $e_X(Z)^*=e_X(Z)^{-1}$ we obtain
\[ \id_Z\otimes s^*\mcirc e_Y(Z)= e_X(Z)\mcirc s^*\otimes\id_Z\quad\forall Z, \]
thus $s^*\in\Hom_{\2Z_1(\2C)}((Y,e_Y),(X,e_X))$. \qed\\

In the applications of the quantum double to operator algebras, like to the asymptotic
subfactor \cite{iz2} or quantum field theory \cite{klm}, one is mainly interested in the
unitary quantum double. In order for the results of Theorem \ref{main1} to remain valid for
$\2Z_1^*(\2C)\subset \2Z_1(\2C)$ one must show $\2Z_1^*(\2C)$ that is equivalent to $\2Z_1(\2C)$
as a tensor category. Given an isomorphism $s: X\rarr Y$ in a $W^*$-category we can use
polar decomposition \cite{glr} to obtain a unitary morphism $\tilde{s}: X\rarr Y$. But
we cannot construct a unitary half-braiding in this way since it is not clear that the
unitaries $e_X(Z), Z\in\2C$ can be chosen such that naturality (\ref{hb-i}) and the braid
relation (\ref{hb-ii}) hold. Therefore a global approach is needed, which we develop using
our machinery from Section \ref{sec-morita}. 

Let $\2C$ be a $*$-category with conjugates, simple unit and finitely many simple objects.
All dimensions $d(X)$ are positive, and we choose the square roots of the latter and of
$\dim\2C$ to be positive. Reconsidering the constructions of Section \ref{sec-morita} we
now choose the bases $\{ t_{ij}^{k\alpha}\}$ in $\Hom_\2C(X_k,X_iX_j)$ to be orthonormal,
i.e.\ ${t'}_{ij}^{k\alpha}={t_{ij}^{k\alpha}}^*$, and similarly $v_i'=v_i^*$. Then
$v'=v^*$ and $w'=w^*$, such that the considerations of \cite[Subsection 5.3]{mue09} apply.  
We thus obtain a $*$-bicategory $\2E_*\subset\2E$ which is equivalent to $\2E$. The
considerations in Sections \ref{sec-morita} and \ref{sec-modular} of this paper remain
essentially unchanged except for replacing $\ve(X), \ol{\ve}(X)$ by standard solutions
$r_X, \ol{r}_X$ of the conjugate equations \cite{lro} everywhere.

\blemma \label{orthp}
Let $\2C$ be a $*$-category. Let $(X,e_X)\in\2Z_1(\2C)$ and 
$F(X,e_X)=(\oj(X\boxtimes\11)J,p_X)$. Then the idempotent 
$p_X\in\End_\2E(\oj(X\boxtimes\11)J)$ satisfies $p_X=p_X^{\#}$ iff $e_X(Z)$ is unitary for
all $Z$.
\elemma 
\prf We recall from Proposition \ref{funcf} that $p_X$ is given by
\[ p_X=\lambda^{-1} \sum_i v_i\otimes\id_{X\boxtimes\11} \mcirc
   e_X(X_i)\boxtimes\id_{X^\op_i} \mcirc\id_{X\boxtimes\11} \otimes v'_i.
\]
In view of the definition \cite[Subsection 5.3]{mue09} of the involution $\#$ on 
$\End_\2E(\oj (X\boxtimes\11)J)\equiv\Hom_\2A(\oj (X\boxtimes\11)JQ,Q\oj (X\boxtimes\11)J)$ 
we have 
\[ p_X^\#= \lambda^{-1} \sum_i \quad
\begin{tangle}
\object{Q}\step[1.5]\object{X\boxtimes\11}\\
\id\step\id\step\hcoev\mobj{r^*}\\
\hh\id\step\id\step\mfrabox{v_i}\id\\
\hh\id\step\id\step\id\obj{\hat{X}_i}\step\id\\
\hh\id\step\frabox{e^*_i}\step\id\\
\hh\id\step\id\obj{\hat{X}_i}\step\id\step\id\\
\hh\id\step\mfrabox{v_i^*}\id\step\id\\
\mobj{r}\hev\step\id\step\id\\
\step[1.5]\object{X\boxtimes\11}\step[1.5]\object{Q}
\end{tangle}
\]
where $e_i^*\equiv e_X(X_i)^*\boxtimes\id_{X_i^\op}$ and $r=w\circ v: \11\rarr Q^2$.
In view of (\ref{gamma}) it is clear that there are uniquely determined  
$\hat{r}_i: \11\rarr \hat{X}_{\ol{\imath}}\otimes\hat{X}_i,\ i\in \Gamma$ such that
\[ r = \sum_i v_{\ol{\imath}}\otimes v_i\mcirc\hat{r}_i. \]
Using $\id_Q\otimes r^*\mcirc r\otimes\id_Q=\id_Q$ one easily shows
\[ \id_{\hat{X}_i}\otimes\hat{r}^*_i\mcirc\hat{r}_{\ol{\imath}}\otimes\id_{\hat{X}_i} =
  \id_{\hat{X}_i}. \]
(This amounts to the identification $\ol{\hat{r}_i}=\hat{r}_{\ol{\imath}}$ which is
possible since all the self-conjugate $\hat{X}_i, i\in \Gamma$ are orthogonal.) Thus
\[ v_i^*\otimes\id_Q\mcirc r = \id_{\hat{X}_i}\otimes v_{\ol{\imath}}\mcirc
  \hat{r}_{\ol{\imath}}, \quad \quad
   \id_Q\otimes v_i^*\mcirc r = v_{\ol{\imath}}\otimes\id_{\hat{X}_i}\mcirc\hat{r}_i, \]
and we obtain
\[ p_X^\#= \lambda^{-1} \sum_i \quad
\begin{tangle}
\object{Q}\step[1.5]\object{X\boxtimes\11}\\
\id\step\id\step\hcoev\mobj{\hat{r}^*_{\ol{\imath}}}\\
\hh\id\step\id\step\id\step\mfrabox{v^*_{\ol{\imath}}}\\
\hh\id\step\id\step\id\obj{\hat{X}_i}\step\id\\
\hh\id\step\frabox{e^*_i}\step\id\\
\hh\id\step\id\obj{\hat{X}_i}\step\id\step\id\\
\hh\mfrabox{v_{\ol{\imath}}}\id\step\id\step\id\\
\mobj{\hat{r}_i}\hev\step\id\step\id\\
\step[1.5]\object{X\boxtimes\11}\step[1.5]\object{Q}
\end{tangle}
\]
This equals $p_X$ iff
\[ \id_{X_{\ol{\imath}}X}\otimes \ol{r}^*_i \mcirc 
  \id_{X_{\ol{\imath}}}\otimes e_X(X_i)^*\otimes\id_{X_{\ol{\imath}}} \mcirc
  r_i\otimes\id_{XX_{\ol{\imath}}}=e_X(X_{\ol{\imath}}) \quad\forall i\in \Gamma. \]
Considering Lemma \ref{l3}, this is the case iff $e_X(X_i)^*=e_X(X_i)^{-1}$ for all
$i\in \Gamma$. In view of Lemma \ref{halfbr} and the fact that the $x_i^\alpha$ occurring in
its proof are automatically isometries, this is equivalent to unitarity of $e_X(Z)$ for
all $Z$.
\qed

\btheor Let $\2C$ be a tensor $*$-category with simple unit, finitely many simple objects,
conjugates, direct sums and subobjects. Then $\2Z_1^*(\2C)$ is monoidally equivalent to
$\2Z_1(\2C)$, thus modular. \etheor  
\prf Let $(X,e_X)\in\2Z_1(\2C)$ and $(\oj(X\otimes\11)J,p_X)=F(X,e_X)$. Since
$\End_\2E(\oj(X\otimes\11)J)$ is a finite dimensional von Neumann algebra it contains an
orthogonal projection $q_X=q_X^2=q_X^*$ and an invertible element $s$ such that
$sp_Xs^{-1}=q_X$. It is clear that
$(\oj(X\boxtimes\11)J,q_X)\cong(\oj(X\boxtimes\11)J,p_X)$, and
by the lemma there is a unitary half braiding $\tilde{e}_X(\cdot)$ such that 
$F(X,\tilde{e}_X)=(\oj(X\boxtimes\11)J,q_X)$. Thus 
\[
\2Z_1^*(\2C)\ \stackrel{\otimes}{\simeq}\ \2B_*\ \stackrel{\otimes}{\simeq}\ \2B\
  \stackrel{\otimes}{\simeq} \2Z_1(\2C), \]
where 
$\2B_*\equiv\End_{\2E_*}(\6B)$. (The equivalence $\2B_*\cong\2B$ has already been
demonstrated in \cite[Proposition 5.6]{mue09}.)
\qed


\section{The quantum double of a braided category}\label{moddoub}
For the moment, let $\2C$ be any (strict) braided monoidal category. Given such a category
$\2C$ we denote by $\tilde{\2C}$ the braided monoidal category which coincides with $\2C$
as a monoidal category, but has the braiding 
\[ \tilde{c}(X,Y)= c(Y,X)^{-1}. \]

It is well known (e.g., \cite[Proposition XIII.4.3]{ka}) that for a braided monoidal category 
$\2C$ there is a strict braided monoidal functor $I: \2C\rightarrow \2Z_1(\2C)$ given by 
\[ I(X)= (X, e_X) \quad \mbox{with}\quad e_X(\cdot)=c(X,\cdot), \]
\[ I(f)= f \]
on the objects and morphisms, respectively. $I$ is full, faithful and injective on the
objects, thus an embedding of $\2C$ into $\2Z_1(\2C)$.

Now, also $\tilde{\2C}$ embeds into $\2Z_1(\2C)$ via the functor $\tilde{I}$ defined by
\[ \tilde{I}(X)= (X, \tilde{e}_X) \quad \mbox{with}\quad \tilde{e}_X(\cdot)=
   \tilde{c}(X,\cdot),
\]
\[ \tilde{I}(f)= f. \]

\blemma $I(\2C)$ and $\tilde{I}(\2C)$ are replete full subcategories of $\2Z_1(\2C)$.
\elemma
\prf By definition, $(Y,e_Y)\in\2Z_1(\2C)$ being isomorphic to $I(X)=(X,e_X)$ (where
$e_X(Z)=c(X,Z)$) means that there is an isomorphism $u: X\rightarrow Y$ in $\2C$ such 
that $e_Y(Z)= id_Z\otimes u\circ e_x(Z)\circ u^{-1}\otimes\id_Z$. With $e_X(Z)=c(X,Z)$
and naturality of the braiding $c$ this implies $e_Y(Z)=c(Y,Z)$ and thus
$(Y,e_Y)=I(Y)$. Thus $I(\2C)$ is replete. The proof for $\tilde{I}(\tilde{\2C})$
clearly is the same. That $I(\2C)$ and $\tilde{I}(\2C)$ are full subcategories of 
$\2Z_1(\2C)$ follows from naturality of the braiding in $\2C$, which implies that every
morphism $u: X\rightarrow Y$ in $\2C$ automatically satisfies the condition (\ref{hb-i})
in Definition \ref{hb} and thus is a morphism from $(X,c(X,\cdot))$ to $(Y,c(Y,\cdot))$. 
\qed

\bdefin Two subcategories $\2A, \2B$ of a braided tensor category are said to commute iff
$c(X,Y)\circ c(Y,X)=\id_{YX}$ for all $X\in\2A, Y\in\2B$. For a braided monoidal category
$\2C$ and a subcategory $\2A$ the relative commutant $\2C\cap\2A'$ is the full subcategory
defined by 
\[ \obj\,\2C\cap\2A'=\{ X\in\obj\,\2C\ | \ c(X,Y)\circ c(Y,X)=\id_{YX} 
   \ \forall Y\in\obj\,\2A \}. \]
\edefin

The properties of the braiding imply that $\2C\cap\2A'$ is monoidal and stable under
isomorphisms (thus replete), direct sums, retractions and two-sided duals. When there is
no danger of confusion about the ambient category $\2C$ we write also simply $\2A'$. Note
that $\2Z_2(\2C)=\2C\cap\2C'$, which justifies the terminology `center'. 

\bprop Let $\2C$ be braided monoidal. Then
\bean \2Z_1(\2C)\cap I(\2C)' &=& \tilde{I}(\tilde{\2C}), \\
     \2Z_1(\2C)\cap \tilde{I}(\tilde{\2C})' &=& I(\2C). \eean
\label{64}\eprop
\prf By definition, $\2Z_1(\2C)\cap I(\2C)'$ is the full subcategory of $\2Z_1(\2C)$ 
whose objects $(X, e_X)$ satisfy 
\[ c((X,e_X),(Y,e_Y))\circ c((Y,e_Y),(X,e_X))=id_{(Y,e_Y)(X,e_X)} \quad 
   \forall (Y,e_Y)\in \Gamma(\2C). \]
Using the definition of $I$ and the definition of the braiding in $\2Z_1(\2C)$ by 
$c((X,e_X),(Y,e_Y))=e_X(Y)$ we obtain
\[ \obj\,\2Z_1(\2C)\cap I(\2C)'= \{ (X,e_X)\in \2Z_1(\2C) \ | \ 
   e_X(Y)\circ c(Y,X)=\id_{YX} \}. \]
But this amounts to 
\[ \obj\,\2Z_1(\2C)\cap I(\2C)'= \{ (X,e_X) \ | \ X\in\2C,\ e_X(Y)=c(Y,X)^{-1} \}, \]
which is nothing but $\obj\tilde{I}(\tilde{\2C})$. The second equality is proven in 
the same way. \qed

\brem As an obvious consequence we see that the subcategories $I(\2C)$ and
$\tilde{I}(\tilde{\2C})$ of $\2Z_1(\2C)$ are equal to their second commutants:
$I(\2C)''=I(\2C)$. Note that this holds without any technical assumptions on $\2C$.
See Remark \ref{rem-dc} below for remarks on a general double commutant theorem.
\erem

The next observation provides another link between the centers $\2Z_1$ and $\2Z_2$ (besides
the triviality of $\2Z_2(\2Z_1(\2C))$ stated by Theorem \ref{main1}). It can be interpreted 
as saying that $I(\2C)\vee\tilde{I}(\tilde{\2C})$, the monoidal subcategory of $\2Z_1(\2C)$
generated by $I(\2C)\cong\2C$ and $\tilde{I}(\tilde{\2C})\cong \tilde{\2C}$, is an
amalgamated product over their intersection $\2Z_2(\2C)$.

\blemma Let $\2C$ be braided. Then in $\2Z_1(\2C)=\2Z_1(\2C)$ we have
\[ I(\2C)\cap\tilde{I}(\tilde{\2C})=I(\2Z_2(\2C)). \]
\label{63}\elemma
\prf Obviously,  $I(X)=\tilde{I}(Y)$ is equivalent to $X=Y$ and
$c(X,\cdot)=\tilde{c}(X,\cdot)$ and thus to $X\in\2Z_2(\2C)$. \qed\\

The following results are stated in somewhat greater generality than needed here since we
have other applications in mind, cf.\ Remark \ref{rem-dc}.

\blemma \label{ff}
Let $\2C$ be monoidal and semisimple with two-sided duals. Let $\2K, \2L$ be
full monoidal subcategories which are semisimple (i.e.\ closed under direct sums and
retractions) and have trivial intersection $\2K\cap\2L$ in the sense that every object
contained both in $\2K$ and $\2L$ is a multiple of the tensor unit. If $K_1, K_2\in\2K$
and $L_1, L_2\in\2L$ then   
\begin{equation}\Hom_\2C(K_1L_1,K_2L_2)\cong \Hom_\2K(K_1,K_2)\otimes_\7F \Hom_\2L(L_1,L_2). 
\label{isom5}\end{equation}
More precisely, the linear maps
\[ \otimes: \Hom_\2K(K_1,K_2)\otimes_\7F \Hom_\2L(L_1,L_2) \longrightarrow
   \Hom_\2C(K_1L_1,K_2L_2) \]
induced by $(k, l)\mapsto k\otimes l$ are isomorphisms for all $K_1, K_2, L_1, L_2$. 
If $K_1,K_2\in\2K,\ L_1,L_2\in\2L$ are simple then $K_1L_1, K_2L_2\in\2C$ are simple. They
are isomorphic iff $K_1\cong K_2, L_1\cong L_2$. 
\elemma
\prf By duality we have
\[ \Hom_\2C(K_1L_1,K_2L_2)\cong \Hom_\2C(\ol{K}_2K_1,L_2\ol{L}_1). \]
Now $\ol{K}_2K_1\in\2K$ and $L_2\ol{L}_1\in\2L$, and since $\2K, \2L$ are {\it monoidal}
subcategories and closed w.r.t.\ retractions, all subobjects of $\ol{K}_1K_1,\
L_2\ol{L}_1$ are in $\2K$ and $\2L$, respectively. Since the monoidal unit $\11$ is (up to
isomorphism) the only simple object common to $\2K$ and our categories are semisimple, all
morphisms  $f: \ol{K}_2K_1\rightarrow L_2\ol{L}_1$ thus factorize  through the monoidal
unit: 
\[ \Hom_\2C(\ol{K}_2K_1,L_2\ol{L}_1) \cong\Hom_\2C(\ol{K}_2K_1,\11)\otimes_\7F
   \Hom_\2C(\11,L_2\ol{L}_1). \]
Using duality again we obtain (\ref{isom5}).
Thus $\Hom_\2C(K_1L_1,K_2L_2)$ and $\Hom_\2K(K_1,K_2)\otimes_\7F \Hom_\2L(L_1,L_2)$ 
have the same dimension and the $\otimes$-product on 
$\Hom(K_1,K_2)\times \Hom(L_1,L_2)$ extends to an isomorphism. The remaining claims are
obvious consequences.
\qed

\bprop \label{65}
Let $\2C$ be braided monoidal and semisimple with two-sided duals. Let $\2K, \2L$
be semisimple full monoidal subcategories which commute and have trivial intersection. 
Then the full monoidal subcategory $\2K\vee\2L$ of $\2C$ generated by $\2K$ and $\2L$ (by
tensor products and direct sums) is equivalent as a braided monoidal category to 
$\2K\boxtimes\2L$. If $\2C$ is spherical then this is an equivalence of spherical
categories.
\eprop 
\prf 
Consider the functor 
$T: \2K\otimes_\7F\2L\rarr \2K\vee\2L$ defined by $X\boxtimes Y\mapsto X\otimes Y$. By the
above it is full and faithful. In order to prove that $T$ is strong monoidal we compute
\bean T(X\boxtimes Y)\otimes T(Z\boxtimes W) &=& X\otimes Y\otimes Z\otimes W, \\
   T((X\boxtimes Y)\otimes(Z\boxtimes W)) &=& T((X\otimes Z)\boxtimes(Y\otimes W))=
      X\otimes Z\otimes Y\otimes W. \eean
Now the family 
\[ F_2(X\boxtimes Y,Z\boxtimes W) = \id_X\otimes c(Y,Z)\otimes\id_W \]
of morphisms
$T(X\boxtimes Y)\otimes T(Z\boxtimes W)\rarr T((X\boxtimes Y)\otimes(Z\boxtimes W))$
clearly is natural and makes $T$ strong monoidal. The easy proof of the coherence
condition is left to the reader. In order to show that $T$ is a braided tensor functor we
must prove that the diagram
\[
\begin{diagram}
  T(X\boxtimes Y)\otimes T(Z\boxtimes W) & \rTo^{c_{\2C}} & T(Z\boxtimes W)\otimes T(X\boxtimes Y) \\
  \dTo^{F_2} && \dTo_{F_2} \\
  T(X\boxtimes Y\otimes Z\boxtimes W) & \rTo_{T(c_{\2K\boxtimes\2L})} & T(Z\boxtimes W\otimes X\boxtimes Y)
\end{diagram}
\]
commutes, where $c_{\2K\boxtimes_\7F\2L}=c_\2K\boxtimes c_\2L$ is the direct product
braiding. Using the definition of $T$ and $F_2$ and taking into account that $\2K$ and
$\2L$ commute this is an easy exercise. 
Now the functor $T$ extends uniquely (up to natural isomorphism) to 
$\2K\boxtimes\2L=\ol{\2K\otimes_\7F\2L}^\oplus$, remaining braided monoidal by naturality
of the braiding. This extension is essentially surjective, thus an equivalence of braided
spherical categories. That the equivalence respects spherical structures (if present) is
obvious. 
\qed 

\bcoro \label{cor-fact}
Let $\2C$ be braided monoidal and semisimple with two-sided duals (and spherical
structure). Let $\2K\subset\2C$ be a semisimple full monoidal subcategory which has
trivial center $\2Z_2(\2K)=\2K\cap\2K'$. Then we have the equivalence
\[ \2K\boxtimes (\2C\cap\2K') \ \stackrel{\otimes}{\simeq}_{br} \ \2K\vee
  (\2C\cap\2K')\subset\2C \]
of braided (spherical) categories.
\ecoro
\prf The subcategory $\2L=\2C\cap\2K'$ commutes with $\2K$. Furthermore, 
$\2K\cap\2L=\2K\cap\2K'=\2Z_2(\2K)$ is trivial by assumption. Thus the proposition
applies. \qed

\brem \label{rem-dc}
If we knew that $\2K\vee(\2C\cap\2K')=\2C$, we could conclude that $\2C$ is equivalent, as
a braided tensor category, to the direct product $\2K\boxtimes (\2C\cap\2K')$.
In \cite{mue11} we will prove that this is indeed the case if $\2C$ is a modular
$*$-category. Thus whenever a modular category $\2C$ contains a modular category $\2K$ as
a full tensor subcategory then $\2C\stackrel{\otimes}{\simeq}_{br}\2K\boxtimes\2L$, where 
$\2L=\2C\cap\2K'$ is also modular. As a consequence, every modular category is a (finite)
direct product of prime (or simple) ones, usually in a non-unique way. 
The proof relies on the following double commutant theorem: If $\2C$ is a modular
$*$-category and $\2K$ is a semisimple monoidal subcategory closed under duality then (i)
$\2K''=\2K$ and  (ii) $\dim\2K\cdot\dim\2K'=\dim\2C$. Here we are interested only in the
full inclusion $I(\2C)\subset\2Z_1(\2C)$ where $\2C$ is modular, which can be treated
without the full strength of the double commutant theorem. (Recall that we have proven
$\2Z_1(\2C)\cap(\2Z_1(\2C)\cap I(\2C)')'=I(\2C)$.)
\erem

\btheor Let $\2C$ be a modular semisimple spherical category. Then there is a canonical
equivalence
\[ \2Z_1(\2C)\stackrel{\otimes}{\simeq}_{br} \2C\boxtimes\tilde{\2C} \]
of braided tensor categories.
\label{main2}\etheor
\prf We apply the corollary to $\2Z_1(\2C)$ and the subcategory $I(\2C)$. The latter is
braided isomorphic to $\2C$, thus has trivial center. Therefore, as braided spherical
categories 
\[ \2C\boxtimes\tilde{\2C}\cong    I(\2C)\boxtimes\tilde{I}(\tilde{\2C})=
   I(\2C)\boxtimes \2Z_1(\2C)\cap I(\2C)' \simeq I(\2C)\vee I(\2C)'. \]
Thus we are done if we can show that the full subcategory
$I(\2C)\vee\tilde{I}(\tilde{\2C})$ exhausts $\2Z_1(\2C)$. If we assume $\2C$ to be a
$*$-category then also $\2Z_1^*(\2C)\simeq\2Z_1(\2C)$ is. By the above, we have
$\dim(I(\2C)\vee\tilde{I}(\tilde{\2C}))=\dim(\2C\boxtimes\tilde{\2C})=(\dim\2C)^2$, which
coincides with the dimension of $\2Z_1(\2C)$ by our main theorem. Since 
$\dim(I(\2C)\vee\tilde{I}(\tilde{\2C}))$ is a full semisimple subcategory of $\2Z_1(\2C)$
the categories must coincide. This argument does not work if $\2C$ is not a
$*$-category. We therefore give another proof which works in generality.

To this purpose we show that the minimal central idempotents of $\Xi_L$
corresponding to the simple objects $I(X_k)\tilde{I}(X_l)\in\2Z_1(\2C), \ k,l\in \Gamma$ sum up
to the unit of $\Xi_L$. By the definitions of $I, \tilde{I}$ we have
$I(X_k)\tilde{I}(X_l) = (X_kX_l, e_{X_kX_l}(\cdot))$ with
\[ e_{X_kX_l}(Z)= c(X_k,Z)\otimes\id_{X_l}\mcirc \id_{X_k}\otimes c(Z,X_l)^{-1}. \]
Thus according to Proposition \ref{isotyp} the sum over the corresponding minimal central
idempotents in $\Xi_L$ is given by
\[ \left( \sum_{k,l} z(I(X_k)\tilde{I}(X_l))\right)[i,j,n]= 
   \frac{\delta_{in}}{d_i}\sum_{k,l}\sum_{\alpha=1}^{N_{kl}^i} d_kd_l \quad  \
\begin{tangle}
\object{X_j}\step[1.5]\object{X_i}\\
\id\step\hcd\\
\hxx\step\id\\
\step[-1]\mobj{X_k}\step\id\step[0.3]\mobj{X_l}\step[0.7]\hx \\
\hcu\step\id\\
\hstep\object{X_i}\step[1.5]\object{X_j}
\end{tangle}
\]
Computations which are identical to those in Lemma \ref{condexp1} (except for turning an
over- into an under-crossing) show that the right hand side equals
\[
= \sum_k d_k \quad
\begin{tangle}
\object{X_j}\step[3]\object{X_i}\\
\id\step\hcoev\step\id\\
\hxx\step\id\step\id\\
\step[-1]\mobj{X_k}\step\id\step[0.1]\mobj{X_{\ol{k}}}\step[0.9]\hx\step\id \\
\hev\step\hx\\
\Step\object{X_i}\step\object{X_j}
\end{tangle}
\quad = \ \ \sum_k d_k \quad
\begin{tangle}
\object{X_j}\step\object{X_i}\\
\hxx\step\hcoev\\
\id\step\hxx\step\id\\
\id\mobj{X_k}\step\id\step[0.1]\mobj{X_{\ol{k}}}\step[0.9]\hx \\
\hh\id\step\hev\step\id\\
\object{X_i}\step[3]\object{X_j}
\end{tangle}
\]
But this is nothing else than 
\begin{equation} \label{slide}
   = \left(\sum_k d_k S(X_k,X_j)\right)c(X_j,X_i)^{-1} = \left(\sum_k d_k
  S(X_k,X_j)\right) c(X_i,X_j). 
\end{equation}
Since $\2C$ is assumed modular we have $\sum_k d_k S(X_k,X_j)=\delta_{j,0}\dim\2C$ and
thus 
\[ \left( \sum_{k,l} z(I(X_k)\tilde{I}(X_l))\right)[i,j,n]= 
    \delta_{in}\delta_{j,0}\dim\2C\,\id_{X_i}. \]
The reader is invited to convince himself that this is the unit of the tube algebra
$\Xi_L$ by plugging it into (\ref{tube2}). 
\qed

The method used in the proof allows to prove the following characterization of modular
categories, which appeared in \cite{bebl}.

\bcoro  \label{rem-mod} ({\it of proof}) 
Let $\2C$ be a $\7F$-linear semisimple spherical braided tensor category with finitely
many simple objects. Then $\2C$ is modular iff $\dim\2C\ne 0$ and $\2Z_2(\2C)$ is trivial.
\ecoro
\prf If $\2C$ is modular then $\dim\2C\ne 0$ \cite{t} and $\2Z_2(\2C)$ is trivial. 
If, conversely, $\2Z_2(\2C)$ is trivial then for every $j\ne 0$ there exists $i$ such that  
$c(X_i,X_j)\ne c(X_j,X_i)^{-1}$. But then (\ref{slide}) implies 
$\sum_k d_k S(X_k,X_j)=0\ \forall j\ne 0$, which is known to be equivalent to
invertibility of $S$.
\qed

\brem It is well known that a braided tensor category $\2C$ is monoidally isomorphic to
its reverse $\2C^\rev$ which coincides with $\2C$ as a category but has the tensor product 
reversed: $X\otimes_\rev Y:=Y\otimes X$. On the other hand the duality functor
$X\mapsto\ol{X}$ provides a monoidal equivalence
$\2C^\op\stackrel{\otimes}{\simeq}\2C^{rev}$. Putting this together we have
$\tilde{\2C}\stackrel{\otimes}{\simeq}\2C\stackrel{\otimes}{\simeq}\2C^\rev\stackrel{\otimes}{\simeq}\2C^\op$.
Thus for modular $\2C$ we actually have an equivalence
$\2Z_1(\2C)\stackrel{\otimes}{\simeq}\2C\boxtimes\2C^\op$ of tensor categories, not just
weak monoidal Morita equivalence. 
\erem


\section{Applications} \label{sec-appl}
\subsection{An adjoint for the forgetful functor $\2Z_1(\2C)\rarr\2C$} \label{ss-adj}
In Section \ref{sec-morita} we proved that the functor $F: \2Z_1(\2C)\rarr\2B$ is fully 
faithful and essentially surjective. By \cite[Theorem IV.4.1]{cwm} this implies that $F$
has a two-sided adjoint $G: \2B\rarr\2Z_1(\2C)$. Together with Lemma\ \ref{adj} this implies
\[ \Hom_\2C(X,Y)\cong\Hom_\2B(F(X,e_X),\oj (Y\boxtimes\11)J)\cong
   \Hom_{\2Z_1(\2C)}((X,e_X),G(\oj (Y\boxtimes\11)J)), \]
where $(X,e_X)\in\2Z_1(\2C), Y\in\2C$. With the forgetful functor 
$H: \2Z_1(\2C)\rarr\2C, (X,e_X)\mapsto X$ this becomes
\begin{equation}\Hom_\2C(H(X,e_X),Y)\cong \Hom_{\2Z_1(\2C)}((X,e_X),G(\oj (Y\boxtimes\11)J)).
\label{adjt}\end{equation}
We thus have  
\bprop The forgetful functor $H: \2Z_1(\2C)\rarr\2C, (X,e_X)\mapsto X$ has a two-sided
adjoint $K: \2C\rarr\2Z_1(\2C),\ X\mapsto G(\oj (X\boxtimes\11)J)$. On the objects one has
\begin{equation}K(Y)\cong \bigoplus_{(X,e_X)} \dim\Hom_\2C(X,Y) \, (X,e_X), 
\label{kf}\end{equation}
where the summation is over the isomorphism classes of simple objects in $\2Z_1(\2C)$. 
\eprop
\prf Eq.\ (\ref{adjt}) just says that $K$ is a right adjoint of $H$. That $K$ is also a
left adjoint of $H$ is proven in the same way. One must also show that the isomorphisms in
(\ref{adjt}) are natural w.r.t.\ $(X,e_X)$ and $Y$. We leave this to the reader. For
$Y=X_i$ and $(X,e_X)$ simple, (\ref{adjt}) implies that $K(X_i)$ contains $(X,e_X)$ with 
multiplicity $\dim\Hom_\2C(X,X_i)$. For general $Y$ we have
\[ K(Y)\cong \bigoplus_{i\in \Gamma} \dim\Hom(X_i,Y) \bigoplus_{(X,e_X)} 
   \dim\Hom_\2C(X,X_i) \, (X,e_X), \]
and (\ref{kf}) follows by semisimplicity of $\2C$. \qed

\brem By the general theory \cite{mue09} there is a dual Frobenius algebra 
$\hat{\5Q}=(\hat{Q},\ldots)$ in $\2B$, where $\hat{Q}=\oj J$. Under the equivalence 
$\2Z_1(\2C)\stackrel{\otimes}{\simeq}\2B$, $\hat{Q}$ corresponds to
\begin{equation} G(\oj J)=K(\11)\cong \bigoplus_{(X,e_X)} \dim\Hom_\2C(X,\11) \, (X,e_X).
\label{k_1}\end{equation}
Thus this object is part of the Frobenius algebra in $\2Z_1(\2C)$ which establishes the weak
monoidal Morita equivalence $\2Z_1(\2C)\approx\2C\boxtimes\2C^\op$. $K(\11)$ clearly
contains the unit $(\11,\id)$ of $\2Z_1(\2C)$ with multiplicity one. The reader might
find it amusing to identify explicitly the morphisms $w, w'$ in $\2Z_1(\2C)$ which come with
the canonical Frobenius algebra.
\erem


\subsection{Invariants of 3-manifolds}\label{TQFT}
There are two classes of invariants of 3-manifolds associated with a modular tensor
category $\2C$, cf.\ \cite{t}. On the one hand we have the surgery invariants $RT(M,\2C)$
of Reshetikhin and Turaev \cite{rt2} which are based on the fact that every connected
oriented closed 3-manifold can be obtained from $S^3$ by surgery along a framed link. It
turned out \cite{brug1} that modularity of the category $\2C$ is not really necessary,
since it suffices that $\2C$ be `modularizable'. Yet, the invariant of the manifold being
defined in terms of a link invariant, the existence of a (non-symmetric) braiding is
essential. On the other hand there are the state sum invariants based on a triangulation
of the manifold. Generalizing on \cite{tv}, an invariant $TV(M,\2C)$ associated with any
modular category $\2C$ was defined in \cite{t}. Later it was understood that in fact no
braiding is necessary for the construction of a triangulation invariant, cf.\ \cite{bw1,
gk}, provided $\2C$ has two-sided duals. (This had been anticipated in \cite{ocn4}, which
was never published.) We denote the corresponding invariant by $Tr(M,\2C)$. Gelfand and
Kazhdan formulated a conjecture \cite[Conjecture 1]{gk} pointing towards a link between
the two invariants being provided by the quantum double. Our results on the quantum double
of semisimple tensor categories allow us to prove this conjecture.

\bprop 
Let $\2C$ satisfy the assumptions of Theorem \ref{main1} and consider the state-sum TQFT
associated with $\2C$, as defined in \cite{gk}. Then the dimension of the vector space
$\2H_{S^1\times S^1}$ associated to the two dimensional torus equals the number
$\#\2Z_1(\2C)$ of isomorphism classes of simple objects of $\2Z_1(\2C)$. 
\eprop 
{\it Sketch of Proof.} By the considerations of Subsection \ref{ss-tube}, $\#\2Z_1(\2C)$
coincides with the dimension of the center of the tube algebra $\Xi_L$. But this center is 
isomorphic to $\2H_{S^1\times S^1}$, as discovered by Ocneanu \cite{ocn3} and explained in
more detail in \cite[Theorem 3.1]{ek1}. \qed 

\brem The above argument is only a sketch because the triangulation TQFT in 2+1 dimensions
considered in \cite{ocn3, ek1} is derived from a subfactor, see \cite{kosu} for a detailed
exposition. Here as in \cite[Section 7]{mue09} we use the fact that the latter is
equivalent to the invariant defined in \cite{bw1, gk}. This is more or less clear, but
certainly deserves being made precise, as we plan to do in \cite{mue15}. Note also that in
order for a spherical category to give rise to a triangulation TQFT -- as opposed to just
the invariant -- one must assume that it does not contain symplectic self-dual simple
objects. This is done in \cite{t} and \cite[p.\ 4018]{bw2}, but unfortunately ignored in
the bulk of the literature on the subject.
\erem

By Theorem \ref{main1}, $\2Z_1(\2C)$ is modular, thus gives rise to a surgery TQFT in 2+1
dimensions, cf.\ \cite{t}. For these TQFTs it is known (by construction) that the
dimension of the vector space $\2H_{S^1\times S^1}$ associated with the torus equals the 
number of isomorphism classes of simple objects in the category. Thus the above result
provides support for the conjecture that the triangulation and surgery TQFTs associated with
$\2C$ and $\2Z_1(\2C)$, respectively, are isomorphic. (This conjecture, while very natural,
seems to have appeared in print only in \cite[Question 5]{ke}.) In particular, the
corresponding invariants of closed oriented 3-manifolds should coincide:
\begin{equation} \label{conject}
  RT(M,\2Z_1(\2C))=Tr(M,\2C)\quad\forall M. 
\end{equation}
Presently, we have no proof for this, but we note that Sato and Wakui recently announced a
proof in the setting of unitary categories arising from a subfactor. If $\2C$ is modular,
the braided equivalence $\2Z_1(\2C)\stackrel{\otimes}{\simeq}_{br}\2C\boxtimes\tilde{\2C}$
proven in Section \ref{moddoub} implies
\[ RT(M,\2Z_1(\2C))=RT(M,\2C\boxtimes\tilde{\2C})=RT(M,\2C)\cdot RT(M,\tilde{\2C})
  =RT(M,\2C)\cdot RT(-M,\2C), \]
and (\ref{conject}) follows from \cite[Theorem VII.4.1.1]{t}, according to which
$TV(M,\2C)=RT(M,\2C)\cdot RT(-M,\2C)$. For non-modular $\2C$ we only have the following
weaker result:

\bprop 
Let $\2C$ be as in Theorem \ref{main1} and $M$ an oriented closed 3-manifold. Then
\[ RT(M,D(\2C))\cdot RT(-M,\2Z_1(\2C))=Tr(M,\2C)\cdot Tr(-M,\2C). \]
If $\2C$ is unitary then $|RT(M,\2Z_1(\2C))|=|Tr(M,\2C)|$.
\eprop
\prf We compute
\bean RT(M,\2Z_1(\2C))\cdot RT(-M,\2Z_1(\2C)) &=& TV(M,\2Z_1(\2C)) \\
  &=& Tr(M,\2Z_1(\2C)) \\
  &=& Tr(M,\2C\boxtimes\2C^\op) \\
  &=& Tr(M,\2C)\cdot Tr(M,\2C^\op) \\
  &=& Tr(M,\2C)\cdot Tr(-M,\2C). \eean
Here the first equality is due to Turaev's theorem, which applies since $\2Z_1(\2C)$ is
modular. The second is the equality \cite{bw2} of $TV$ and $Tr$ for spherical $\2C$. The
third equality follows from the weak monoidal Morita equivalence
$\2Z_1(\2C)\approx\2C\boxtimes\2C^\op$ together with the Morita invariance of the invariant
$Tr$, cf.\ \cite[Theorem 7.1]{mue09}. The last two equalities follow from general
properties of the invariant $Tr$ \cite{bw2}. 

If $\2C$ is unitary we have $RT(-M,\2Z_1(\2C))=\ol{RT(M,\2Z_1(\2C))}$ and
$Tr(-M,\2C)=\ol{Tr(M,\2C)}$, and we are done.
\qed


\subsection{Subfactor theory: The Longo-Rehren subfactor}\label{ss-subf}
As stated in the introduction, the present project originated in the author's observation 
that the quantum double of monoidal categories appears implicitly in Izumi's preprint
\cite{iz1}. Therefore it seems reasonable to briefly comment on the subfactor setting.

Let $M$ be a type III factor with separable predual. Then the tensor category $\End_f(M)$
of (normal unital $*$-) endomorphisms $\rho$ of $M$ such that $[M:\rho(M))]<\infty$ is a
$*$-category with duals, direct sums and subobjects. (Here one uses that every orthogonal
projection $p=p^2=p^*$ in $M$ is equivalent to $\11$, i.e.\ there is $V\in M$ such that
$V^*V=\11, VV^*=p$.) Let $\2C\subset\End_f(M)$ be a full monoidal subcategory with the
same completeness properties and finite dimension. Choosing representers 
$\{ \rho_i, i\in \Gamma\}$ for the classes of simple objects, defining
\[ A = M\ol{\otimes} M^\op \]
and picking a direct sum
\[  \gamma = \bigoplus_{i\in \Gamma} \rho_i \otimes \rho_i^\op \]
one shows \cite{lre} $\gamma$ to be part of a Frobenius algebra (`$Q$-system')
$(Q,v,v^*,w,w^*)$ in $\End_f(A)$. At this point one applies a beautiful and fundamental
result due to Longo \cite{lo}, which implies that there is a subfactor $B\subset A$ such
that $\gamma$ is a canonical endomorphism for the inclusion $B\subset A$. This means that
there is a normal morphism $\ol{\iota}: A\rarr B$ which is a dual (in the 2-category of
factors, morphisms and intertwiners) of the embedding morphism $\iota=\id: B\rarr A$, such
that $\gamma=\iota\circ\ol{\iota}$. The subfactor $B$ is simply given by 
\begin{equation}B= w^*\gamma(A)w. \label{ba}\end{equation}
(The verification that this really gives a subalgebra is easy.) We call the subfactor thus  
obtained from $M$ and $\2C$ the Longo-Rehren subfactor.

Among the objects of interest in subfactor theory are the monoidal subcategories 
$\Hom_{B\subset A}(A,A)\subset\End_f(A)$ and $\Hom_{B\subset A}(B,B)\subset\End_f(B)$
generated by $\gamma=\iota\circ\ol{\iota}$ and $\hat{\gamma}=\ol{\iota}\circ\iota$,
respectively, and the categories $\Hom_{B\subset A}(A,B)$, $\Hom_{B\subset A}(B,A)$ of
morphisms which are contained in $\ol{\iota}\circ(\iota\circ\ol{\iota})^n$ and
$\iota\circ(\ol{\iota}\circ\iota)^n$, respectively, for some $n\in\7Z_\ge$. The reader
should appreciate that in this way every subfactor with finite index provides us with a
$\7C$-linear $*$-2-category with two objects and with non-strict spherical structure, thus
in particular with a Morita context for the tensor categories $\Hom_{B\subset A}(A,A)$ and
$\Hom_{B\subset A}(B,B)$. (The dimension of the four categories of 1-morphisms is finite
iff the subfactor has finite depth.) Our painful construction in \cite{mue09} just models
the categorical structure implicit in subfactor theory, where thanks to the inbuilt
structure one just needs the simple formula (\ref{ba})!

Alas, the above construction does not necessarily yield (tensor) categories which are 
equivalent to the $\Hom_\2E(\6A,\6A)$, $\Hom_\2E(\6B,\6B)$, $\Hom_\2E(\6A,\6B)$,
$\Hom_\2E(\6B,\6A)$ of Section \ref{ss-mcata}. This becomes clear already by comparing our 
$\2A=\2C\boxtimes\2C^\op$ with 
\[ \Hom_{B\subset A}(A,A)= \{ \rho\in\End_f(A) \ | \ \rho\prec\gamma^n, n\in\7Z_\ge \}. \]
The latter obviously is (equivalent to) a full subcategory of $\2C\boxtimes\2C^\op$, but
they coincide only if every $\rho_i\otimes\rho_j, i,j\in \Gamma$ is contained in $\gamma^n$ for
some $n$. This condition can be shown to be equivalent to connectedness of a certain
graph, the fusion graph of $\2C$. 

With these preparations a short inspection of Izumi's work \cite{iz2}, where the
categorical double does not appear explicitly, shows that essentially he has proven the
following theorem 
\btheor Let $M$ be a type III factor with separable predual and let $\2C$ be a full
monoidal subcategory of $\End_f(M)$ which is closed under duals, direct sums,
subobjects and is finite dimensional, and let $B\subset A$ be the corresponding LR
subfactor. If the fusion graph of $\2C$ is connected then we have the following
equivalences of tensor categories: 
\bean \Hom_{B\subset A}(A,A) &\simeq& \2C\boxtimes\2C^\op, \\
   \Hom_{B\subset A}(B,B) &\simeq& \2Z_1(\2C). \eean
\etheor
\prf The fusion graph is connected iff the objects $X\boxtimes X^\op$ generate all of 
$\2C\boxtimes\2C^\op$. Thus the statement on $\Hom_{B\subset A}(A,A)$ is contained in 
\cite[Theorem 4.1]{iz2}. Under the connectedness assumption Izumi's `quantum double of
$\Delta$' coincides with the $B-B$-morphisms $\Hom_{B\subset A}(B,B)$. Then our second
claim follows from \cite[Theorem 4.6]{iz2}, where the quantum double appears in only slightly
disguised form. Instead of half braidings $Z\mapsto e_X(Z)$ satisfying the braid relation
and naturality Izumi uses maps $I\ni i\mapsto e_X(X_i)$ satisfying the braiding fusion
relation. These two pictures are equivalent by our Lemma \ref{l17}. 
\qed

\brem The above theorem is the precise formulation of Ocneanu's remarkable intuitive
insight \cite{ocn3} that his asymptotic subfactor \cite{ocn2} (which is strongly related 
\cite{mas} to the Longo-Rehren subfactor $B\subset A$) is ``the subfactor analogue of
Drinfel'd's quantum double''. In view of the fact that irreducible depth-two subfactors are 
precisely the subfactors arising from outer actions of a Hopf algebra, the most natural
way to make Ocneanu's claim precise would be the following: The asymptotic subfactor of
$M^H\subset M$ is isomorphic to $P^{D(H)}\subset P$ or its dual. Yet, this clearly cannot
be the case since the index of the asymptotic inclusion coincides with the global index of
the original subfactor, which for depth two coincides with the index. Thus for
$N=M^H\subset M$, $[A:B]=[M:N]$ and $B\subset A$ cannot arise from a $D(H)$-action. 
\erem

Using the results of Subsection \ref{ss-adj} we can remove the connectedness 
condition: 
\bcoro Let $M, \2C$ be as in the theorem, but with possibly disconnected fusion
graph. Then $\Hom_{B\subset A}(A,A)$ is equivalent to the monoidal subcategory of
$\2C\boxtimes\2C^\op$ generated by the $X\boxtimes X^\op$ where $X$ runs through the
simple objects of $\2C$. $\Hom_{B\subset A}(B,B)$ is braided equivalent to the sub-tensor
category of $\2Z_1(\2C)$ generated by those simple objects $(X,e_X)\in\2Z_1(\2C)$ for which
$X$ contains the tensor unit $\11$ of $\2C$. 
\ecoro
\prf The first statement is well known.
By definition, $\Hom_{B\subset A}(B,B)$ is generated by the dual canonical object 
$\hat{Q}$. For the LR-subfactor this is the $K(\11)$ given in (\ref{k_1}), from which the 
second claim follows. \qed

\vspace{1cm}
\noindent{\it Acknowledgments.} During the long time of preparation of this work I was
financially supported by the European Union, the NSF and the NWO and hosted by the 
universities ``Tor Vergata'' and ``La Sapienza'', Rome, the IRMA, Strasbourg, the School
of Mathematical Sciences, Tel Aviv, the MSRI, Berkeley, and the Korteweg-de Vries
Institute, Amsterdam, to all of which I wish to express my sincere gratitude. 

The results of this paper and of \cite{mue09} were presented at an early stage at the
conference {\it Category Theory 99} at Coimbra, July 1999, at the conference {\it
$C^*$-algebras and tensor categories} at Cortona, August 1999, and at the workshop {\it
Quantum groups and knot theory} at Strasbourg, September 1999.

On these and other occasions I received a lot of response and encouragement. The following
is a long but incomplete list of people whom I want to thank for their interest and/or
useful conversations:
J.\ Baez, J.\ Bernstein, A.\ Brugui\`{e}res, D.\ E.\ Evans, J.\ Fuchs, F.\ Goodman, 
M.\ Izumi, V. F. R. Jones, C.\ Kassel, Y.\ Kawahigashi, G. Kuperberg, R.\ Longo, 
G.\ Maltsiniotis, G.\ Masbaum, J.\ E.\ Roberts, N.\ Sato, V.\ Turaev, L. Tuset, P.\ Vogel,
A.\ Wassermann, H.\ Wenzl and S. Yamagami.


\appendix

\section{On Quantum Doubles of Finite Dimensional Hopf Algebras}
The core of this paper was the proof that the quantum doubles of certain tensor
categories are modular. That $\Rep\, D(H)$ is modular has been proven for $H=\7CG$
\cite{ac}, where $G$ is a finite group, and for semisimple $H$ over an algebraically
closed field $k$ of characteristic zero \cite{eg}. A proof which covers also weak Hopf
algebras (or finite quantum groupoids) is given in \cite{ntv}. Our aim in this appendix 
is to give a proof which uses the ideas of Lyubashenko \cite{ly1,ly2} and Majid
\cite{lyma} and therefore is more in the spirit of our proof in the categorical situation.
In the sequel $H$ will always be a finite dimensional Hopf algebra. Since the main
application will be to quantum doubles the following will be useful. 

\blemma \label{l-double}
Let $H$ be a finite dimensional Hopf algebra $H$ over the field $k$ and let $D(H)$ be the 
quantum double. The following are equivalent: 
\begin{itemize}
\item[(i)] $H$ is semisimple and cosemisimple.
\item[(ii)] The antipode of $H$ is involutive and $char\,k\nmid\dim\,H$.
\item[(iii)] $D(H)$ is semisimple.
\item[(iv)] $D(H)$ is cosemisimple.
\end{itemize}
\elemma
\brem If the characteristic of $k$ is zero then $H$ is semisimple iff it is cosemisimple,
and the second condition in (ii) is vacuous.
\erem
\prf For the equivalences (i)$\Leftrightarrow$(iii)$\Leftrightarrow$(iv) see \cite{rad}
and for (i)$\Leftrightarrow$(ii) see \cite{eg2}.
\qed

In order for the category $\Rep\, D(H)$ to be modular it must be semisimple, which by
the lemma reduces us to the case where $H$ satisfies (i) and (ii). 
\blemma Let $H$ satisfy the (equivalent) conditions of Lemma \ref{l-double}.
Then there are two-sided integrals $\Lambda\in H, \mu\in \hat{H}$ which are traces in the
sense that
\[ \langle \mu,ab\rangle = \langle\mu, ba\rangle, \quad\quad
   \langle \alpha\beta,\Lambda\rangle = \langle \beta\alpha,\Lambda\rangle
\]
for all $a,b\in H, \alpha\beta\in\hat{H}$. The category $\Rep\, D(H)$
is a spherical category.
\elemma
\prf Semisimple Hopf algebras are unimodular \cite{ls}, which by definition means that
there are two-sided integrals. By \cite[Proposition 8]{ls}, unimodular Hopf algebras
satisfy
\[ \langle\mu, ab\rangle=\langle\mu, bS^2(a)\rangle \quad \forall a,b\in H, \]
thus $\langle\mu,\cdot\rangle$ is tracial by involutivity of $S$. Sphericity of
$\Rep\, D(H)$ is now an obvious consequence of \cite{bw2} where it was shown under the
weaker assumption that $S^2$ is inner.
\qed
 
We briefly recall some results on quasitriangular Hopf algebras. As shown be Drinfel'd
\cite{drin1}, the antipode of a finite dimensional quasitriangular Hopf algebra $H$ is
inner, i.e.\ there is an invertible $u\in H$ such that $S^2(A)=uAu^{-1}$. One has the
explicit formulae
\bean u &=& m\circ(S\otimes\id)(R_{21}), \\
 u^{-1} &=& m\circ(\id\otimes S^2)(R_{21}),
\eean
(I.e., $u=\sum_i S(f_i)e_i$ if $R=\sum_i e_i\otimes f_i$). Furthermore, Drinfel'd proved 
\[ uS(u)=S(u)u\in Z(H),\quad\quad\ve(u)=1, \quad\quad \Delta(u)=(R_{21}R)^{-1}(u\otimes u). \]

Recall \cite{t} that a ribbon \ha\ is a \qtha\ $H$ together with $\theta\in Z(H)$
satisfying
\begin{equation} \theta^2=uS(u), \quad\quad S(\theta)=\theta,\quad\quad
  \ve(\theta)=1,\quad\quad \Delta(\theta)=(R_{21}R)^{-1}(\theta\otimes\theta). 
\label{r-elem}\end{equation} 

\bprop \label{ribbon}
Let $H$ be a quasitriangular semisimple and cosemisimple Hopf algebra. Then $H$ is a
ribbon \ha\ with $\theta=u$.
\eprop
\prf Since $S^2=\id$ it follows that Drinfel'd's $u$ is central. Now \cite[Proposition 4]{kr}
implies that $u$ is a ribbon element. 
\qed

\brem For a ribbon Hopf algebra $H$ with $\theta=u$ the quantum trace, which for a
representation $\pi$ on a vector space $V$ is defined by
\[ \mbox{Tr}^q_\pi(X):=\mbox{Tr}\circ\pi(u\theta^{-1}\,X), \]
coincides with the usual trace $\mbox{Tr}$ on $\mbox{End}\,V$. In particular, all quantum 
dimensions $d(\pi)$ coincide with the classical dimensions $\dim V_\pi$.Therefore,
\[ \dim\Rep\, H=\sum_i d(\pi_i)^2=\sum_i (\dim V_{\pi_i})^2=\dim_k H. \]

\erem

In order to conclude that $\Rep\, H$ is modular it remains to prove that the $S$-matrix of
the ribbon category $\Rep\, H$ is invertible. A related notion of non-degeneracy was
introduced in \cite{rst}, where a \qtha\ was called {\it factorizable} if the map 
\[ \hat{H}\rightarrow H:\quad z\mapsto\langle z\otimes\id,I\rangle, \quad I=R_{21}R \]
is injective, thus invertible. Furthermore, it was shown that every \qd\ 
$D(H)$ is factorizable. The notion of factorizability plays an important role in the
works \cite{lyma,ke} where an action of $SL(2,\7Z)$ on ribbon \ha s is defined 
and studied.

\bdefin[\cite{lyma}] 
For a \qtha\ the selfdual Fourier transforms $\2S_+, \2S_-$ are defined by the linear
endomorphisms of $H$: 
\bean \2S_+(b) &=& (id\otimes\mu)(R_{21}(\11\otimes b)R_{12}), \\
\2S_-(b) &=& (id\otimes\mu)(R^{-1}_{12}(\11\otimes b)R^{-1}_{21}), \eean
where $\mu$ is a left integral in $\hat{H}$. If $H$ is ribbon the map $\2T: H\rightarrow H$ is
defined by $\2T(b)=\theta b$. 
\edefin 

\btheor[\cite{lyma}] For a factorizable ribbon \ha\ the following modular relations hold:
\begin{equation} \label{sl2z}
  \2S_+\circ\2S_- = \id = \2S_-\circ\2S_+,\quad (\2S_+\circ\2T)^3=\lambda\2S_+^2,   \quad
   \2S_+^2=\underline{S},  
\end{equation}
where $\underline{S}(x)=R^{(2)}S(Ad\,R^{(1)} (x))$ with $Ad\,Y(x)=Y_{(1)}xS(Y_{(2)})$
is the braided antipode,
and $\lambda\ne 0$ is defined by $\2S_+(\theta)=\lambda\theta^{-1}$. \etheor

\blemma The following decompositions hold in every finite dimensional \ha.
\[ a\otimes\11 = \sum_i (\11\otimes x_i)\Delta(y_i)=\sum_i \Delta(y_i)(\11\otimes
   S^2(x_i)) \] 
where
\[ \sum_i y_i\otimes x_i = (\id\otimes S^{-1})\Delta(a). \]
\elemma
\prf Inserting $\sum_i y_i\otimes x_i =a_{(1)}\otimes S^{-1}(a_{(2)})$ into 
$\sum_i (\11\otimes x_i)\Delta(y_i)$ we obtain
\begin{equation} (\11\otimes S^{-1}(a_{(2)}))\Delta(a_{(1)})=(\11\otimes S^{-1}(a_{(3)}))
   (a_{(1)}\otimes a_{(2)})=a_{(1)}\otimes S^{-1}(a_{(3)})a_{(2)}=a\otimes\11, \end{equation}
and the other equality is verified similarly. \qed

\bprop \label{cent}
Let $H$ be a quasitriangular semisimple, cosemisimple \ha\ and $\mu\in\hat{H}$ a left
integral. Then the selfdual Fourier transforms $\2S_\pm$ map the center of $H$ into itself:
\eprop
\begin{equation} \2S_\pm(Z(H))\subset Z(H). \end{equation}
\prf By \cite[Proposition 8]{ls} unimodularity is equivalent to the identity
\begin{equation} \mu(ab)=\mu(b S^2(a)) \ \forall a,b\in H, \label{trac}\end{equation}
which will be used in the sequel. Let $a\in H, \ b\in Z(H)$. Then
\bea a \, \2S_+(b) &=& a \, (id\otimes\mu)(R_{21}(\11\otimes b)R_{12}) \nn\\
  &=& (id\otimes\mu)((a\otimes\11) R_{21}R_{12}(\11\otimes b)) \nn\\
  &=& \sum_i (id\otimes\mu)((\11\otimes x_i)\Delta(y_i)R_{21}R_{12}(\11\otimes b)) \nn\\
  &=& \sum_i (id\otimes\mu)((\11\otimes x_i)R_{21}R_{12}(\11\otimes b)\Delta(y_i)) \\
  &=& \sum_i (id\otimes\mu)(R_{21}R_{12}(\11\otimes b)\Delta(y_i)(\11\otimes S^2(x_i)))\nn\\
  &=& (id\otimes\mu)(R_{21}R_{12}(\11\otimes b)(a\otimes\11))\nn\\
  &=& (id\otimes\mu)(R_{21}R_{12}(\11\otimes b)) \, a = \2S_+(b) \, a, \nn\eea
thus $\2S_+(b)\in Z(H)$. 
We have used $b\in Z(H),\ [ \Delta(\cdot), R_{21}R_{12} ] =0$, (\ref{trac}) and the 
lemma. \qed

\brem 1. In restriction to the center, the braided antipode $\underline{S}$ appearing in
(\ref{sl2z}) equals the antipode $S$.\\
2. For a ribbon algebra $H$ it is trivial that $\2T$ maps $Z(H)$ into itself, since 
$\theta\in Z(H)$.
\erem

The modularity condition requires invertibility of the matrix
\begin{equation} S_{i,j}=(Tr^q_{\pi_i}\otimes Tr^q_{\pi_j})(R_{21}R_{12}),
\label{modmat}\end{equation} 
where $i,j$ range over the equivalence classes of irreducible representations of $H$.
Now, in the semisimple case the representations are in one-to-one correspondence with the 
minimal projections in $Z(H)$, which leads to the following result.

\btheor \label{mod}
Let $H$ be a factorizable quasitriangular semisimple and cosemisimple \ha. 
Then the category $\Rep\, H$ is modular. 
\etheor
\prf We have already proven that $\Rep\, H$ is a spherical ribbon category. Thus
by Proposition \ref{cent} the center of $H$ is stable under the Fourier transform $\2S_+$. By
factorizability $\2S_+$ is invertible, and the same holds for the restriction 
$\2S_+\restr Z(H)$. By the remark after Lemma \ref{ribbon} the quantum traces on
$H$-modules $V$, in terms of which the $S$-matrix is defined, coincide with the usual
traces on $\mbox{End} V$. In terms of the basis for $Z(H)$ given by the minimal
idempotents $P_i$ we obtain 
\[ \2S_+(P_j)=\sum_i \2S_{ij}\, P_i ,\]
where the matrix $\2S=(\2S_{ij})$ is invertible. But $\2S$ is nothing but the modular
matrix (\ref{modmat}) as we have 
\bean \2S_{ij} &=& d_i \, \mu(P_i\,\2S_+(P_j)) 
   = d_i \, (\mu\otimes\mu)(R_{21}R_{12}\,(P_i\otimes P_j)) \\
  &=& \frac{1}{d_j}(\mbox{Tr}^q_i\otimes\mbox{Tr}^q_j)(R_{21}R_{12}) \eean
We have used $\mu(P_i)=1/d_i$ and $\mu(x\,P_i)=1/d_i \ \mbox{Tr}_i(x)$ \cite{w}. 
\qed

Note that the proof is conceptually quite similar to our proof of modularity for general
semisimple spherical categories. In view of Lemma \ref{l-double} the following is now
immediate: 

\bcoro 
Let $H$ be semisimple and cosemisimple. Then $\Rep\, D(H)$ is modular.
\ecoro

We close the appendix with a remark which is meant to aid the reader in appreciating
the `self-duality' of a quantum double $D(H)$, in particular since in general it is not
self-dual in the sense of \ha s: $D(H)^*\not\simeq D(H)$. (For a finite abelian group
$G$ we in fact have $D(G)\simeq D(G)^*\simeq \7CG\otimes\7C(G)$.) For any finite
dimensional Hopf algebra one can use the integrals to define `Fourier transforms'
$H\rightarrow\hat{H}$. In \cite{n} Fourier transforms
$\2F_{\sigma,\sigma'}, \sigma, \sigma'=\pm$, defined as linear maps $H\rightarrow\hat{H}$
which intertwine certain actions of $H$ on $H$ and $\hat{H}$ by multiplication and
translation, respectively, were studied systematically and used to give a new proof of 
the invertibility of the antipode. There is a beautiful relation between these more 
conventional Fourier transforms, relating $H$ and $\hat{H}$, and the selfdual 
Fourier transforms \cite{lyma}, which map $D(H)$ onto itself. For simplicity we restrict 
ourselves to FDKA's, where things are easier since the Haar measures are two-sided 
invariant traces and since there are unique Fourier transforms 
$\2F: H\rightarrow \hat{H}$ and $\hat{\2F}: \hat{H}\rightarrow H$:
\bea \langle \2F(x),y\rangle &=& \langle \mu, xS(y)\rangle, \quad\forall x,y\in H, \\
  \langle \alpha, \hat{\2F}(\beta)\rangle &=& \langle \hat{S}(\alpha)\beta, 
  \Lambda\rangle,  \quad\forall \alpha,\beta\in \hat{H}, \eea
where $\Lambda, \mu$ are the integrals in $H, \hat{H}$, respectively. If
$\iota: H\rarr D(H), \hat{\iota}: \hat{H}\rarr D(H)$ are the canonical embedding maps then
the following diagram commutes:
\[\begin{diagram}
\hat{H}\otimes H & \rTo^{\hat{\iota}\otimes\iota} & D(H)\otimes D(H) & \rTo^{m} & D(H)\\
  \dTo^{\hat{\2F}\otimes\2F} &&&& \dTo_{\2S_-} \\
  H\otimes \hat{H} & \rTo_{\iota\otimes\hat{\iota}} & D(H)\otimes D(H) & \rTo_{m} & D(H)
\end{diagram}\]
This nice observation is due to Kerler \cite[Proposition 9]{ke} (with different conventions),
who, however, did not emphasize the interpretation of $\2F, \hat{\2F}$ as conventional
Fourier transforms.

\end{document}